\documentclass[draft]{amsart}
\usepackage{times}
\usepackage{amsfonts}
\usepackage{verbatim}
\usepackage{amsmath}
\usepackage{amscd}
\usepackage{amssymb}
\usepackage{latexsym}
\usepackage{amsthm}

\input{xy}
\xyoption{all}

\theoremstyle{plain}

\newtheorem{conjectuur}{Conjecture}

\swapnumbers
\newtheorem{theorem}[subsection]{Theorem}
\newtheorem*{theorem-unnumbered}{Theorem}
\newtheorem{corollary}[subsection]{Corollary}
\newtheorem{lemma}[subsection]{Lemma}
\newtheorem{proposition}[subsection]{Proposition}

\theoremstyle{definition}
\newtheorem{definition}[subsection]{Definition}
\newtheorem{example}[subsection]{Example}

\theoremstyle{remark}
\newtheorem{notation}[subsection]{Notation}

\newtheorem{remark}[subsection]{Remark}

\newtheorem*{question-unnumbered}{Question}

\swapnumbers

\newcommand\NYCCT{\address{Department of Mathematics\\
NYC College of Technology (CUNY)\\
NY, NY 11201 (USA)}
}

\newcommand{\emptyprop}{q}

\newcommand \after{\circ}
\newcommand \algpow[2]{#1[[#2]]^{\text{alg}}}

\newcommand \complet[1]{\widehat {#1}}

\newcommand \id{\mathfrak a}

\newcommand \iso{\cong}

\newcommand \map[1]{{\newcommand{\tmpprop}{#1q} 
\if\tmpprop\emptyprop \to\else 
\xrightarrow{{\phantom{i}{#1}\phantom{i}}}\fi}}
\newcommand \maxim{\mathfrak m}
\newcommand \nat{\mathbb N}
\newcommand \norm[1]{\left|#1\right|}

\newcommand \pol[2]{#1[#2]}
\newcommand \pow[2]{#1[[#2]]}

\newcommand \pr{\mathfrak p}

\newcommand \range [2]{#1,\dots,#2}

\newcommand \rij[2]{(#1_1,\dots,#1_{#2})}

\newcommand \tensor{\otimes}

\newcommand \op\operatorname

\newcommand \acf{algebraically closed field}

\newcommand \ch{characteristic}
\newcommand \homo{homomorphism}
\newcommand \CM{Coh\-en-Mac\-au\-lay}
\renewcommand\iff{if and only if}
\newcommand \DVR{discrete valuation ring}

\newcommand \infal[1]{\op{Inf}(#1)}

\newcommand \zet{\mathbb Z}

\newcommand \hens[1]{#1^\sim}
\newcommand \ultraset{ultraset}
\newcommand \pc{$\mathfrak B$-closure}
\newcommand \remove[1]{\ \complet{#1}\ }

\newcommand \tc[1]{\op{cl}(#1)}
\newcommand \ac[2]{{\mathbb{#1}}_{#2}^{\operatorname{alg}}}
\newcommand \BS{Brian\c{c}on-Skoda}

\DeclareMathOperator*{\UP}{ulim}
\newcommand \up[1]{\UP_{#1}}

\newcommand \ul[1]{\seq{#1}\infty}
\newcommand \seq[2]{#1\mathstrut_{#2}}
\newcommand \seqaff[2]{#1\mathstrut^\text{aff}_{#2}}
\newcommand \sr{approximation}
\newcommand \Sr{Approximation}

\newcommand  \BCM[1]{\mathfrak B(#1)}
\newcommand  \hull[1]{\mathfrak D(#1)}

\newcommand  \los{\L os' Theorem}

\newcommand \frob[1]{\mathbf{F}_{#1}}
\newcommand \ulfrob{\frob\infty}

\renewcommand\frak{\mathfrak}
\newcommand\filtered{nested}
\newcommand\Filtered{Nested}
\newcommand\Frac{\operatorname{Frac}}

\numberwithin{equation}{subsection}

\newcommand \fff{\text{\it\underline{rin}g}}

\newcommand \loccat{\operatorname{\bf Loc}}
\newcommand \lefcat{\operatorname{\bf Lef}}
\newcommand \ancat{\operatorname{\bf An}_K}
\newcommand \cohcat{\operatorname{\bf Coh}_K}
\newcommand \extcohcat{\operatorname{\bf Coh}^*_K}

\title {Lefschetz Extensions,
Tight 
Closure, and big Cohen-Macaulay Algebras}
\author{Matthias Aschenbrenner}
\thanks{Partially supported  by a  grant from the National Science 
Foundation and by the Mathematical Sciences Research
Institute, Berkeley, CA}
\address {
Department of Mathematics, Statistics, and Computer Science \\
University of Illinois at Chicago \\
851 S. Morgan St. (M/C 249) \\
Chicago, IL 60607}
\email {maschenb@math.uic.edu}

\author{Hans Schoutens}
\thanks{Partially supported by a  grant from the National Science 
Foundation and by visiting positions at Universit\'e Paris VII and at 
the Ecole Normale Sup\'erieure.}
\NYCCT
\email{hschoutens@citytech.cuny.edu}
\date{November 15, 2004}
\begin{document}

\begin{abstract}
We associate to every equicharacteristic zero Noetherian local ring $R$
a faithfully flat ring extension which is an 
ultraproduct of rings of various prime \ch{s}, in a weakly functorial way. 
Since such 
ultraproducts carry naturally a non-standard  Frobenius, we can 
define a new tight closure operation on $R$ by mimicking the positive 
\ch\  functional definition of tight closure. This approach avoids 
the use of generalized N\'eron Desingularization and only relies 
on Rotthaus' result on Artin Approximation in \ch\ zero. If $R$ is moreover 
equidimensional and universally catenary, then we can also associate to it in a 
canonical, weakly functorial way a balanced big \CM\ algebra.
\end{abstract}

\maketitle

\setcounter{tocdepth}{1}
\tableofcontents

\pagebreak

\section*{Introduction}

In this paper, we investigate when a ring of \ch\ zero  can be
embedded  in an ultraproduct of rings of positive \ch. Recall that an
ultraproduct of a family of rings is a sort of `average' of its members; see
\S\ref{s:UP} for more details. To facilitate the discussion, let us
call a ring of \ch\ zero a \emph{Lefschetz ring} if it is realized as
an ultraproduct of rings of prime \ch. The designation alludes to  an
old heuristic principle in algebraic geometry regarding transfer
between positive and zero \ch, which Weil \cite{Weil} attributes to
Lefschetz. A \emph{Lefschetz field} 
is a Lefschetz ring which happens to be a field.
To model-theorists it is well-known that the field $\mathbb C$ of complex
numbers is Lefschetz. Moreover,
any field  of characteristic zero embeds
into a Lefschetz field. It follows that
\emph{any} domain of \ch\ zero embeds into a Lefschetz ring, but in
doing so, we loose the entire ideal theory of the domain. It is
therefore natural to impose that the embedding preserves enough of
the ideal structure, leading to:

\begin{question-unnumbered}
Given a Noetherian ring  $R$ of \ch\ zero, can we find a \emph{faithfully
flat} ring extension of $R$ which is Lefschetz?
\end{question-unnumbered}

Suppose that $R$ is a ring of characteristic zero
which admits a faithfully flat Lefschetz
extension  $D$. Hence $D$ is an ultraproduct of a family $(D_w)$ of
rings $\seq Dw$ of prime characteristic; infinitely many
different prime characteristics must occur. Each $\seq Dw$ can be
viewed as a kind of `reduction modulo $p$',  or \emph{\sr}, of $R$.
Faithful flatness guarantees that the $D_w$ retain enough properties
of the original ring. (See \S\ref{s:app} below.)
For an easy example consider the following criterion for ideal
membership in $R$: \textsl{given $f_0,\dots,f_s\in R$ and given
$\seq{f_i}w\in\seq Dw$ whose ultraproduct is equal to the image of
$f_i$ in $D$, we have $f_0\in\rij fsR$ \iff\
$\seq{f_0}w\in(\seq{f_1}w,\dots,\seq{f_s}w)\seq Dw$ for almost all
$w$}.

The main motivation for posing the above question stems from the
following observations. Any ring   of prime \ch\ $p$  admits an
endomorphism which is at the same time algebraic and canonical, to
wit, the Frobenius $\frob p\colon x\mapsto x^p$. This has an immense
impact on the homological algebra of a prime \ch\ ring, as is
witnessed by a myriad of papers exploiting this fact. To mention just
a few: Peskine-Szpiro \cite{PS72} on homological conjectures,
Hochster-Roberts  \cite{HR} on the \CM\ property of rings of
invariants,   Hochster \cite{Ho78} on big \CM\ algebras and
Mehta-Ramanathan \cite{MR} on Frobenius splitting of Schubert
varieties. This approach has  found  its culmination in the
\emph{tight closure theory} of  Hochster-Huneke \cite{HH,HHTC,HuTC}.
(For a more extensive history of the subject, see \cite[Chapter
0]{HuTC}; the same book is also an excellent introduction to tight
closure theory.)

Hochster and Huneke also developed tight closure in \ch\ zero (see
\cite{HHZero} or \cite[Appendix 1]{HuTC}), but without any appeal to
an endomorphism and relying on deep theorems about Artin
Approximation and N\'eron desingularization.
Any Lefschetz ring  $D$, however, is endowed with a \emph{non-standard
Frobenius} $\ulfrob$, obtained by taking the ultraproduct of the
Frobenii on the $\seq Dw$. The endomorphism $\ulfrob$ acts on the
subring $R$ of $D$, and although it will in general not leave $R$
invariant, its presence makes it possible to generalize the \ch\ $p$
functional definition of tight closure to any Noetherian ring $R$
admitting a faithfully flat Lefschetz extension.  This was carried
out in \cite{SchNSTC} for the case where $R$ is an algebra of finite
type over $\mathbb C$. Here we had a canonical choice for a
faithfully flat Lefschetz extension, called the \emph{non-standard
hull} of $R$. The resulting closure operation was termed
\emph{non-standard tight closure}. Variants and further results can
be found in \cite{SchRatSing,SchMixBCM,SchBCM,SchLogTerm,SchMixBCMCR}.

Let us briefly recall the construction of the non-standard hull of a
finitely generated algebra $A$ over a Lefschetz field $K$,
%$\mathbb C$, or more generally,
%over any uncountable \acf\ $K$ of \ch\ zero, 
and at the same time indicate
the problem in the non-affine case.   
For ease of exposition assume that
$K$ is an ultraproduct of fields $\seq Kp$ of \ch\ $p$,
with $p$ ranging over the set of 
prime numbers. (See also Proposition~\ref{P:acf} below.)
If $A$ is of the form $\pol
KX/I$, where $I$ is an ideal of $\pol KX=K[X_1,\dots,X_n]$, 
and we have already constructed   a faithfully flat Lefschetz
extension $D$ of $\pol KX$, then $D/ID$ is  a
faithfully flat Lefschetz extension of $A$. So we may assume $A=\pol
KX$. %(The argument is slightly more complicated if $A$ is a
%locali\-zation of a finitely generated $K$-algebra.)  
There is an obvious candidate for a Lefschetz ring, namely
the ultraproduct $\ul {\pol KX}$ of the $\pol{\seq Kp}X$. In a
natural way $\ul {\pol KX}$ is a $K$-algebra. Taking the ultraproduct
of the constant sequence $X_i$ in $K_p[X]$ yields an element in $\ul
{\pol KX}$, which we continue to write as $X_i$. By \los\ (see
Theorem~\ref{T:los} below), the elements $X_1,\dots,X_n\in\ul {\pol
KX}$ are algebraically independent over $K$ and hence can be viewed
as indeterminates over $K$. This yields a canonical embedding of
$\pol KX$ into $\ul {\pol KX}$. Van den Dries observed that this
embedding is faithfully flat \cite{vdD79,SvdD}, thus giving a
positive answer to the question above for finitely generated
$K$-algebras.

In \cite[\S3.3]{SchSymPow}, the Artin-Rotthaus Theorem \cite{AR} was 
used to extend  results from the finitely generated case to the 
complete case. This ad hoc application will be replaced   in this 
paper by constructing a faithfully flat Lefschetz extension  for 
{\em every}\/ Noetherian local ring  of equal \ch\ zero. However, for the 
proof, a stronger form of Artin Approximation is needed, to wit 
\cite{Rot}.  By the Cohen  Structure Theorem,  any equi\ch\ zero 
Noetherian local ring has a faithfully flat extension which is a 
homomorphic image of a power series ring $\pow KX$ (where  $K$ is as 
before), so the problem is essentially reduced to $\pow  KX$. There 
is again a natural candidate   for a faithfully flat  Lefschetz 
extension, namely the ultraproduct $\ul {\pow KX}$ of the  $\pow{\seq 
Kp}X$. Since $\ul  {\pol KX}$ is a subring of $\ul  {\pow  KX}$, so 
is $\pol KX$. Moreover, one easily verifies that   $\ul  {\pow KX}$ 
with the  $X$-adic topology is complete, that is to say,  each Cauchy 
sequence in $\ul  {\pow KX}$ has a limit in $\ul  {\pow  KX}$. The 
obstruction in extending the above argument from $\pol KX$ to $\pow 
KX$  is that the  $X$-adic topology on $\ul  {\pow KX}$ is not 
Hausdorff, and hence these limits are not unique. Therefore,  to send 
$f\in\pow KX$ to an element in $\ul  {\pow KX}$, we must pick a limit 
in $\ul  {\pow KX}$ of the Cauchy sequence $(f_n)$, where $f_n\in\pol 
KX$ is the truncation of $f$ at degree $n$. It is not at all obvious 
how to do this systematically  in order to get a ring  \homo{} 
$h\colon\pow KX \to \ul  {\pow KX}$. (It is not hard to prove that 
such an $h$, once defined, must be faithfully flat.) An example 
exhibits some of the subtleties encountered: Let us say that a
power series $f\in L[[X]]$, where $L$ is a field, \emph{does not involve
the variable $X_i$} if $f\in L[[X_1,\dots,X_{i-1},X_{i+1},\dots,X_n]]$. 
Similarly we say that
an element of $\ul  {\pow KX}$ \emph{does not involve the variable $X_i$} 
if it is the ultraproduct of power series in $\pow{\seq Kp}X$ not 
involving $X_i$. Using an example of from \cite{RobSol} we explain in 
\S\ref{s:uc} why there can be no \homo\ $h\colon\pow KX\to \ul  {\pow 
KX}$ with the property that for arbitrary $i$, if $f\in\pow KX$ does 
not involve the variable $X_i$, then neither does $h(f)$. (Of course 
there is always a limit of the  $f_n$ in $\ul{\pow KX}$ which has this property.)  To 
circumvent these kinds of problems, we  use  Artin Approximation 
to derive the following positive answer to the question posed at the 
beginning:

\begin{theorem-unnumbered}
For each equi\ch\ zero Noetherian local ring $R$, we can construct a
local Lefschetz ring $\hull R$ and a   faithfully flat embedding
$\eta_R\colon R\to \hull R$.
\end{theorem-unnumbered}

In fact, the result also holds for semi-local rings. 
More importantly,
$\mathfrak D$ can be made \emph{functorial}
 in a certain way, which is crucial for applications.
% on the category of  equi\ch\ zero Noetherian local 
%rings of cardinality at most $\kappa$, for any  choice of uncountable 
%cardinal $\kappa$. 
See Theorem~\ref{T:maincan} for the precise statement.

\subsection*{Organization of the paper}
%The paper is organized as follows. 
Sections~\ref{s:UP}--\ref{s:LP} of Part~1 contain a proof of  
the theorem above. The proof will be 
self-contained except for the use of Rotthaus' result \cite{Rot}. 
We also discuss further connections with Artin Approximation 
and cylindrical approximation. In Section~\ref{s:app} we investigate which 
algebraic properties are carried
over from $R$ to  the rings whose ultraproduct is $\hull R$.
The reader who is mostly interested in the applications of the
theorem might skip this section at first reading and immediately proceed
to Part~2 (referring back to Section~\ref{s:app} whenever necessary). 

We apply our main theorem in two ways.
First, in Section~\ref{s:tc} we define (\emph{non-standard})
tight closure in any  equi\ch\ zero Noetherian local ring and prove the
basic facts (such as its triviality on regular rings, Colon Capturing
and \BS). In contrast with the Hochster-Huneke version from
\cite{HHZero} or \cite[Appendix 1]{HuTC} we do not have to invoke
generalized N\'eron desingularization. In order for this paper not to
become too long,   issues such as the existence of test elements,
persistence of tight closure, detailed comparison with other tight closure
operations, F-rationality and F-regularity will be postponed to a future
publication. 

Our second application is a direct construction of a
balanced big \CM\ algebra for each equi\ch\ zero Noetherian local
ring, simpler than the one given in \cite{HHbigCM2}. 
This construction is weakly functorial on the subcategory of
equidimensional and universally catenary  rings of bounded
cardinality. 
Using non-standard hulls, the second author gave 
a similar construction for finitely generated algebras
over a field \cite{SchBCM}. The method, which itself relies on a
result of \cite{HHbigCM}, easily extends to the present situation, at
least for complete domains with algebraically closed residue field.
%The third application consists of some model-theoretic consequences
%of our main theorem for the existential theory of power series rings.

\subsection*{Conventions}
Throughout, $m$ and $n$ range over the set
$\mathbb{N} := \{0,1,2,\dots\}$
of natural numbers.
%We put ${\mathbb{N}}^{>0}:=\mathbb{N}\setminus\{0\}$.
By `ring' we always mean `commutative ring with
multiplicative identity $1$'.

\part{Faithfully Flat Lefschetz Extensions}

After some preliminaries on ultraproducts in \S\ref{s:UP} and on 
\filtered\ rings in \S\ref{s:EET}, in \S\ref{s:LP}
we prove the theorem from the introduction (in the form of
Theorem~\ref{T:maincan}). 
The construction of the desired Lefschetz extensions is achieved via
cylindrical approximation in equicharacteristic zero, which is a
corollary of Rotthaus' theorem \cite{Rot}, 
as we explain in \S\ref{s:UAA}. In \S\ref{s:app} we then discuss the
relationship between $R$ and the components of $\hull R$.

\section{Ultraproducts}\label{s:UP}

Let $\mathcal W$ be an infinite set. A \emph{non-principal 
ultrafilter} on $\mathcal W$ is a collection of infinite subsets of 
$\mathcal W$ which is closed under finite intersections and has the 
property that for any $W\subseteq \mathcal W$, either $W$ or its 
complement $\mathcal{W}\setminus W$ 
belongs to the collection. (One should think of the 
subsets $W$ which are in the ultrafilter as `big' and those not in it 
as `small'.) Given an infinite set $\mathcal W$, any collection of 
infinite subsets of $\mathcal W$ which is closed under finite 
intersections can be enlarged to a non-principal ultrafilter on 
$\mathcal W$. (See for instance \cite[Theorem~6.2.1]{Hod}.) 
Applying this to the collection of co-finite subsets of $\mathcal W$ implies
that on every infinite set there exists at least one non-principal ultrafilter.
With a few exceptions %(notably in \S\ref{Model-theoretic Applications}), 
we will always consider a fixed ultrafilter on a given infinite set, 
so there is no need to name the ultrafilter. Henceforth we call a set 
$\mathcal W$ endowed with some non-principal ultrafilter an \emph{\ultraset}.

In the remainder of this section we let $\mathcal W$ be an ultraset, and
we let $w$ range over $\mathcal W$. For each $w$ let $\seq Aw$ be a ring.  
The \emph{ultraproduct}  of the family $(A_w)$
(with respect to 
$\mathcal W$) is by definition the quotient of the product 
$\prod_w\seq Aw$ modulo the ideal $\mathcal I_{\text{null}}$ 
consisting of the sequences almost all of whose entries are zero. 
Here and elsewhere, a property is said to hold \emph{for almost all 
indices} if the subset of all $w$ for which it holds lies in the 
ultrafilter. We will often denote the ultraproduct of the family $(\seq Aw)$ by
        \begin{equation}\label{e:defup}
         \up {w\in\mathcal W}\seq Aw:=  \prod_{w\in\mathcal W} 
\seq Aw \bigg/\ \mathcal I_{\text{null}}.
        \end{equation}
Sometimes we denote such an ultraproduct simply by $\ul A$,
and we also speak, somewhat imprecisely, of 
`the ultraproduct of the $\seq Aw$' (with respect to $\mathcal W$). 
Given a 
sequence $a=(\seq aw)$ in $\prod_w \seq Aw$ we call its canonical 
image in $\ul A$ the \emph{ultraproduct} of the $a_w$ and 
denote it by 
        \begin{equation*}
        \ul a := \up{w\in\mathcal W} \seq aw.
        \end{equation*}
Similarly if $\seq{\mathbf{a}}w=(\seq{a_1}w,\dots,\seq{a_n}w)\in 
(\seq Aw)^n$ for each $w\in\mathcal W$ and $\ul{a_{i}}$ is the 
ultraproduct of the $\seq{a_i}w$ for $i=1,\dots,n$, then 
$\ul{\mathbf{a}}:=(\ul{a_1},\dots,\ul{a_n})\in (\ul A)^n$ is called 
the \emph{ultraproduct} %in $\ul A$ 
of the  $n$-tuples $\seq{\mathbf a}w$. 
If all $\seq Aw$ are the same, say equal to the ring $A$, then the 
resulting ultraproduct is called an \emph{ultrapower} of $A$ (with 
respect to $\mathcal W$), denoted by
        \begin{equation*}
        A^{\mathcal W}:=\up {w\in\mathcal W} A.
        \end{equation*}
The map $\delta_A\colon A\to A^{\mathcal W}$ which sends $a\in A$ to
the ultraproduct of the constant sequence with value $a$ is a ring embedding,
called the \emph{diagonal embedding} of $A$ into $A^{\mathcal W}$. We 
will always view $A^{\mathcal W}$ as an $A$-algebra via $\delta_A$.
Hence if $A$ is an $S$-algebra (for some ring $S$), then so is 
$A^{\mathcal W}$ in a natural way.

Let $\ul A$ and $\ul B$ be ultraproducts, with respect to  the same 
ultraset $\mathcal W$, of rings $\seq Aw$ and $\seq Bw$ respectively. 
If for each $w$ we have a map $\seq\varphi w\colon\seq Aw\to\seq 
Bw$, then we obtain a map $\ul\varphi\colon \ul A\to\ul B$,
called the \emph{ultraproduct} of the $\seq\varphi w$ (with respect to
$\mathcal W$),
by the rule
        \begin{equation*}
        a=\up w\seq aw\mapsto \ul\varphi(a):=\up w\seq\varphi w(\seq aw).
        \end{equation*}
(The right-hand side is 
independent of the choice of the $a_w$ such that $a=\up w\seq aw$.) Almost all
$\seq\varphi w$ are homomorphisms \iff\ $\ul\varphi$ is a homomorphism, and 
the $\seq \varphi w$ 
are injective (surjective) \iff\  $\ul\varphi$ is injective (surjective, 
respectively).  
%Let $\pi_{\seq Aw}\colon \seq Aw\times \seq Bw\to \seq Aw$ and 
%$\pi_{\seq Bw}\colon \seq Aw\times \seq Bw \to \seq Bw$ be the 
%natural projection maps and let $\ul{(\pi_A)}$ and $\ul{(\pi_B)}$ be 
%their respective ultraproducts. It is straightforward to verify that
        %\begin{equation*}
        %\ul{(\pi_A)}\times \ul{(\pi_B)}\colon \up w \seq Aw\times 
%\seq Bw\to \up w\seq Aw\times \up w\seq Bw
        %\end{equation*}
%is an isomorphism. 

These definitions apply in particular to ultrapowers, 
that is to say, the case where all $\seq Aw$ and $\seq Bw$ are equal to 
respectively $A$ and $B$. In fact, we then can extend them
to arbitrary $S$-algebras, for some base ring $S$. For instance, let 
$A$ and $B$ be $S$-algebras, and let $\varphi\colon A\to B$ be an 
$S$-algebra \homo. The \emph{ultrapower} of $\varphi$ (with respect 
to $\mathcal W$), denoted $\varphi^{\mathcal W}$, is the ultraproduct 
of the  $\seq\varphi w:=\varphi$. One easily verifies that 
$\varphi^{\mathcal W}\colon A^{\mathcal W}\to B^{\mathcal W}$ 
is again an $S$-algebra \homo. 
%Similarly, we 
%get an $S$-algebra isomorphism
%       \begin{equation}\label{eq:produp}
%       (\pi_A)^{\mathcal W}\times (\pi_B)^{\mathcal W}\colon 
%(A\times B)^{\mathcal W}\to A^{\mathcal W}\times B^{\mathcal W}
%       \end{equation}
%where $\pi_A$ and $\pi_B$ denote the natural projection maps.

The main model-theoretic fact about
ultraproducts is called \los. For most of our purposes
the following equational version suffices.

\begin{theorem}[Equational \los]\label{T:los}
%Let $\ul A$ be the ultraproduct of the rings $\seq Aw$ and let $\ul 
%{\mathbf a}$ be the ultraproduct of the $n$-tuples $\seq {\mathbf 
%a}w$ in $\seq Aw$. 
Given a system $\mathcal S$ of equations and 
inequalities
        \begin{equation*}
        f_1 =f_2=\dots=f_s=0,\  g_1\neq0, \ g_2\neq0,\dots,g_t\neq 0
        \end{equation*}
with $f_i,g_j\in\pol {\zet}{X_1,\dots,X_n}$, the tuple $ \ul {\mathbf 
a}$ is a solution of $\mathcal S$ in $A_\infty$ \iff\ almost all 
tuples $\seq{\mathbf a}w$ are solutions of $\mathcal S$ in $A_w$.
\end{theorem}

In particular it follows that any ring-theoretic property that can be 
expressed ``equationally'' holds for $\ul A$ \iff\ it holds for
almost all the rings $\seq Aw$. For example,
the ring $\ul A$ is reduced (a domain, a field) \iff\ 
almost all the rings $\seq Aw$ are reduced (domains, fields, respectively). 
All these statements can deduced from \los\ using appropriately chosen
systems $\mathcal S$. For instance, a ring $B$ is reduced \iff\ 
the system $X^2=0, X\neq 0$ (in the single indeterminate $X$) has no
solution in $B$. We leave the details of 
these and future routine applications of \los\ to the reader.
An example of a property which \emph{cannot} 
be transferred between $\ul A$ and the
$\seq Aw$ in this way is Noetherianity. (However,
$\ul A$ is Artinian of length $\leq l$ \iff\ almost all $\seq Aw$ are
Artinian of length $\leq l$, see \cite[Proposition 9.1]{JL}.)
Also note that if almost all $\seq Aw$ are algebraically closed fields, then
$\ul A$ is an algebraically closed field; the converse is false in general,
as \cite[Example~2.16]{JL} shows.

We refer to \cite{Chang-Keisler}, \cite{EkUP} or \cite{Hod} 
for in-depth discussions of ultraproducts. 
A brief review by the second author, adequate for 
our present needs, can be found in \cite[\S2]{SchNSTC}.
Using induction on the quantifier complexity of a formula, 
Theorem~\ref{T:los} readily 
implies the ``usual'' version of \los, stating that in
$\ul A$, the tuple $\ul {\mathbf a}$
satisfies a given (first-order) formula in the language of rings
\iff\ almost all $\seq {\mathbf a}w$ satisfy the same formula (in $\seq Aw$).
In particular, a sentence in the language of rings holds in $\ul A$ \iff\ it 
holds in almost all $\seq Aw$. 
Similarly, if for each $w$ we are given an endomorphism $\varphi_w\colon
A_w\to A_w$ of $A_w$, 
then its ultraproduct $\varphi_\infty$ is an endomorphism of
$A_\infty$, and a formula in the language of difference rings (= rings with
a distinguished endomorphism) holds for
the tuple $\ul {\mathbf a}$ in $(A_\infty,\varphi_\infty)$
\iff\ it holds for almost all $\seq {\mathbf a}w$ 
in $(A_w,\varphi_w)$.
On occasion, we invoke these stronger forms of \los.
(See for instance, \cite[Theorem 
9.5.1]{Hod} for a very general formulation.)

The ultraproduct construction also extends to more
general algebraic structures than rings. 
For example, if for each $w$ we are given
an $A_w$-module $M_w$, we may define 
$$M_\infty :=   \up {w\in\mathcal W}\seq Mw:= \prod_{w\in\mathcal W} 
\seq Mw\bigg/\ \mathcal M_{\text{null}}
$$
where $\mathcal M_{\text{null}}$ is the submodule of $\prod_w M_w$ consisting
of the sequences almost all of whose entries are zero. Then $M_\infty$ is
a module over $A_\infty$ in a natural way. If the $A_\infty$-module
$M_\infty$ is generated by $m_{1\infty},\dots,m_{s\infty}$, then the
$A_w$-module $M_w$ is generated by $m_{1w},\dots,m_{sw}$, for almost all $w$.
It is possible to formulate
a version of \los\ for modules. Since this will not be needed in the present 
paper,
let us instead
illustrate the functoriality inherent in the ultraproduct construction 
by establishing a fact which will be useful in \S\ref{s:LP}.
Suppose that for each $w\in\mathcal W$ we are
given an $A_w$-algebra $B_w$ and an $A_w$-module $M_w$. 

\begin{proposition}\label{P:Tor}
If $M_\infty$ has a resolution 
$$\cdots \xrightarrow{} (A_\infty)^{n_{i+1}}\xrightarrow{\varphi_i}
(A_\infty)^{n_i} \xrightarrow{\varphi_{i-1}} (A_\infty)^{n_{i-1}} \xrightarrow{} \cdots
\xrightarrow{\varphi_0} (A_\infty)^{n_0} \xrightarrow{} M_\infty\to 0$$
by finitely generated free $A_\infty$-modules $(A_\infty)^{n_i}$
and $B_\infty$ is coherent,
then as $B_\infty$-modules
\begin{equation}\label{e:tor}
\operatorname{Tor}_i^{A_\infty}(B_\infty,M_\infty) \iso
\bigl(\operatorname{Tor}_i^{A_w}(B_w,M_w)\bigr)_\infty
\end{equation}
for every $i\in\mathbb N$.
\end{proposition}

Here the module on the right-hand side of \eqref{e:tor} is the
ultraproduct of the $B_w$-modules $\operatorname{Tor}_i^{A_w}(B_w,M_w)$.
Before we begin the proof, first note that we may identify the free
$A_\infty$-module
$(A_\infty)^n$ with the ultraproduct $(A_w^n)_\infty$
of the free $A_w$-modules $A_w^n$ 
in a canonical way. Under this identification, if
$\mathbf{a}_{1w},\dots,\mathbf{a}_{mw}$ are elements of $A_w^n$, then
the $A_\infty$-submodule of $(A_\infty)^n$ generated by the
ultraproducts
$\mathbf{a}_{1\infty},\dots,\mathbf{a}_{m\infty}\in (A_\infty)^n$
of the $\mathbf{a}_{1w},\dots,\mathbf{a}_{mw}$, respectively,
corresponds to the ultraproduct $N_\infty$ of the $A_w$-submodules
$N_w := A_w\mathbf{a}_{1w}+\cdots+A_w\mathbf{a}_{mw}$ of $A_w^n$
(an $A_\infty$-submodule of $(A_w^n)_\infty$).
The canonical surjections $\pi_{w}\colon A_w^n\to A_w^n/N_w$
induce a surjection $\pi_{\infty}\colon 
(A_\infty)^n=(A_w^n)_\infty\to (A_w^n/N_w)_\infty$
whose kernel is  $N_\infty$. Hence we may identify $(A_\infty)^n/N_\infty$ and
$(A_w^n/N_w)_\infty$.

\begin{proof}[Proof \textup{(Proposition~\ref{P:Tor})}]
The $A_\infty$-linear maps $\varphi_i$ are given by certain
$n_{i-1}\times n_i$-matrices with entries in $A_\infty$. 
Hence each $\varphi_i$ is an ultraproduct $\varphi_i=\up w \varphi_{i,w}$ of
$A_w$-linear maps $\varphi_{i,w}\colon A_w^{n_{i+1}}\to
A_w^{n_i}$ with $\ker\varphi_i=(\ker\varphi_{i,w})_\infty$ and
$\operatorname{im}\varphi_i=(\operatorname{im}\varphi_{i,w})_\infty$.
Hence for given $i>0$ the complex
$$A_w^{n_{i+1}}\xrightarrow{\varphi_{i,w}}
A_w^{n_i} \xrightarrow{\varphi_{i-1,w}} A_w^{n_{i-1}} \xrightarrow{} \cdots
\xrightarrow{\varphi_{0,w}} A_w^{n_0} \xrightarrow{} M_w\to 0$$
is exact for almost all $w$, by \los.
On the other hand, tensoring the
free resolution of $M_\infty$ from above with $B_\infty$ we obtain the complex
$$\cdots \xrightarrow{} (B_\infty)^{n_{i+1}}\xrightarrow{\psi_i}
(B_\infty)^{n_i} \xrightarrow{} \cdots
\xrightarrow{\psi_{0}} (B_\infty)^{n_0} \xrightarrow{} B_\infty\tensor_{A_\infty}
M_\infty\to 0$$
where $\psi_i:=1\tensor\varphi_i$. (We identify $(B_\infty)^{n_{i}}$ and 
$B_\infty\tensor_{A_\infty}(A_\infty)^{n_{i}}$ as usual, for each $i$.)
%where each $S_\infty$-linear map $B_i$ is described by the matrix
%obtained by applying $\varphi_\infty$ to
%each entry of the matrix describing $A_i$. 
Writing each $\psi_i$ as an ultraproduct $\psi_i=\up w \psi_{i,w}$ of
$B_w$-linear maps $\psi_{i,w}\colon B_w^{n_{i+1}}\to
B_w^{n_i}$ yields, for given $i>0$, that
$\psi_{i,w}=1\tensor\varphi_{i,w}$ for almost all $w$, hence
$\operatorname{Tor}_i^{A_w}(B_w,M_w)\iso\ker \psi_{i-1,w}/\operatorname{im} \psi_{i,w}$
for almost all $w$. Since $B_\infty$ is coherent, the $B_\infty$-module 
$\ker\psi_{i-1}$ is finitely generated, and we get
$$\operatorname{Tor}_i^{A_\infty}(B_\infty,M_\infty) = \ker \psi_{i-1}/
\operatorname{im} \psi_i \iso \bigl(\ker \psi_{i-1,w}/\operatorname{im} \psi_{i,w}\bigr)_\infty.$$
This proves the case $i>0$ of the proposition. Using the remarks
preceding the proof it is easy to show that
$B_\infty\tensor_{A_\infty} M_\infty \iso 
\bigl( B_w \tensor_{A_w} M_w \bigr)_\infty$, proving the case $i=0$.
\end{proof}

%The following fact will be useful later.
%
%\begin{lemma}\label{Flatness}
%Given rings $\seq Aw$ indexed by an ultraset $\mathcal W$, the canonical homomorphism $$ \prod_{w\in\mathcal W} 
%A_w \to \up {w\in\mathcal W}{\seq Aw}$$ sending a sequence $(a_w)$ to its ultraproduct $a_\infty$, is flat.
%%\rom{(In particular, if $R$ is $(\alpha,\beta)$-super coherent, so is $R^*$.)}
%\end{lemma}
%\begin{proof}
%Let $A$ (respectively, $\ul A$) be the (ultra-) product of the $A_w$. We have $A_\infty = A/\mathcal{I}_{\operatorname{null}}$, where 
%$\mathcal{I}_{\operatorname{null}}$
%denotes the ideal of $A$ consisting of all
%sequences $a=(a_w)\in A$ such that $a_w=0$ for almost all $w$.
%(See \eqref{e:defup}.)
%By \cite[Theorem~1.2.15]{Glaz}, the quotient $A/{\mathcal I}_{\operatorname{null}}$ is a flat 
%$A$-module if, and only if,
%for every $a\in \mathcal{I}_{\operatorname{null}}$ there is some
%$c\in \mathcal{I}_{\operatorname{null}}$ with $(1-c)a=0$.
%To see that this condition is satisfied, 
%let $a=(a_w)\in \mathcal{I}_{\operatorname{null}}$, so
%$\Delta:=\{w:a_w=0\}$ is a member of the ultrafilter of $\mathcal W$. Define
%$c=(c_w)\in A$ by $c_w:=0$ if $w\in\Delta$ and $c_w:=1$ if
%$w\notin\Delta$. Then $c\in \mathcal{I}_{\operatorname{null}}$ and $(1-c)a=0$ 
%as required. 
%\end{proof}

\subsection{Lefschetz rings}\label{s:lefschetz}
An ultraproduct $\ul A$ of rings $A_w$ with respect to an ultraset $\mathcal W$ 
will be called \emph{Lefschetz} (with respect to $\mathcal W$) if
almost all of the $A_w$ are of prime characteristic and $\ul A$ is
of characteristic zero. 
(The condition on $\op{char}(\ul A)$ holds precisely if for each prime number
$p$, the set $\{w: \op{char}(\seq Aw)=p\}$ does not belong to
the ultrafilter of $\mathcal W$.)
A Lefschetz field is a Lefschetz ring that happens to be a field; in this
case almost all $A_w$ are fields. 
The following proposition is a well-known consequence of
\los. We let $p$ range over the set of prime numbers.
As usual $\mathbb F_p$ denotes the field with $p$ elements
and $\ac Fp$ its algebraic closure.
%An ultraset $\mathcal W$ whose underlying set is countable is called
%$\aleph_0$-regular if there exists an infinite subset
%$\mathcal V$ of the ultrafilter of $\mathcal W$ with the property
%that each $w$ 
%belongs to only finitely
%many sets in $\mathcal V$. Any countably infinite set $\mathcal W$
%can be equipped with a
%non-principal ultrafilter such that $\mathcal W$ is
%$\aleph_0$-regular, see \cite[Proposition 4.3.5]{Chang-Keisler}.

\begin{proposition}\label{P:acf}
There is a \emph{(}non-canonical\emph{)} isomorphism between the 
field of complex numbers $\mathbb C$ and an ultraproduct of the $\ac 
Fp$. 
\end{proposition}
\begin{proof}
%Let $F$ be an \acf\  of \ch\ zero having cardinality 
%$2^\lambda$, where $\lambda\geq\aleph_0$. 
%For each prime number $p$ let $\seq Fp$ be an \acf\ of 
%\ch\ $p$ and of  cardinality $\lambda$. In case $F=\mathbb C$ we can 
%choose $\seq Fp:=\ac Fp$.  
Equip 
the set $\mathcal P$ 
of prime numbers with a non-principal ultrafilter and let $\ul F$ be the 
ultraproduct  of the $\ac Fp$ with respect to the ultraset $\mathcal P$. 
%We can express in an equational way 
%that a ring is a field and that this field is algebraically closed.
By the remarks following Theorem~\ref{T:los}, we see that
$\ul F$ is an \acf.  
Since $l$ is a unit in $\seq 
Fp$, for every prime $l$ distinct  from $p$, it is a unit in $\ul F$, 
by \los. Consequently, $\ul F$ has \ch\ zero. 
The cardinality of $\ul F$ 
%has again 
%cardinality $2^\lambda$; 
is that of the continuum; see 
\cite[Proposition~4.3.7]{Chang-Keisler}.   Any 
two \acf{s} of \ch\ zero, of the same uncountable cardinality, are 
isomorphic, since they have the same transcendence degree over 
$\mathbb Q$.
%(Steinitz' Theorem; for a model theoretic proof using quantifier 
%elimination, or, equivalently, Chevalley's Theorem, see 
%\cite[Corollary 4.5.7]{Hod}.)
Hence $F\iso \ul F$.
\end{proof}

\begin{remark}
Note that that the particular choice of non-principal ultrafilter
on $\mathcal P$ used in the proof above is irrelevant.
The same argument may also be employed to show, more generally: 
every \acf\ of \ch\ zero of uncountable
cardinality $2^\lambda$ \textup{(}for some infinite cardinal 
$\lambda$\textup{)} is isomorphic
to a Lefschetz field $F_\infty$ with respect to
$\mathcal P$ all of whose components $F_p$ are
algebraically closed fields of characteristic $p$.
It follows that every field of characteristic zero can be embedded into
a Lefschetz field all of whose components are
algebraically closed fields. Moreover,
under the assumption of the Generalized Continuum Hypothesis ($2^\lambda=
\lambda^+$ for all infinite cardinals $\lambda$) 
every uncountable algebraically closed field of characteristic zero is
Lefschetz.
\end{remark}

The following class of Lefschetz rings will be of special interest to us:

\begin{definition}\label{D:analytic-Lefschetz}
A Lefschetz ring $\ul A$ (with respect to the ultraset $\mathcal W$) will be
called an \emph{analytic Lefschetz ring}
(with respect to $\mathcal W$) if
almost all of the $A_w$ are complete Noetherian local rings
of prime equicharacteristic
with algebraically closed residue field. 
Let $A_\infty$ and $B_\infty$ be analytic Lefschetz rings.
An ultraproduct  $\ul\varphi\colon \ul A\to \ul B$
of local ring homomorphisms $\varphi_w\colon A_w\to B_w$ will be called 
a \emph{homomorphism of analytic 
Lefschetz rings}\/ (with respect to $\mathcal W$).
\end{definition}

By \los\ every analytic Lefschetz ring is a local ring, and every homomorphism
of analytic Lefschetz rings is a local homomorphism of local rings.
If $A=A_\infty$ is a Lefschetz ring (an analytic Lefschetz ring)
with respect to $\mathcal W$
and $I$ a finitely generated proper ideal of $A$,
then $A/I$ is isomorphic to a
Lefschetz ring (an analytic Lefschetz ring, respectively)
with respect to the same ultraset $\mathcal W$. 
Hence if the maximal ideal of the analytic
Lefschetz ring $A_\infty$ is finitely
generated, then
the residue field of  $A_\infty$ may be identified with 
the ultraproduct $\ul K$ of the residue fields $K_w$ of $A_w$ 
in a natural way. In particular $\ul K$ is
itself Lefschetz and algebraically closed.

%\begin{proof}
%Suppose $A=\up{w\in\mathcal W} \seq Aw$, where almost all the rings 
%$A_w$ are
%of prime characteristic, and $I=(a_1,\dots,a_n)$. For
%each $i$ choose $a_{iw}\in A_w$ whose ultraproduct is $a_i$. We put
%$B_w:= A_w/(a_{1w},\dots,a_{nw})$ for each $w$. Since $I\neq A$, by
%\los\ almost all $B_w$ are non-trivial rings of
%prime characteristic. Hence  $B:=\up{w\in\mathcal W} B_w$ is a Lefschetz ring
%with respect to $\mathcal W$.
%The canonical surjections $\pi_w\colon A_w\to B_w$
%induce a surjection $\pi_\infty\colon A\to B$. By \los\ it follows that
%$I=\ker\pi_\infty$, hence $A/I\cong B$. If almost all $A_w$ are
%complete Noetherian local with algebraically residue field of
%prime characteristic, then the same is true for the $B_w$, and
%almost all $\pi_w$ are local homomorphisms.
%\end{proof}

\begin{example}\label{E:analytic-Lefschetz}
For fixed $n$ let $A_w := K_w[[X_1,\dots,X_n]]$ be the ring of formal
power series in indeterminates $X_1,\dots,X_n$
over an algebraically closed field $K_w$
of characteristic $p(w)>0$. If for every integer $p>0$, almost all $p(w)$
are $>p$, then the ultraproduct $\ul A$ of the $A_w$ has characteristic zero
and hence is an analytic
Lefschetz ring.
In this example, $A_\infty$ is a 
$K_\infty$-algebra in a natural way. 
In general, if $K$ is a
Lefschetz field (with respect to $\mathcal W$)
and $K\to A$ is a homomorphism
of analytic Lefschetz rings (with respect to $\mathcal W$), then we call $A$
an \emph{analytic Lefschetz $K$-algebra} (with respect to $\mathcal W$). 
The analytic Lefschetz $K$-algebras with respect to $\mathcal W$
form a category whose morphisms are the
homomorphisms of analytic
Lefschetz rings with respect to $\mathcal W$ 
that are also $K$-algebra homomorphisms.
\end{example}

We will on occasion use the following construction.

\subsection{Ultraproducts of polynomials of bounded degree}\label{s:ulpol}
%Let $\ul A$ be an ultraproduct of rings $\seq  Aw$. 
Let 
$X=\rij Xn$ be a tuple of indeterminates and let $\ul B$ be the 
ultraproduct of the polynomial rings
$\pol{\seq Aw}X$. Taking the ultraproduct of the 
natural \homo{s} $\pol \zet X\to \pol{\seq Ap}X$ gives a canonical 
\homo\ $\pol \zet X^{\mathcal W}\to \ul B$. We will continue to write 
$X_i$ for the image of $X_i$ under this \homo. On the other 
hand, $\ul A$ is a subring of $\ul B$. Using \los, we see that
$X_1,\dots,X_n$ remain algebraically independent over $\ul A$, so 
that we have  in fact a canonical embedding $\pol{\ul A}X\subseteq \ul 
B$.
Suppose now we are given, for some $d\in\mathbb{N}$ and each $w$, a  polynomial
        \begin{equation*}
        \seq Qw=\sum_\nu \seq{a_\nu}w  X^\nu \in\pol{\seq Aw}X
\qquad (\seq{a_\nu}w\in \seq Aw)
        \end{equation*}
of degree at most $d$. Here the sum ranges 
over all multi-indices $\nu=\rij\nu n \in {\mathbb N}^n$   with $d\leq 
|\nu|:=\nu_1+\cdots+\nu_n$, and as usual $X^\nu$ is shorthand 
for $X_1^{\nu_1}\cdots X_n^{\nu_n}$.  Let $\ul{a_\nu}\in \ul A$ be 
the ultraproduct of the $\seq{a_\nu}w$ and put
        \begin{equation*}
        \ul Q:=\sum_\nu \ul{a_\nu} X^\nu,
        \end{equation*}
a polynomial in $\pol {\ul A}X$ of degree  $\leq d$.  
(The polynomial $\ul Q$ has degree $d$ \iff\ almost all $\seq Qw$ have degree 
$d$.) We call $\ul Q$ the 
\emph{ultraproduct} of the $\seq Qw$. This is justified by the fact 
that the image of $\ul Q$ under the canonical embedding $\pol{\ul 
A}X\subseteq \ul B$ is the ultraproduct of the $\seq Qw$. In contrast, 
ultraproducts of polynomials of unbounded degree do no longer belong 
to the subring $\pol{\ul A}X$.

\section{Embeddings and Existential Theories}\label{s:EET}

In this section, we want to address the following question:
\emph{given $S$-algebras $A$ and $B$,  when does there  exist an
$S$-algebra \homo\ $A\to B$?} If one is willing to replace $B$  by
some ultrapower, then a simple criterion exists (Corollary~\ref{C:LP}
below). Although this does not solve the question raised above, it
suffices for showing that a faithfully flat Lefschetz extension
exists (see \S\ref{s:LP}). To obtain the desired 
functoriality, we need a
\filtered\ version of this result, which we now explain.

\subsection{\Filtered\ rings}\label{s:fil}
A \emph\filtered\ ring  is a ring $R$ together with a \emph{nest}\/ of
subrings, that is, an ascending chain of subrings
        \begin{equation*}
        R_0\subseteq R_1\subseteq\cdots\subseteq R_n\subseteq\cdots
        \end{equation*}
of $R$ whose union equals $R$. We agree that whenever $R$ is a 
\filtered\ ring, we denote the  subrings in the nest by $R_n$, and we 
express this by saying that  $R=(R_n)$ is a \filtered\ ring.  Every 
ring $R$ can be made into a  \filtered\ ring using the nest with 
$R_n:=R$ for all $n$. (We say that $R$ is \emph{trivially \filtered}.)

Let $R=(R_n)$ and $S=(S_n)$ be \filtered\
rings.  A homomorphism $\varphi\colon S\to R$ is called a
\emph{homomorphism of \filtered\ rings} if $\varphi(S_n)\subseteq  R_n$
for all $n$. Alternatively, we say that  $R$ is a
\emph{\filtered\ $S$-algebra} (via $\varphi$).
Note that in this case,  $R_n$ is naturally an
$S_n$-algebra, for every $n$. An $S$-algebra homomorphism $R\to R'$
between \filtered\ $S$-algebras $R$ and $R'$ which is a homomorphism
of \filtered\ rings is called a \emph{homomorphism of \filtered\
$S$-algebras.} If $R\to R'$ is injective, we may identify $R$ with a
subalgebra of $R'$, and we refer to this situation by calling $R$ a
\emph{\filtered\ $S$-subalgebra} of $R'$. A bijective homomorphism
of \filtered\ rings (\filtered\ $S$-algebras)
is called an \emph{isomorphism}\/ of \filtered\ rings
(\filtered\ $S$-algebras, respectively).

\begin{example}\label{e:nested}
Let $L$ be a field and $Y_0,Y_1,\dots$ an infinite sequence of
finite (possibly empty) tuples $Y_n=(Y_{n1},\dots,Y_{nk_n})$ of 
distinct indeterminates,
$k_n\in\mathbb N$. For each $n$ put 
\begin{align*}
S_n &:= L[[Y_0]][Y_1,\dots,Y_n], \\
R_n &:= L[[Y_0]][Y_1,\dots,Y_n]^{\operatorname{alg}},\\
A_n &:= L[[Y_0,\dots,Y_n]].
\end{align*}
Here and elsewhere, given a domain $D$ and a finite tuple $Y$ of
indeterminates we denote by
$\algpow DY$  the subring of  $\pow DY$ consisting of all elements
which are algebraic over $\pol DY$. (If $D$ is an excellent domain, then
$\algpow DY$ is 
equal to the Henselization $\hens{\pol DY}$ of $\pol DY_{1+YD[Y]}$ 
at the ideal 
generated by the indeterminates, see 
\cite[p.~126]{Raynaud}.)
We make the subrings $S:=\bigcup S_n$,  $R:=\bigcup_n R_n$ and
$A:=\bigcup_n A_n$ of $L[[Y_0,Y_1,\dots]]$ into nested rings with nests
$(S_n)$, $(R_n)$ and $(A_n)$, respectively.
Then $R$ is 
a nested $S$-subalgebra of the nested $S$-algebra $A$.
(This example will play an important role in \S\ref{s:UAA}.) 
\end{example}

Let $S$ be a \filtered\ ring and  $V$   a \filtered\ $S$-algebra with 
nest $(V_n)$. We
say that $V$ is \emph{of finite type} (over $S$) if each $V_n$ is a 
finitely generated
$S_n$-algebra, and for some $n_0$, each $V_n$ with $n\geq
n_0$ is the $S_n$-subalgebra of $V$ generated by $V_{n_0}$, that is to say,
$V_n=\pol{S_n}{V_{n_0}}$. Choose $n_0$ minimal  with this property.
Clearly, all the knowledge about $V$ is already contained in the
initial chain $V_0\subseteq V_2\subseteq\cdots\subseteq
V_{n_0}$,  and consequently, we  refer to it as the \emph{relevant
part} of $V$, and to  $n_0$ as its \emph{length}.

\subsection{\Filtered\ equations and \filtered\ algebras of finite type}\label{s:ft}
In the following let $S$ be a nested ring.
A \emph{\filtered\ system of polynomial equations} with coefficients
from $S$ is a finite sequence $\mathcal S$ of systems of polynomial
equations
\begin{equation}\label{e:nestedeq}
\begin{array}{lclclcl}
   P_{00}(Z_0)     &=& \cdots& = &P_{0k}(Z_0)&=&0, \\
   P_{10}(Z_0,Z_1) &=& \cdots& = &P_{1k}(Z_0,Z_1)&=&0, \\
   \hfill\vdots\hfill & & & & \hfill\vdots\hfill  \\
   P_{n0}(Z_0,\dots,Z_n)&=& \cdots& = &P_{nk}(Z_0,\dots,Z_n)&=&0,
\end{array}
\end{equation}
for some $n$ and $k\in\mathbb{N}$, some tuples $Z_i=(Z_{i1},\dots,Z_{ik_i})$ of
indeterminates over $S$ (where $k_i\in\mathbb N$)
and some polynomials $P_{ij}\in \pol{S_i}{Z_0,\dots,Z_i}$.
%The tuple $Z:=(Z_0,\dots,Z_n)$ is called the
%\emph{tuple of variables} of $\mathcal S$;
% often we will not specify the particular partitioning of $Z$.
Given a nested $S$-algebra $A$, a  
tuple $(\mathbf{a}_0,\dots, \mathbf{a}_n)$ with $\mathbf{a}_i\in 
(A_i)^{k_i}$ is called a \emph{\filtered\ solution of $\mathcal S$} 
in $A$ if $P_{ij}(\mathbf{a}_0,\dots,\mathbf{a}_i)=0$ for all 
$i=0,\dots,n$ and $j=0,\dots,k$. Similarly, given an ideal $\id$ of 
the nested $S$-algebra $A$, we call $(\mathbf{a}_0,\dots, \mathbf{a}_n)$ an 
\emph{approximate \filtered\ solution of $\mathcal S$ modulo $\id$} 
in $A$, if $P_{ij}(\mathbf{a}_0,\dots,\mathbf{a}_i)\equiv0\bmod
\id$ for all $i,j$.

Let $V$ be a \filtered\ $S$-algebra of finite type and let $n$ be the 
length of its relevant part. For $i\leq n$, choose  tuples $\mathbf 
a_i\in V_i^{k_i}$, such that each $V_i$ is generated as an 
$S_i$-algebra by $\mathbf a_0,\dots,\mathbf a_i$. Let 
$P_{i1},\dots,P_{ik}$ be generators of the kernel of the 
$S_i$-algebra \homo\ $\pol{S_i}{Z_0,\dots,Z_i}\to V$ given by 
$Z_0\mapsto \mathbf a_0,\dots, Z_i\mapsto\mathbf a_i$. In particular
         \begin{equation*}
         V_i\iso 
\pol{S_i}{Z_0,\dots,Z_i}/(P_{i0},\dots,P_{ik})\pol{S_i}{Z_0,\dots,Z_i}. 
         \end{equation*}
The system of equations  $P_{00}=\cdots =P_{nk}=0$ form a \filtered\ 
system of polynomial equations with coefficients from $S$, called a 
\emph{defining \filtered\ system of equations} for $V$. (It depends on 
the choice of generators $\mathbf a_i$.) Note that the generating 
tuple $(\mathbf a_0,\dots,\mathbf a_n)$ is a \filtered\ solution of 
this system in $V$. Conversely, any \filtered\ system of polynomial 
equations with coefficients from $S$ together with a \filtered\ 
solution in some \filtered\ $S$-algebra $B$ gives rise to a 
\filtered\ $S$-subalgebra of $B$ of finite type.

Given an ultraset $\mathcal U$ and a nested $S$-algebra $B$
we consider the ultrapower
$B_n^{\mathcal U}$ as an $S_n$-subalgebra
of $B^{\mathcal U}$ in the natural way. We
make the $S$-subalgebra $\bigcup_{n}
B_n^{\mathcal U}$ of the ultrapower $B^{\mathcal U}$ into a
\filtered\ $S$-algebra by means of the nest $(B_n^{\mathcal U})$.
We denote this \filtered\ $S$-algebra by $B^{\langle\mathcal U\rangle}$.
The main result of this section is the following criterion for the
existence of a \homo\ of \filtered\ $S$-algebras from a nested $S$-algebra
$A$ to an ultrapower of $B$.

\begin{theorem}\label{T:LPfil}
Let $A$ and $B$ be \filtered\
$S$-algebras.  If each $S_n$ is Noetherian, then the following are
equivalent:
\begin{enumerate}
\item\label{i:eqhomo} every \filtered\ system of polynomial equations
with coefficients from $S$ which has a \filtered\ solution in $A$
also has one in $B$;
\item\label{i:lochomo} for every  \filtered\ $S$-subalgebra of finite
type $V$ of $A$, there exists a \homo\ of \filtered\ $S$-algebras
$\varphi_V\colon V\to B$;
\item\label{i:ulhomo} there exists  a homomorphism of \filtered\
$S$-algebras $\eta\colon A\to B^{\langle\mathcal U\rangle}$, for some
ultraset $\mathcal U$.
\end{enumerate}
%Moreover, if all the $\varphi_V$ in \eqref{i:lochomo} can be taken
%injective, then
%$\eta$ in \eqref{i:ulhomo} can be taken injective.
\end{theorem}
\begin{proof}
Suppose that \eqref{i:eqhomo} holds, and let $V$ be a  \filtered\ 
$S$-subalgebra  of finite type of $A$. Suppose 
$V_0\subseteq\cdots\subseteq V_n$ is the relevant part of $V$  (so 
that $V_m=\pol {S_m}{V_n}$ for all $m\geq n$).  Let $\mathcal S$ be a 
defining \filtered\ system of equations of $V$ and let $(\mathbf 
a_0,\dots,\mathbf a_n)$ with $\mathbf a_i\in (A_i)^{k_i}$ be the 
\filtered\ solution in $A$ arising from a generating set of $V$ over 
$S$ (see \S\ref{s:ft}). By assumption, there exists a \filtered\ 
solution  $(\mathbf{b}_0,\dots, \mathbf{b}_n)$ of $\mathcal S$ with 
$\mathbf{b}_i\in (B_i)^{k_i}$   for all $i$. Hence the
$S_n$-algebra homomorphism  $\pol{S_n}{Z_0,\dots,Z_n}\to B_n$ given
by $Z_i\mapsto \mathbf{b}_i$ for $i=0,\dots,n$ factors through an 
$S_n$-algebra homomorphism $\varphi_V\colon V_n\to B_n$ with 
$\varphi_V(V_i)\subseteq B_i$ for all $i$. Since $V_m=\pol{S_m}{V_n}$ 
for $m\geq n$, we can extend this to a \homo\ of \filtered\ 
$S$-algebras $V\to B$, proving implication 
\eqref{i:eqhomo}~$\Rightarrow$~\eqref{i:lochomo}.

Assume next that \eqref{i:lochomo} holds. Let $\mathcal U$ be the
collection of all \filtered\ $S$-subalgebras  of finite type of $A$
(an infinite set). For each finite subset 
$E=\bigl\{(a_1,n_1),\dots,(a_k,n_k)\bigr\}$ of  $A\times\nat$ let 
$\langle E\rangle$ be the subset of $\mathcal U$ consisting of all 
\filtered\ $S$-subalgebras $V=(V_n)$  of finite type of $A$ with
$a_i\in V_{n_i}$ for all $i$. Any finite intersection of sets of the
form $\langle E\rangle$ is again of that form. Hence we can find a
non-principal ultrafilter on $\mathcal U$ containing each $\langle
E\rangle$, where $E$ runs over all finite subsets of $A\times\nat$. 
For each $V\in\mathcal U$, let $\widetilde{\varphi}_V\colon A\to
B$ be the map which coincides with $\varphi_V$ on $V$ and which is
identically zero outside $V$. (This is of course no longer a
homomorphism.) Define $\eta\colon A\to B^{\mathcal U}$ to be the
restriction to $A$ of the ultraproduct of the
$\widetilde{\varphi}_V$. In other words,
         \begin{equation*}
         \eta(a):= \up {V\in\mathcal U} \widetilde{\varphi}_V(a) \qquad
\text{for $a\in A$.}
         \end{equation*}
It remains to verify that the image of $\eta$ lies inside
$B^{\langle\mathcal U\rangle}$ and that the induced \homo\ $A\to
B^{\langle\mathcal U\rangle}$ is a \homo\ of \filtered\ $S$-algebras.
For $a,b\in A_n$, we have for each $V\in \langle
\{(a,n),(b,n)\}\rangle$ that $\widetilde{\varphi}_V(a)=\varphi(a)$
and $\widetilde{\varphi}_V(b)=\varphi(b)$ lie in $B_n$ and
         \begin{equation*}
         \widetilde{\varphi}_V(a+b)=\varphi_V(a+b)=\varphi_V(a)
+\varphi_V(b) = \widetilde{\varphi}_V(a) + \widetilde{\varphi}_V(b).
         \end{equation*}
Since this holds for almost all $V$, we get that $\eta(a),\eta(b)\in
B_n^{\mathcal U}$ and $\eta(a+b)=\eta(a)+\eta(b)$. In particular,
the image of $\eta$ lies inside $B^{\langle\mathcal U\rangle}$. By a
similar argument, one also shows that $\eta(ab)=\eta(a)\eta(b)$ and
$\eta(sa)=s\eta(a)$ for $s\in S$. 
%Note that if all $\varphi_V$ are
%injective, then $\eta$ is injective. 
We have shown
\eqref{i:lochomo}~$\Rightarrow$~\eqref{i:ulhomo}.

Finally, suppose that $\eta\colon A\to B^{\langle\mathcal U\rangle}$
is a homomorphism of \filtered\ $S$-algebras, for some ultraset
$\mathcal U$. Suppose moreover that we are given a \filtered\ system
$\mathcal S$ of polynomial equations with coefficients from $S$ as above,
which has a \filtered\ solution $(\mathbf{a}_0,\dots,\mathbf{a}_n)$ in $A$.
Then $\bigl(\eta(\mathbf{a}_0),\dots,\eta(\mathbf{a}_n)\bigr)$ is a
\filtered\ solution of $\mathcal S$ in the \filtered\ $S$-algebra
$B^{\langle\mathcal U\rangle}$. Using \los\ it follows that $\mathcal
S$ has a \filtered\ solution in $B$. This shows 
\eqref{i:ulhomo}~$\Rightarrow$~\eqref{i:eqhomo}.
\end{proof}

Applying the theorem to trivially \filtered\ rings we obtain the
following partial answer to the question raised at the beginning of this
section. It is an incarnation of a
model-theoretic principle (originating with Henkin \cite{Henkin})
which has proven to be useful in other situations related to Artin
Approximation; for instance, see \cite[Lemma~1.4]{BDL83}
and \cite[Lemma~12.1.3]{vdD90}.

\begin{corollary}\label{C:LP}
Let $S$ be a Noetherian ring and let $A$ and $B$ be  $S$-algebras.
The following are equivalent:
\begin{enumerate}
\item\label{i:E+} every \textup{(}finite\textup{)} system of
polynomial equations with coefficients from $S$ which is solvable in
$A$, is solvable in $B$;
\item\label{i:locLP} for each finitely generated $S$-subalgebra $V$
of $A$, there exists an $S$-algebra \homo\ $\varphi_V\colon V\to B$;
\item\label{i:LP} there exists an \ultraset\ $\mathcal U$  and an
$S$-algebra \homo\ $\eta\colon A\to B^{\mathcal U}$. \qed
\end{enumerate}
%Moreover, if all the $\varphi_V$ in \eqref{i:locLP} can be taken
%injective, then so can $\eta$   in \eqref{i:LP}.\qed
\end{corollary}

We finish this sections with some remarks 
about Theorem~\ref{T:LPfil} and its corollary above.

\begin{remark}\label{R:Noe}
Only the proof of the implications
\eqref{i:eqhomo}~$\Rightarrow$~\eqref{i:lochomo} and
\eqref{i:E+}~$\Rightarrow$~\eqref{i:locLP}
used the assumption that each $S_n$ (respectively, $S$) is
Noetherian. These implications do hold without the Noetherian
assumption, provided we allow for infinite systems  (in finitely many
variables)  in  \eqref{i:eqhomo}  and \eqref{i:E+} respectively.
\end{remark}

\begin{remark}\label{R:cofinal}
In the proof of \eqref{i:lochomo}~$\Rightarrow$~\eqref{i:ulhomo} we may
replace the underlying set of the ultraset 
$\mathcal U$ by any cofinal collection
of nested $S$-subalgebras of finite type of $A$.
\end{remark}

\begin{remark}\label{R:can}
We can strengthen \eqref{i:ulhomo} and \eqref{i:LP} by making $\eta$
canonical, that is to say, independent of the choice of $S$-algebra
\homo{s} $\varphi_V$. Let us just give the argument in the
non-\filtered\ case. Replace the above index set $\mathcal U$ by the
set $\mathcal A$ of all $S$-algebra \homo{s} $\varphi\colon V\to B$
whose domain $V$ is a  finitely generated $S$-subalgebra of $A$.
Given a finite subset $E$ of $A$, let  $\langle E\rangle$ be the
subset of all $\varphi\in\mathcal A$ whose domain contains $E$. If we
assume \eqref{i:locLP} and $A$ is not finitely generated,  then
$\mathcal A$ is infinite and no  $\langle E\rangle$ is empty, so that
we can choose a non-principal ultrafilter on $\mathcal A$ which
contains all the  $\langle E\rangle$, for $E$ a finite subset of $A$.
The remainder of the construction is now the same. Namely,    define
$\eta\colon A\to B^{\mathcal A}$ to be the restriction to $A$ of the
ultraproduct of all $\widetilde{\varphi}$, where for each
$\varphi\in\mathcal A$ we let  $\widetilde{\varphi}\colon A\to B$ be
the extension by zero of $\varphi$. The same argument as above then
yields that $\eta$ is an $S$-algebra \homo.
\end{remark}

\begin{remark}\label{r:emb}
We also have criteria 
for $A$ to \emph{embed}\/ into an ultrapower of $B$: under the same
assumptions as in Corollary~\ref{C:LP},
the following are equivalent:
\begin{enumerate}
\item\label{i:emb1} \textit{every \textup{(}finite\textup{)} system of
polynomial equations and inequalities
with coefficients from $S$ which is solvable in
$A$, is solvable in $B$;}
\item\label{i:emb2} 
\textit{given a finitely generated $S$-subalgebra $V$
of $A$ and finitely many non-zero elements $a_1,\dots,a_n$ of $V$ there exists an 
$S$-algebra homomorphism $V\to B$ sending each $a_i$ to a non-zero element of
$B$;}
\item \textit{there exists an \ultraset\ $\mathcal U$  and an
embedding $A\to B^{\mathcal U}$ of $S$-algebras.}
\end{enumerate}
In particular, 
if all the $\varphi_V$ in \eqref{i:locLP} can be taken
injective, then so can $\eta$   in \eqref{i:LP}. Similar criteria may be
formulated in the general nested case.
We leave the proof (which is analogous to the proof of Theorem~\ref{T:LPfil})
to the reader.
\end{remark}

In the next remarks (not essential later)
we assume that the reader is familiar with basic notions of
model theory; see \cite{Chang-Keisler} or \cite{Hod}.

\begin{remark}\label{R:emb}
The language $\mathcal L(S)$ of $S$-algebras (in the sense of
first-order logic) consists of the language $\mathcal L=\{ {0}, {1},
{+}, {-}, {\cdot} \}$ of rings augmented by a unary function symbol
$s^\times$, for each $s\in S$. We construe each $S$-algebra as an
$\mathcal L(S)$-structure by interpreting the ring symbols as usual
and $s^\times$ as multiplication by $s$. We can then reformulate
\eqref{i:E+} in more model-theoretic terms as:
\begin{itemize}
\item[\eqref{i:E+}$'$] \textit{$B$ is a model of the positive
existential theory of $A$ in the language $\mathcal L(S)$}.
\end{itemize}
Similarly \eqref{i:emb1} may be replaced by
\begin{itemize}
\item[\eqref{i:emb1}$'$]\textit{$B$ is a model of the 
\textup{(}full\textup{)} existential $\mathcal L(S)$-theory of $A$.}
\end{itemize}
\end{remark}

\begin{remark}\label{R:sat}
Suppose that $B$ is $|A|$-saturated (as an $\mathcal L(S)$-structure).
Then to \eqref{i:E+}--\eqref{i:LP} in Corollary~\ref{C:LP} we may add
the equivalent statement
\begin{itemize}
\item[\eqref{i:LP}$'$] \textit{There exists an $S$-algebra homomorphism
$A\to B$.}
\end{itemize}
For a proof see for instance \cite[Theorem
10.3.1]{Hod}. The assumption on $B$ is satisfied
if $S$ (and hence $\mathcal L(S)$) is countable, $A$ has
cardinality at most $\aleph_1$, and $B$ is an ultraproduct of
a countable family of $S$-algebras with respect 
to a non-principal ultrafilter. 
(See \cite[Theorem 6.1.1]{Chang-Keisler}.)
If, on the other hand, $B$ is $\aleph_0$-saturated,
then in Remark~\ref{r:emb} we may replace \eqref{i:emb2} with
\begin{itemize}
\item[\eqref{i:emb2}$'$] \textit{For every finitely generated $S$-subalgebra
$V$ of $A$
there exists an embedding $V\to B$ of $S$-algebras.}
\end{itemize}
\end{remark}

\section{Artin Approximation and Embeddings in Ultraproducts}\label{s:UAA}

In this section, $K$ is a field which is the 
ultraproduct of fields $\seq Kp$ (not necessarily algebraically 
closed nor of different \ch s) with respect to an ultraset $\mathcal 
P$. In most applications, $\mathcal P$ will have as underlying  set 
the set of prime numbers and each $\seq Kp$ will have \ch\ $p$.
%Recall that any  uncountable \acf\ $K$ of \ch\ zero occurs in this 
%way, by Proposition~\ref{P:acf}. 
For a finite tuple 
$X=\rij Xn$  of indeterminates, we put
        \begin{equation*}
        \ul{\pow KX} := \up {p\in\mathcal P} \pow{\seq Kp}X.
        \end{equation*}
We start with an important fact about ultraproducts of powers series
rings taken from \cite[Lemma 3.4]{BDDL}; since we will need a similar
argument below (Proposition~\ref{P:infal}), we indicate the proof.
The \emph{ideal of infinitesimals} of a local ring $(S,\mathfrak m)$
is the ideal $\infal S := \bigcap_{d\in\mathbb N} {\mathfrak m}^d$ of
$S$. The $\maxim$-adic topology on $S$ is separated \iff\ $\infal
S=0$, and this is the case if $S$ is Noetherian by Krull's
Intersection Theorem.

\begin{proposition}\label{P:sep}
There is a surjective $\pol KX$-algebra \homo\ \[\pi\colon
\ul {\pow KX}\to \pow KX\] whose kernel is   $\infal{\ul {\pow KX}}$.
\end{proposition}
\begin{proof}
We start by defining $\pi$. Let $\ul f\in\ul {\pow KX}$ and choose
$\seq fp\in\pow{\seq Kp}X$,  for $p\in\mathcal P$,  whose
ultraproduct is $\ul f$. Write each $\seq fp$ as
        \begin{equation*}
        \seq fp:=\sum_\nu \seq{a_\nu}pX^\nu
        \end{equation*}
with $\seq{a_\nu}p\in\seq Kp$. Here the sum ranges over all
multi-indices $\nu\in {\mathbb N}^n$. Let $\ul
{a_\nu}\in K$ be the ultraproduct of the $\seq{a_\nu}p$ and define
        \begin{equation*}
        \pi(\ul f):=\sum_\nu \ul{a_\nu}X^\nu \in \pow KX.
        \end{equation*}
It follows from \los\ that $\pi$ is  a well-defined $\pol KX$-algebra
\homo. Its surjectivity is clear. So it remains to show that the
kernel of $\pi$ is  $\infal{\ul {\pow KX}}$. If $\ul f\in\infal{\ul
{\pow KX}}$, then by \los, for each $d\in\mathbb N$,
there is a member $U_d$ of
the ultrafilter such that $\seq fp\in\rij Xn^d\pow{\seq Kp}X$ for all
$p\in U_d$. In particular, for each $\nu\in {\mathbb N}^n$
we have that $\seq{a_\nu}p=0$, for all $p\in U_{|\nu|+1}$. Therefore
$\ul{a_\nu}=0$, and since this holds for all $\nu$, we see that
$\ul f\in\ker\pi$. The converse holds by reversing the
argument.
\end{proof}

\begin{remark}\label{R:sep}
In fact, we may replace in the above $\ul{\pow KX}$ by its subring
$\ul{\pol KX}$, given as the ultraproduct of the $\pol{\seq Kp}X$.
That is to say, $\pi$ induces a surjective $\pol KX$-algebra \homo\
$\ul{\pol KX}\to \pow KX$ with kernel equal to the intersection of
all $\rij Xn^d\ul{\pol KX}$ for $d\in\mathbb N$.
Indeed, it suffices to show that $\pi$
maps  $\ul{\pol KX}$ onto $\pow KX$. Let us explain this just in case
the underlying set of $\mathcal P$ is countable and hence, after
identification, we may think of it as a subset of $\nat$.
Given $f=\sum_\nu a_\nu X^\nu\in \pow KX$, choose
$\seq{a_\nu}p\in \seq Kp$ so that their ultraproduct is $a_\nu$ and
put
        \begin{equation*}
        \seq fp:=\sum_{|\nu|\leq p} \seq{a_\nu}p X^\nu \in K_p[X].
        \end{equation*}
Then $\pi(\ul f)=f$, where $\ul f\in\ul{\pol KX}$
is the
ultraproduct of the polynomials $\seq fp$.
\end{remark}

\subsection{Artin Approximation.}\label{R:AA}
Recall  that a Noetherian local ring $(R,\maxim)$ is said to satisfy 
\emph{Artin Approximation} if every system of polynomial equations
over $R$ which is solvable in the completion $\complet R$ of $R$ is 
already solvable in $R$. In view of Corollary~\ref{C:LP}, 
this is equivalent with the existence of an ultraset 
$\mathcal U$ and  an $R$-algebra \homo\
        \begin{equation}\label{eq:AAup}
        \complet R\to R^{\mathcal U}.
        \end{equation}
In fact, if $R$ satisfies Artin Approximation, then $R$ is existentially
closed in $\complet R$, that is to say,
every system of polynomial equations and inequalities
over $R$ which is solvable in $\complet R$ has a solution in $R$.
(Since $R$ is dense in $\complet R$, 
inequalities, and also congruence conditions,
can be incorporated in a system of equations.)
Artin proved (in \cite{Art68} and \cite[Theorem 1.10]{Art69},
respectively) that the ring of convergent complex power series 
$\mathbb C\{X\}$
and the ring of algebraic power series $\algpow LX$, 
with $L$ an arbitrary field, satisfy
Artin Approximation.

\subsubsection*{Artin's Conjecture.}
A local ring $(R,\frak m)$ satisfying Artin Approximation is necessarily Henselian,
and Artin conjectured that the converse holds if $R$ is excellent. 
This conjecture was eventually confirmed to be true
\cite{Pop,Spi,Swan}. In each of these papers, Artin's
Conjecture is derived from \emph{generalized N\'eron 
Desingularization}, stating that a \homo\ $A\to B$ of Noetherian 
rings is regular \iff\
$B$ is the direct limit of smooth $A$-algebras. In the development of
tight closure in \ch\ zero in the sense of Hochster and Huneke
\cite{HHZero}, this latter theorem plays an essential role. In 
this paper we give an
alternative definition of tight closure 
%(which conjecturally
%yields the same closure operation) 
relying only on a weaker form of 
Artin Approximation, to wit, Rotthaus' result \cite{Rot} on the 
Artin Approximation property for 
rings of the form $\algpow{\pow LX}Y$ with $L$ 
a field of \ch\ zero. (In Theorem~\ref{T:Z} below, which 
is not needed anywhere else, %except in Remark~\ref{R:E+pow}, 
we do  need generalized N\'eron Desingularization.)

\subsubsection*{Strong Artin Approximation.}
We say that a Noetherian local ring $(R,\maxim)$  satisfies \emph{Strong Artin
Approximation}, if any system of polynomial equations over $R$ which
is solvable modulo arbitrary high powers of $\maxim$ is already
solvable in  $R$. By Corollary~\ref{C:LP}, this  amounts to the
existence of an ultraset $\mathcal U$ and an $R$-algebra \homo\
        \begin{equation}\label{eq:SAAup}
        \prod_{n\in\nat} R/\maxim^n\to R^{\mathcal U}.
        \end{equation}
From \eqref{eq:AAup} and \eqref{eq:SAAup} it follows that $R$
satisfies Strong Artin Approximation \iff\ $R$ satisfies Artin
Approximation and $\complet R$ satisfies Strong Artin Approximation. In
\cite{BDDL}, a very quick proof using ultraproducts is given to show
that $\pow LX$ satisfies Strong Artin Approximation, for every
uncountable \acf\ $L$. Using the Cohen Structure Theorem, one then
deduces from this  and the positive solution of Artin's Conjecture,
that every equi\ch, excellent, Henselian local ring with an
uncountable algebraically closed residue field satisfies Strong
Artin Approximation.

\subsubsection*{Uniform Strong Artin Approximation.}
Any version in which the same conclusion as in Strong Artin
Approximation can be reached just from the solvability modulo a
single power $\maxim^N$ of $\maxim$, where $N$ only depends on (some
numerical invariants of) the system of equations, is called
\emph{Uniform Strong Artin Approximation}. In \cite{BDDL}, using
ultraproducts, Uniform Strong Artin Approximation for $R=\algpow LX$
is shown to follow from Artin Approximation for that ring. In 
more general situations,  
additional restrictions have to be imposed on the equations 
(see \cite [Theorem 6.1]{Art69} or \cite[Theorem 3.2]{BDDL}) and 
substantially more work is required \cite{DL80,SchPS}. For 
instance, the proof of the parametric version
in \cite[Theorem 3.1]{DSEx} uses the positive solution  \cite{Rot} of
Artin's Conjecture in the equi\ch\ case.

\subsubsection*{Nested Conditions.}
An even more subtle question regarding (Strong or Uniform Strong) Artin 
Approximation for subrings of $L[[X]]$
is whether one can maintain side  conditions 
on the solutions requiring some of the entries of  a 
solution tuple to depend only on some of the variables, provided  the 
given (approximate) solutions also satisfy such constraints. In 
\cite{BDDL}, several examples are presented to show that this might 
fail in general (see also \S\ref{s:uc} below). 
However, Rotthaus' approximation result \cite{Rot}
implies that  cylindrical approximation does hold,
provided $\operatorname{char} L=0$. (This was first noted in \cite{BDDL}.)
Here, by {\em cylindrical approximation}\/ we mean
Artin Approximation 
for nested systems of polynomial equations
in the context of Example~\ref{e:nested}.
We refer to Theorem~\ref{T:USAAfil} below for a  precise formulation.
%In \cite[Theorem~4.2]{BDDL} the authors prove
%Uniform  Strong Artin Approximation for $\algpow LX$ subject 
%to  a single  condition involving only a single variable. They 
%further observe \cite[Remark following Theorem~4.3]{BDDL} that this 
%can be extended to arbitrary (nested) conditions, once the 
%corresponding Artin Approximation result is proven to hold. For sake 
%of completeness, we have included their argument, see
%Theorem~\ref{T:USAAfil} below.

% (see for instance the argument in \cite[Remark 
%1.5]{BDL83}, which we will repeat here in the form of 
%Lemma~\ref{L:finch} below). In \cite[Theorem 4.2]{BDDL}, the authors 
%prove   

\subsection{Embeddings in ultraproducts}
We now turn to the issue of embedding a power series ring in 
the ultraproduct of power series rings, which is needed for our 
construction of a Lefschetz hull in the next section. The existence 
of a Lefschetz hull is immediate from the following corollary to 
Artin's  original result on the Artin Approximation property for 
algebraic power series. 

\begin{proposition}\label{P:LPK}
For every finitely generated $\pol KX$-subalgebra $V$ of $\pow KX$ 
there exists
a $\pol KX$-algebra \homo\  $V\to\ul{\pow KX}$. In particular, there 
exists an ultraset $\mathcal U$ and a  $\pol KX$-algebra \homo\ $\pow 
KX\to \ul{\pow KX}^{\mathcal U}$.
\end{proposition}
\begin{proof}
Translating the Artin Approximation property for $\algpow KX$ in the 
terminology of Corollary~\ref{C:LP} yields the existence of a $\pol 
KX$-algebra \homo\ $\varphi\colon V\to \algpow KX$.
As the Henselian property can be expressed in terms of the solvability of 
certain systems of polynomial equations, it
follows by \los\ that $\ul{\pow KX}$ is
Henselian. By the universal property of Henselizations 
there exists a unique $K[X]$-algebra homomorphism 
$\algpow KX\to \ul{\pow KX}$. Composition with $\varphi$
yields the desired $\pol 
KX$-algebra \homo\ $V\to\ul{\pow KX}$. The last 
assertion is now clear by Corollary~\ref{C:LP}.
\end{proof}

The remainder of the section is devoted to 
enhancements of this, and in particular, the \filtered\ version 
from Theorem~\ref{T:LP}, which we need to obtain  functoriality of
Lefschetz extensions. In the following $L$ denotes a field and
$S=\bigcup_n S_n$ the nested subring of
the nested ring $A=\bigcup_n A_n$ as
defined in Example~\ref{e:nested}, so
$$S_n = L[[Y_0]][Y_1,\dots,Y_n], \quad
A_n = L[[Y_0,\dots,Y_n]]\qquad\text{for all $n$,}$$ 
where $Y_0,Y_1,\dots$ is an infinite
sequence of finite tuples $Y_n$ of indeterminates.
If we need to emphasize the field, we will 
write $S_L$ and $A_L$ for $S$ and $A$.  
We need some further notations concerning nested rings.

\begin{notation}
Let
$Q=\bigcup_{n} Q_n$ be a \filtered\ ring. 
We denote by $Q(1)$ 
the ring $Q$ considered as a \filtered\ ring with nest $Q(1)_n:=Q_{n+1}$. 
A \homo\ $\psi\colon Q\to R$ of 
\filtered\ rings is then also a \homo\ 
$Q(1)\to R(1)$ of \filtered\ rings.

If $I$ is an ideal of $Q$, then we construe $Q/I$ as a nested ring with
nest given by $(Q/I)_n := Q_n/I\cap Q_n$ for all $n$. 
If $\pr$ is a prime ideal of $Q$, then the localization
$Q_{\pr}$ is a nested ring with nest given by
$(Q_{\pr})_n:=(Q_n)_{\pr\cap Q_n}$ for all $n$. 
If each $Q_n$ is a local ring with maximal ideal $\maxim_n$, then
$Q$ is local with maximal ideal $\maxim:=\bigcup_n\maxim_n$ and residue field
$Q/\maxim=\bigcup_n Q_n/\maxim_n$. In this case we say that $(Q,\maxim)$ is
a \emph{nested local ring.}\/
If moreover $Q_n\cap\maxim^k=\maxim_n^k$
for every $k$, then $\infal{Q}\cap Q_n=\infal{Q_n}$ for all $n$.

If $Q$ is a \filtered\ $R$-algebra,
for some \filtered\ ring $R=\bigcup_n R_n$, and $T$ an $R_0$-algebra, then we
consider $Q\tensor_{R_0} T$ as a \filtered\ $R$-algebra by means of the
nest $(Q_n\tensor_{R_0} T)$.

Given a Henselian local ring 
$(H,\mathfrak n)$ and a \homo\ $\psi\colon Q\to H$ 
we denote the Henselization of $Q_{\mathfrak n\cap 
Q}$ by $\hens Q$ and we let $\hens\psi\colon \hens Q\to H$ be the 
unique extension of $\psi$ given by the universal property of 
Henselizations. 
(Often, $H$ is to be understood from the context.) 
Note that then  $\hens Q$ is a \filtered\ local ring with
nest $(\hens Q)_n:=\hens{(Q_n)}$.
For instance, applied to $Q:=S$, $H:=A$ and $\psi$ the natural inclusion 
$S\to A$, we get the nested $S$-subalgebra 
$\hens S = R$ of $A$ (see Example~\ref{e:nested}).
\end{notation}

The argument in the proof of the following lemma was inspired by \cite[Remark 
1.5]{BDL83}.

\begin{lemma}[Cylindrical Approximation]\label{L:finch}
If $V$ is a \filtered\ $S$-subalgebra of finite type of $A$, then 
there exists a \homo\ of \filtered\  $S$-algebras $\varphi\colon V 
\to \hens S$.
\end{lemma}
\begin{proof}
We proceed by induction on the length $n$ of the relevant part 
$V_0\subseteq V_1\subseteq \dots \subseteq V_n$ of $V$, 
where the case $n=0$ holds 
trivially since  then $V=S$. 
Consider the nested ring $T:=S(1)\tensor_{S_1}A_1$ 
with nest $(T_n)$ given by
$$T_n:=S(1)_{n}\tensor_{S_1}A_1 = L[[Y_0,Y_1]][Y_2,\dots,Y_{n+1}].$$ 
In particular, $T$ is a \filtered\ 
$S(1)$-subalgebra of $A(1)$. 
Let $W$ be the 
image of the \homo\ of \filtered\ $T$-algebras 
$V(1)\tensor_{S_1}A_1\to A(1)$ induced by the inclusion $V(1)\subseteq A(1)$, 
so $W$ is a \filtered\ $T$-subalgebra of 
finite type of $A(1)$. Since its relevant part has length $<n$, we 
may apply our induction hypothesis to conclude that there is a \homo\ 
of \filtered\ $T$-algebras  $W\to \hens T$. Via the natural \homo\ 
$V(1)\to W$, we get a \homo\ of \filtered\ $S(1)$-algebras $V(1)\to 
\hens T$. Let $W'$ be its image, so that $W'$ is a   \filtered\ 
$S(1)$-subalgebra  of finite type of $\hens T$.

For our purposes in \S\ref{s:LP}, we only have to deal with the case that the base 
field  $L$ has \ch\ zero. In that case, we can use 
\cite[Theorem~4.2]{Rot}, which  implies that 
$\hens {S_1}$ has 
the Artin Approximation property. In case the \ch\ of $L$ is 
positive, we require the positive solution of Artin's Conjecture by 
\cite{Pop,Spi}.  In any case, by \eqref{eq:AAup},  there exists an 
ultraset $\mathcal U$ and an $S_1$-algebra \homo\ $A_1\to 
(\hens{S_1})^{\mathcal U}$. For each $n$,  this $S_1$-algebra homomorphism
extends to an 
$S_n$-algebra \homo\
        \begin{equation*}
        T_n=\pol{A_1}{Y_2,\dots,Y_{n+1}} \to 
\pol{(\hens{S_1})^{\mathcal U}} {Y_2,\dots,Y_{n+1}}\to 
(\hens{S_{n+1}})^{\mathcal U}.
        \end{equation*}
Since the right hand side is Henselian, we may replace the ring on 
the left by its Henselization. Gathering these \homo{s} for all $n$ 
yields a \homo\ of \filtered\ $S(1)$-algebras $\hens T\to 
(\hens{S(1)})^{\langle\mathcal U\rangle}$. Applying \eqref{i:lochomo} 
to the \filtered\ $S(1)$-subalgebra of finite type $W'\subseteq\hens 
T$, yields the existence of a \homo\ of \filtered\ $S(1)$-algebras 
$W'\to \hens{S(1)}$. Composition with $V(1)\to W'$ gives a \homo\ 
$V(1)\to \hens{S(1)}$ of \filtered\ $S(1)$-algebras. Since 
$S_0=V_0=A_0$, this is in fact a \homo\ of \filtered\ $S$-algebras 
$V\to \hens S$, as required.
\end{proof}

\emph{From now on we always assume that $Y_0$ is the empty tuple.}\/
(The more general case was only needed for inductive reasons, in
the proof of the previous lemma.)

In the following corollary we 
specialize to  $L=K$. Then, in a natural way,
$B_n:=\ul{\pow K{Y_1,\dots,Y_n}}$ is an 
$S_n$-algebra, and we may identify $B_n$ with a subring of
$B_{n+1}$, for all $n$; 
hence the subring $B:=\bigcup_n B_n$ of
$\up {p} \pow{\seq Kp}{Y_1,Y_2,\dots}$ is a
\filtered\ $S$-algebra with nest $(B_n)$.

\begin{corollary}\label{C:LPfil}
There exists an ultraset $\mathcal U$ and a \homo\  
$\eta\colon A\to B^{\langle\mathcal U\rangle}$ of \filtered\ 
$S$-algebras.
\end{corollary}
\begin{proof}
We only need to verify that condition~\eqref{i:lochomo} in 
Theorem~\ref{T:LPfil} is fulfilled. To this end, let $V$ be a 
\filtered\ $S$-subalgebra of $A$ of finite type. By 
Lemma~\ref{L:finch}, there exists a \homo\ of \filtered\ $S$-algebras 
$V\to \hens S$. Since $B$ is Henselian, the \homo\ of \filtered\ 
rings   $S\to B$ extends to a \homo\ of \filtered\ rings   $\hens 
S\to B$, and the composition $V\to \hens S\to B$ proves 
\eqref{i:lochomo}.
\end{proof}

We denote by $\maxim$ the ideal of $S$ generated by all the indeterminates
$Y_{ni}$, for all $n$ and $i=0,\dots, k_n$.

\begin{remark}\label{R:piLPfil}
For each $n$, let $\pi_n$ be the canonical epimorphism $B_n\to A_n$ 
given by Proposition~\ref{P:sep}  and let $\pi\colon B\to A$ be the 
induced \homo\ of \filtered\ $S$-algebras (given as the direct limit 
of the $\pi_n$). Then $\pi$ induces an isomorphism between 
$B/\maxim^cB$ and  $S/\maxim^cS$, for all $c\in\nat$. On the other hand, for 
a fixed $c\in\nat$, we can realize $A$ as the union of all \filtered\ 
$S$-subalgebras $V$ of finite type of $A$ such that $V/\maxim^c V\iso 
S/\maxim^c S$. For those $V$, the \homo\ $V\to \hens S$ given by 
Lemma~\ref{L:finch} becomes an isomorphism modulo $\maxim^c$, and
applying Remark~\ref{R:cofinal} with this collection 
we see that we may take $\eta$ in 
Corollary~\ref{C:LPfil} so that its composition with $\pi$ is 
congruent modulo $\maxim^c$ to the diagonal embedding $A\subseteq 
A^{\langle\mathcal U\rangle}$.
Without proof, we mention that one can achieve the even stronger 
condition that $\pi\circ\eta$ is \emph{equal} to the diagonal embedding. (This 
however, even in \ch\ zero, requires  generalized N\'eron 
Desingularization.)
\end{remark}

Applying Corollary~\ref{C:LPfil} with each $Y_n$ for $n\geq 1$ equal 
to a single indeterminate yields 
the following result, needed for the functorial construction of 
faithfully flat Lefschetz extensions in the next section.

\begin{theorem}\label{T:LP}
There exists an ultraset $\mathcal U$ and for each $n$  a
$\pol K{X_1,\dots,X_n}$-algebra \homo\
        \begin{equation*}
        \eta_n\colon \pow K {X_1,\dots,X_n} \to \ul{\pow K
{X_1,\dots,X_n}}^{\mathcal U},
        \end{equation*}
such that for all $n\leq m$, the diagram
$$\xymatrix@R+2em@C+3em{
{\pow K {X_1,\dots,X_n}} \ar[r]^{\eta_n} \ar[d] & 
{\ul{\pow K {X_1,\dots,X_n}}^{\mathcal U}} \ar[d] \\
{\pow K {X_1,\dots,X_m}} \ar[r]^{\eta_m} & 
{\ul{\pow K {X_1,\dots,X_m}}^{\mathcal U}}
}$$
commutes, where the vertical arrows are the natural inclusion maps. \qed
\end{theorem}

Given a nested system of polynomial equations
$\mathcal S$ over $S$ as in \eqref{e:nestedeq} we call the maximum of 
$n,k_0,\dots,k_n$ and
the degrees of the polynomials $P_{ij}$ the \emph{complexity} of $\mathcal S$.
We say that a \filtered\ $S$-algebra $V$ of finite type has
\emph{complexity $\leq d$} if $V$ admits a defining system of
\filtered\ equations of complexity $\leq d$. (See \S\ref{s:ft}.)
The proof of the next theorem is a modification of the argument in 
\cite[Theorem~4.3]{BDDL}.

\begin{theorem}[Uniform Strong Artin Approximation with Nested
Conditions]\label{T:USAAfil}
Given $c,d\in\mathbb N$, there exists $N=N(c,d)\in\mathbb N$ 
with the following property. Let $L$ be a field, let $S:=S_L$  and 
let $V$ be a \filtered\ $S$-algebra of finite type and of complexity
at most $d$. 
If $\psi\colon V\to S/\maxim^NS$ is a \homo\ of $S$-\filtered\ 
algebras, then  there exists a \homo\ $\varphi\colon V\to \hens S$ of 
\filtered\ $S$-algebras  such that
$$\xymatrix@R+2em@C+3em{V \ar[r]^{\varphi} \ar[d]^{\psi} & \hens S \ar[d]^q \\
S/\maxim^NS \ar[r] & S/\maxim^cS}
%        \commdiagramunnumbered V{\varphi}{\hens S} {\psi}{q} 
%{S/\maxim^NS} {} {S/\maxim^cS}
$$
commutes, where $q$ is induced by the canonical 
isomorphism $\hens S/\maxim^c\hens S\iso S/\maxim^cS$.
\end{theorem}
\begin{proof}
Suppose the claim is false for some pair $(c,d)$, so that we 
have counterexamples for increasing powers of $\maxim$. That is to 
say, for each $p\in\nat$ there  is a field $\seq Kp$ and  a 
\filtered\ $S_{\seq Kp}$-algebra $\seq Vp$ of finite type with 
a defining \filtered\ equations 
$\seq{\mathcal S}p$ of complexity at most $d$ and 
a \homo\ of \filtered\ algebras $\seq Vp\to S_{\seq 
Kp}/\maxim^pS_{\seq Kp}$ which is not congruent modulo $\maxim^c$ to
a \homo\ $\seq Vp\to\hens {S_{\seq Kp}}$.
We may assume that there exist
$k\in\nat$ and indeterminates
$Z_i=(Z_{i1},\dots,Z_{id})$, for $i=0,\dots,d$, such that each
$\mathcal S_p$ has the form \eqref{e:nestedeq} with $n=d$, for some
polynomials $P_{ij}\in (S_{K_p})_i[Z_0,\dots,Z_i]$ of degree $\leq d$.
Let $K$ be the ultraproduct of the $\seq Kp$ with respect to some
ultraset with underlying set $\nat$.   
Taking the ultraproduct of the polynomials in 
$\seq{\mathcal S}p$ yields  a \filtered\ system $\mathcal S$ 
of equations with coeffcients in $S_K$. 
Let $V$ be the corresponding \filtered\ 
$S_K$-algebra of finite type. By \cite[(1.8)]{SvdD}, $V$ embeds into
the ultraproduct of the $V_p$.  Taking ultraproducts 
of the \homo{s} $\seq Vp\to S_{\seq Kp}/\maxim^p S_{\seq Kp}\iso 
A_{\seq Kp}/\maxim^p A_{\seq Kp}$ yields a \homo\ $V\to B_K/\infal 
{B_K}\iso A_K$ of \filtered\ $S_K$-algebras, where we used 
Proposition~\ref{P:sep} for the last isomorphism. By 
Lemma~\ref{L:finch}, applied to (the image under) $V\to A_K$, there 
exists a \homo\ $V\to \hens {S_K}$ of \filtered\ $S_K$-algebras, which 
we may assume to be congruent to  $V\to A_K$ modulo $\maxim^cA_K$, by
Remark~\ref{R:piLPfil}. By \los, the 
ultraproduct $\tilde B$ of the $\hens{S_{\seq Kp}}$ is Henselian. 
Since $\tilde B$ is a \filtered\ $S_K$-algebra, it is in fact a 
\filtered\ $\hens {S_K}$-algebra by the universal property of 
Henselizations. Hence we have a composed \homo\ $V\to \tilde B$ of 
\filtered\ $S_K$-algebras which is congruent to $V\to A_K$ modulo 
$\maxim^cA_K$. \los\ then yields for almost all $p$ a \homo\ $\seq 
Vp\to\hens{S_{\seq Kp}}$ of \filtered\ $S_{\seq Kp}$-algebras which 
modulo $\maxim^c$ is equal to the original \homo\ $\seq Vp\to S_{\seq 
Kp}/\maxim^cS_{\seq Kp}$, a contradiction.
\end{proof}

\begin{remark}
Conversely, Lemma~\ref{L:finch} is an immediate consequence of 
Theorem~\ref{T:USAAfil}. Indeed, let  $V$ be a \filtered\ 
$S$-subalgebra of $A$   of finite type. Since $A/\maxim^NA\iso 
S/\maxim^NS$, this induces for each $N$ a \homo\ $V\to S/\maxim^NS$ 
of \filtered\ $S$-algebras. For sufficiently large $N$ this yields by 
Theorem~\ref{T:USAAfil} a \homo\ $V\to \hens S$ of \filtered\ 
$S$-algebras.
\end{remark}

\begin{remark}
Spelling out the previous result in terms of equations 
%when 
%each $Y_n$ for $n\geq 1$ is equal to a single 
%indeterminate 
yields the following equational form of cylindrical 
approximation: For all $c,d\in\nat$ there exists  a 
bound $N=N(c,d)\in\mathbb N$ with the  following property. Let
$L$ be a field and
let $\mathcal S$ be a nested system of
polynomial equations with coefficients from $S_L$, of complexity at most $d$.
If $\mathcal S$
has an approximate nested solution $\mathbf{a}=
(\mathbf{a}_0,\dots,\mathbf{a}_n)$ in $A_L$
modulo $(Y_1,\dots,Y_n)^N$, then $\mathcal S$ has a nested
solution in $\hens{S_L}$ 
which is congruent to $\mathbf{a}$ modulo $(Y_1,\dots,Y_n)^c$.
%and we are given
%polynomials $$P_1(X,Y),\dots,P_s(X,Y)\in\pol 
%L{X,Y}$$ of degree at most $d$, where
%$X=\rij Xn$ and $Y=\rij Ym$, a partition $0\leq 
%i_1<\dots<i_{n-1}<i_n=m$   and an $m$-tuple of polynomials 
%$\mathbf f=\rij fm$ with $f_1,\dots,f_{i_j}\in\pol L{X_1,\dots,X_j}$ 
%for all $j=1,\dots, n$ and  $$P_1(X, \mathbf f)\equiv \cdots\equiv
%P_s(X,\mathbf f)\equiv0\mod\rij Xn^N.$$ 
%Then we can find  $\mathbf g=\rij gm$ with $g_1,\dots,g_{i_j}\in\algpow 
%L{X_1,\dots,X_j}$ for all $j=1,\dots,n$ such that $$P_1(X, \mathbf g)=\cdots=
%P_s(X, \mathbf g)=0$$ and
%$\mathbf f\equiv\mathbf g\bmod\rij Xn^c$.
\end{remark}

\begin{remark}\label{r:e-transfer}
Let $\mathcal{L}(n)$ be the language of rings
augmented by unary predicate symbols $R_0,\dots,R_n$. We construe a 
formal power series ring $S[[X_1,\dots,X_n]]$ over a ring $S$ as
an $\mathcal{L}(n)$-structure by interpreting $R_i$ by the subring
$S[[X_1,\dots,X_i]]$. The previous remark yields the 
following \emph{existential Lefschetz principle}
for nested power series rings: 
An existential $\mathcal{L}(n)$-sentence holds in
$\mathbb C[[X_1,\dots,X_n]]$ \iff\ it holds in $\ac Fp[[X_1,\dots,X_n]]$ 
for all but finitely many
primes $p$.  For existential sentences not involving
the $R_i$, this has  already been noted elsewhere, see
\cite[Proposition~1]{Becker-Lipshitz}.
For $n=1$ a much stronger transfer principle 
holds  in which not only  existential sentences are carried 
over, but any sentence. (This follows from the
Ax-Kochen-Ershov Principle.)
\end{remark}

We finish this section by indicating a strengthening of
Theorem~\ref{T:LP} (not needed later).
%(only needed in \S\ref{Model-theoretic Applications}). 
Given a power series $f\in\pow \zet X=\zet[[X_1,\dots,X_n]]$, let $\seq
fp$ be its image in $\pow{\seq Kp}X$ and let $\ul f$ be the
ultraproduct of the $\seq fp$ in $\ul {\pow KX}$.  One verifies that
the map $f\mapsto \ul f$ is an injective $\pol{\zet}X$-algebra \homo\
$\pow\zet X\to \ul {\pow KX}$ which extends to an injective
$\pol{\zet}X$-algebra \homo\ $\pow\zet X \tensor_{\zet} K\to \ul
{\pow KX}$. We will view $Z:=\pow\zet X \tensor_{\zet} K$ as a
subring of $\ul {\pow KX}$ via this embedding. Write $\maxim_n$ for the
ideal in $\pol{\zet}X$ generated by $X_1,\dots,X_n$. Since
$Z/\maxim_n^kZ\iso \pow KX/\maxim_n^k \pow KX$, for all $k$, we see that
$\pow KX$ is the $\maxim_n Z$-adic completion of $Z$. In particular,
$Z$ is a dense subring of $\pow KX$, equal to the $\pol
KX$-subalgebra of $\pow KX$ generated by all power series with
integral coefficients.  Inspecting the proof of
Proposition~\ref{P:sep}, we see that $\pi$ is in fact a $Z$-algebra
\homo. Let $\hens Z$ be the Henselization of $Z$ at the maximal ideal
$\maxim_n Z$.  By the universal
property of Henselizations, the embedding $Z\subseteq \ul {\pow KX}$
extends to a unique embedding $\hens Z\to\ul {\pow KX}$.
Henceforth we will view $\hens Z$ as a subring of $\ul {\pow KX}$.
Note that since $\pol KX$ is a subring of $\hens Z$, so is its
Henselization $\algpow KX$ at $\maxim_nK[X]$.

\begin{theorem}\label{T:Z}
The ultraset $\mathcal U$ and the $K[X_1,\dots,X_n]$-algebra homomorphisms
$\eta_n$ in Theorem~\textup{\ref{T:LP}} can be chosen so that in addition 
each $\eta_n$ is a
$\hens Z$-algebra \homo.
\end{theorem}
\begin{proof}
This theorem follows as above from the corresponding extension of
Lemma~\ref{L:finch}. To this end, replace the filtered ring $S$ in
Lemma~\ref{L:finch} and its proof
by the \filtered\ ring $T=\bigcup_nT_n$,  where $T_n$ is
the Henselization of  $\pow\zet {Y_0,\dots,Y_n} \tensor_{\zet} K$
at $(Y_0,\dots,Y_n)Z$.
Whenever we invoked Rotthaus' result, we now use the positive 
solution of Artin's Conjecture due to \cite{Pop,Spi,Swan}
instead.  Note that  each  $T_n$ is excellent. (Use for instance the
Jacobian Criterion \cite[Theorem 101]{Ma70}.) Details are left to the reader.
\end{proof}

\begin{remark}
By the same argument as in Remark~\ref{R:piLPfil}, we  can choose
the $\eta_n$ moreover so that its composition with the canonical
epimorphism $\pi^{\mathcal U}\colon\ul {\pow KX}^{\mathcal U}\to
\pow KX^{\mathcal U}$ is equal to the diagonal embedding $\pow
KX\subseteq \pow KX^{\mathcal U}$, for each $n$.
\end{remark}

\begin{remark}
By Theorem~\ref{T:Z}, the existential Lefschetz principle
from Remark~\ref{r:e-transfer}
remains true upon augmenting $\mathcal{L}(n)$ 
by additional constant symbols, one for each power series
in $\pow{\zet}X=\pow{\zet}{X_1,\dots,X_n}$, to be interpreted in
the natural way in $S[[X]]$.
\end{remark}

\section{Lefschetz Hulls}\label{s:LP}

Our objective in this section is to prove the theorem stated in the 
introduction, in a more precise form. {\it Throughout, we fix a Lefschetz 
field $K$ with respect to some ultraset with underlying set equal to 
the set of the prime numbers, whose components $K_p$ are \acf{s} of 
\ch\ $p$.}\/ See the remark following Proposition~\ref{P:acf}  on how to 
obtain such $K$, of arbitrarily large cardinality.

In obtaining a functorially defined Lefschetz extension,
we face the following complication: not every 
automorphism of $K$ is an ultraproduct of automorphisms of its 
components $\seq Kp$. The simplest counterexample is complex 
conjugation on $\mathbb C$, for no \acf\ of positive \ch\ has a 
subfield of index $2$. 
In fact, each subfield of $K$ has an automorphism which cannot be extended to
an automorphism of $K$ that is an ultraproduct of automorphisms of
the $K_p$. 
Therefore there cannot exist a functor 
from the category of equi\ch\ zero Noetherian local 
rings $R$ whose residue field is contained in $K$ to a category of 
analytic Lefschetz rings. 
%since certain automorphisms of $R$ 
%(obtained 
%for instance by a twist on the base field) 
%cannot extend to an 
%automorphism  of $\hull R$. 
The way around this problem is to fix 
some additional data of $R$, as we will now explain.

\subsection{Quasi-coefficient fields}\label{s:qcf}
Let $(R,\maxim)$ be a Noetherian local ring which 
contains the rationals (that is to say, $R$ has equi\ch\ zero).
A subfield $k$ of $R$ is called a \emph{quasi-coefficient field} of $R$ 
if $R/\maxim$ is algebraic over the image of $k$ under the residue
homomorphism $R\to R/\maxim$. Every maximal 
subfield of $R$ is a quasi-coeffcient field. A quasi-coefficient 
field is called a \emph{coefficient field} if the natural map $k\to 
R/\maxim$ is an isomorphism. In general, coefficient fields may not exist.
If $R$ is Henselian  then a subfield of $R$ is a coefficient field \iff\ it is
maximal. In particular, if $R$ is complete, then $R$ has a coefficient field. 
Every  
quasi-coefficient field $k$ of $R$ is contained in
a unique coefficient field of $\complet R$, namely,
the algebraic closure   of $k$ in 
$\complet R$. For proofs and more details, see \cite[\S28]{Mats}.

\subsection{The category $\cohcat$}
In order to state a refined version of the theorem from the 
introduction, we introduce a category $\cohcat$ (for ``Cohen''). 
Its objects are 
quadruples $\Lambda=(R,\mathbf x,k,u)$ where 
\begin{itemize}
\item[(a)] $(R,\maxim)$ is a
Noetherian local ring (the \emph{underlying ring} of $\Lambda$), 
\item[(b)] $\mathbf x$ is a (finite) tuple of elements of $R$ which generate 
$\maxim$, 
\item[(c)] $k$ is a quasi-coefficient field of $R$, and 
\item[(d)] $u\colon R\to K$ is a 
local \homo\ (that is to say, $u$ is a ring homomorphism with 
$\ker u=\maxim$).
\end{itemize} 
A morphism $\Lambda\to\Gamma$ from $\Lambda$ to another such quadruple 
$\Gamma=(S,\mathbf y,l,v)$ is given by a local ring \homo\ 
$\alpha\colon R\to S$ such that
\begin{itemize}
\item[(a)] $\alpha(\mathbf x)$ is an 
initial segment of $\mathbf y$ (if $\mathbf x=\rij 
xn$ and $\mathbf y=\rij ym$, then $n\leq m$ and $y_i=\alpha(x_i)$ for 
$i=\range 1n$),
\item[(b)] $\alpha(k)\subseteq l$, and
\item[(c)] $v\after\alpha=u$.
\end{itemize}
We will often identify a morphism $\Lambda\to\Gamma$ with its 
\emph{underlying \homo}  $\alpha\colon
R\to S$ and hence denote it also by $\alpha$.

Let $\loccat$ be the category of (not necessarily Noetherian) local 
rings, with  the local ring \homo{s} as morphisms. 
Given an ultraset $\mathcal W$ we denote by 
$\lefcat_{\mathcal W}$ the 
category of analytic Lefschetz rings
with respect to $\mathcal W$ as defined in \S\ref{s:UP}.
(Its objects 
are ultraproducts, with respect to $\mathcal W$, 
of complete local rings with algebraically closed 
residue fields of positive \ch, and its morphisms are ultraproducts 
of local \homo{s}.) We stress once more that $\lefcat_{\mathcal W}$, as a 
subcategory of $\loccat$, is not   full.  
%All analytic Lefschetz 
%rings will actually be $K$-algebras and the morphisms between them 
%will be $K$-algebra \homo{s}, so that we could have confined 
%ourselves to this smaller category.
We will denote the \emph{forgetful functor} with values in $\loccat 
$ always  by $\fff$ (regardless of the source category). If 
$F$ and $G$ are functors from a category $\op{\textbf {C}}$ to 
$\loccat $, then we will say that a natural transformation 
$\eta\colon F\to G$ is \emph{faithfully flat} if the ring \homo\ 
$\eta_\Lambda \colon F(\Lambda)\to G(\Lambda)$ is faithfully flat, 
for each object $\Lambda$ in $\op{\textbf {C}}$.

\begin{theorem}\label{T:maincan}
There exists an ultraset $\mathcal W$,
a functor $\mathfrak D\colon \cohcat \to \lefcat_{\mathcal W}$ and a 
faithfully flat natural transformation $\eta\colon \fff\to 
\fff\after\mathfrak D$.
\end{theorem}

We call $\hull\Lambda$ the \emph{Lefschetz hull} of $\Lambda$.
Let us state in more detail what the above functoriality amounts to. 
Given a morphism $\Lambda \to \Gamma $ in $\cohcat$ 
with underlying \homo\ $\alpha\colon R \to S$, where $R:=\fff(\Lambda)$ 
and $S:=\fff(\Gamma)$, we get a morphism $\hull \alpha\colon \hull 
\Lambda\to \hull\Gamma$ in $\lefcat_{\mathcal W}$ and   
faithfully flat \homo{s} 
$\eta_\Lambda\colon R\to \hull\Lambda $ and $\eta_\Gamma\colon S\to 
\hull\Gamma$  fitting into a commutative 
diagram
        \begin{equation}\label{fun}
        \begin{aligned}
        \mbox{
        \xymatrix@R+2em@C+3em{R \ar[r]^\alpha \ar[d]^{\eta_\Lambda} & 
S       \ar[d]^{\eta_\Gamma} \\
        \hull{\Lambda} \ar[r]^{\hull \alpha} &\hull{\Gamma}.}
        }
        \end{aligned}
        \end{equation}
(Technically speaking we should have written $R\to 
\fff(\hull\Lambda)$, etc., but  we'll always identify 
$\hull\Lambda$ with its underlying ring.)

The proof proceeds in two steps. We first prove the theorem for a 
certain subcategory $\ancat$ of $\cohcat$ given by quotients of 
power series rings over $K$ 
(see \S\ref{s:pow}). The existence of the functor $\mathfrak D$ 
for these rings then follows from Theorem~\ref{T:LP}. The second step 
consists in associating in a functorial way to an object 
$\Lambda=(R,\mathbf x,k,u)$ of $\cohcat $ a complete local $K$-algebra 
which is a faithfully flat $R$-algebra (see 
\S\ref{s:proof}).  This is achieved by making a base change to $K$ 
using   $k$ and $u$, and then  taking completion. By Cohen's 
structure theorem, $\mathbf x$ then determines a unique ring
$C(\Lambda)$ in $\ancat$ isomorphic to the latter.
After the proof of Theorem~\ref{T:maincan} we discuss a construction of
Lefschetz hulls with some additional properties.
We finish the section by pointing out (in \S\ref{s:uc})
another obstacle which prevented us from outright defining a functor
from the category of Noetherian local rings whose residue field embeds into $K$
to a category of analytic Lefschetz rings.

We adopt the following notation for  polynomial and power series rings:
we fix a  countable sequence of indeterminates $X_1,X_2,\dots$, and for each
$n$ and each ring $S$, we let $\pol Sn$ and  $\pow Sn$ be shorthand 
for respectively  $\pol S{X_1,\dots,X_n}$ and $\pow 
S{X_1,\dots,X_n}$.  We  write $\ul{\pow Kn}$ for the 
ultraproduct of the $\pow {\seq Kp}n$.

\subsection{Power series rings}\label{s:pow}
We first describe in more detail the category of quotients of power 
series rings over $K$, which we denote by  $\ancat$. Its 
objects are local rings of the form $\pow Kn/I$, for some $n$ and 
some ideal $I$ of $\pow Kn$. 
A morphism in $\ancat$ is a 
$K$-algebra \homo\ $\alpha\colon\pow Kn/I\to \pow Km/J$ where
$n\leq m$, $I\subseteq J$, and $\alpha$
is induced by the inclusion $\pow 
Kn\subseteq\pow Km$. To each object 
$\pow Kn/I$ of $\ancat $ we  associate the
object  $(\pow Kn/I,\mathbf{x},K,\pi_n)$ in $\cohcat$,  
where $\mathbf{x}=(x_1,\dots,x_n)$ with $x_i:=X_i+I$ for each $i$
and $\pi_n\colon \pow Kn\to K$ is the residue map.
Every $\ancat$-morphism  $\alpha\colon\pow Kn/I\to \pow Km/J$
gives rise to a $\cohcat$-morphism (with underlying homomorphism $\alpha$)
between the objects
corresponding to $\pow Kn/I$ and $\pow Km/J$, respectively.
It is easily verified that via this identification, $\ancat$ becomes 
a full subcategory of $\cohcat$. 
%Indeed, if $\alpha\colon \pow Kn/I\to \pow Km/J$ is the underlying
%homomorphism of a $\cohcat$-morphism between the corresponding 
%quadruples, then by 
%definition $\alpha(X_i)=X_i$ for $i=\range 1n$ and 
%$\pi_m\circ\alpha=\pi_n$. The latter equality implies in particular that 
%$\alpha$ leaves $K$ invariant. These two facts combined show that 
%$\alpha$ is a morphism in  $\ancat$. 
%Henceforth, we will always view 
%$\ancat$ as a subcategory of $\cohcat $.

We now embark on the proof of Theorem~\ref{T:maincan}, first for the subcategory 
$\ancat $. Let $\mathcal U$ be the ultraset proclaimed in
Theorem~\ref{T:LP} and set
         \begin{equation*}
         \hull n:=\ul{\pow Kn}^{\mathcal U} \qquad\text{for each $n$.}
         \end{equation*}
By that theorem, there exists, for each $n$, a $\pol Kn$-algebra
\homo\ $\eta_n\colon \pow Kn\to \hull n$ such that for each $n\leq
m$, the diagram
        \begin{equation}\label{mnhull}
        \begin{aligned}
        \mbox{\xymatrix@R+2em@C+3em{{\pow Kn} \ar[r] \ar[d]^{\eta_n} 
& {\pow Km }    \ar[d]^{\eta_m} \\
        {\hull n}  \ar[r]  & {\hull m}}}
        \end{aligned}
        \end{equation}
commutes, where the horizontal maps are the natural inclusions.
We construe $\hull n$ as a $K$-algebra via $\eta_0$;
then each $\eta_n$ is a $K$-algebra homomorphism.

\begin{remark}
If we  are only interested in constructing a Lefschetz extension for a single
$\pow Kn$, then the existence of a $\pol Kn$-algebra \homo\
$\eta_n\colon\pow Kn\to \hull n$ already follows by combining
Theorem~\ref{T:LPfil} with the more elementary Proposition~\ref{P:LPK}.
\end{remark}

\begin{remark}
Suppose that $K=\mathbb C$.
If we are willing to weaken the requirement that $\eta_n$ be
a $K[n]$-algebra homomorphism, then under assumption of the
Continuum Hypothesis $2^{\aleph_0}=\aleph_1$
the passage to the ultrapower $\mathbb C[[n]]_\infty^{\mathcal U}$ is
superfluous: Let $L$ be a countable subfield of $\mathbb C$; then
$S_n=L[n]$ is countable,
and under the assumption $2^{\aleph_0}=\aleph_1$ it follows along the lines of
Remark~\ref{R:sat} that there exists, for each $n$, an
$S_n$-algebra homomorphism $\varrho_n\colon \mathbb C[[n]]\to
\mathbb C[[n]]_\infty$ such that $\varrho_n$ is the restriction of 
$\varrho_m$ to $\mathbb C[[n]]$, for all $n\leq m$.
\end{remark}

Note that $\hull n$, being an ultrapower of the analytic Lefschetz ring
$\ul{\pow Kn}$, is itself an analytic Lefschetz ring. Indeed, we can construct
an ultraset $\mathcal W$ with the following property: for each $n$,
the ring $\hull n$ is isomorphic to the ultraproduct with respect to
$\mathcal W$ of the rings
$\pow{\seq Kw}n$, where $\seq Kw:=\seq K{p(w)}$ for some prime number $p(w)$.
(See \cite[Proposition~6.5.2]{Chang-Keisler}.)
 From now on, we always represent $\hull n$
in this way. The Lefschetz ring $\hull n$ is a $\mathbf K$-algebra,
where $\mathbf K:=\hull{0}=K^{\mathcal U}$, via the natural inclusion
$\hull{0}\to\hull{n}$, and the natural inclusions $\hull{n}\to\hull{m}$
(for $n\leq m$) are $\mathbf{K}$-algebra homomorphisms.
Next we show that $\hull n$  gives the desired faithfully flat 
Lefschetz extension:

\begin{proposition}\label{P:ff}
For each $n$, the \homo\
         $\eta_n\colon \pow Kn\to \hull n$
is faithfully flat.
\end{proposition}

In the proof we use the following variant of
\cite[Corollary 8.5.3]{Bruns-Herzog}.
A module $M$ over a local ring $R$ is called a \emph{big \CM\ module}
over $R$ if there exists a
system of parameters of $R$ which is an $M$-regular sequence.
If every system of parameters of $R$ is an $M$-regular sequence, then $M$
is called a \emph{balanced big \CM\ module} over $R$.
If $(R,\mathfrak m)$ is a regular local ring and $M$ a
balanced big \CM\ module over $R$, then $M$ is
flat, see \cite[proof of Theorem~9.1]{HuTC}.
%\cite[Theorem
%IV.1]{SchFPD} or
%\cite[Lemma 2.1 (d)]{HHbigCM2}.

\begin{lemma} \label{l:balanced}
Let $R$ be a Noetherian local ring and let $M$ be a big
\CM\ module over $R$. If every permutation of an $M$-regular sequence
is again $M$-regular, then $M$ is a balanced big \CM\ module over $R$.
\end{lemma}
\begin{proof}
We proceed by induction on $d:=\dim R$. There is nothing to show if
$d=0$, so assume $d>0$. By assumption, there exists a system of parameters
$(x_1,\dots,x_d)$ of $R$ which is an $M$-regular sequence. Let
$(y_1,\dots,y_d)$ be an arbitrary system of parameters. By prime avoidance we
find $z\in\mathfrak m$ not contained in a minimal prime of 
$(x_1,\dots,x_{d-1})R$
and of $(y_1,\dots,y_{d-1})R$. Hence both $(x_1,\dots,x_{d-1},z)$ and
$(y_1,\dots,y_{d-1},z)$ are systems of parameters of $R$.
Since a power of $x_n$ is a multiple of $z$ modulo $(x_1,\dots,x_{d-1})R$,
the sequence $(x_1,\dots,x_{d-1},z)$ is $M$-regular. Thus, by 
assumption, the
permuted sequence $(z,x_1,\dots,x_{d-1})$ is also $M$-regular. In particular,
the canonical
image of $(x_1,\dots,x_{d-1})$ in $R/zR$ is $M/zM$-regular, showing
that $M/zM$ is a big \CM\ module over $R/zR$. Moreover, every 
permutation of an $M/zM$-regular sequence is again $M/zM$-regular. By 
induction hypothesis, the
canonical image of $(y_1,\dots,y_{d-1})$ in $R/zR$, being a system of
parameters in $R/zR$, is $M/zM$-regular. Hence $(z,y_1,\dots,y_{d-1})$ is
$M$-regular, and therefore, using the assumption once more, so is 
$(y_1,\dots,y_{d-1},z)$. As some power of
$z$ is a multiple of $y_d$ modulo $(y_1,\dots,y_{d-1})R$, we get that
$(y_1,\dots,y_d)$ is $M$-regular, as required.
\end{proof}

\begin{proof}[Proof of Proposition~\ref{P:ff}]
Since $\rij Xn$ is a $\pow{\seq Kw}n$-regular sequence for each
$w$, it is a $\hull n$-regular sequence by \los. It
follows that $\hull n$ is a big \CM\ $\pow Kn$-algebra
via the \homo\ $\eta_n$. Using \los\ once
more, one  shows that every permutation of a $\hull n$-regular
sequence is again $\hull n$-regular (since every permutation of a
$\pow{\seq Kw}n$-regular sequence in $\pow{\seq Kw}n$ remains
$\pow{\seq Kw}n$-regular
by \cite[Proposition 1.1.6]{Bruns-Herzog}). Therefore,
$\hull n$ is a balanced big \CM\ $\pow Kn$-algebra, by the lemma above.
Since $\pow Kn$ is regular,
$\eta_n$ is flat by \cite[proof of Theorem~9.1]{HuTC}, hence faithfully flat.
\end{proof}

Below, we write $I\hull n$ to denote the ideal of $\hull n$ generated
by the image of an ideal $I$ of $\pow Kn$ under $\eta_n$.

\begin{remark}
We have
        \begin{equation*}
        \op{Im}(\eta_n) = \op{Im}(\eta_{n+1}) \cap \hull n\qquad
\text{for all $n$}.
        \end{equation*}
This follows from $X_{n+1}\hull{n+1}\cap\hull{n}=(0)$ and the injectivity
of $\eta_{n+1}$.
\end{remark}

\subsection{Proof of Theorem~\ref{T:maincan} for the category $\ancat 
$}\label{s:quotients}
The construction of $\hull n$ above extends in a natural way to quotients
of $\pow Kn$. Namely, if $I=(a_1,\dots,a_m)\pow Kn$
is an ideal of $\pow Kn$ and $R:=\pow Kn/I$, then we
choose $\seq {b_i}w\in \pow{\seq Kw}n$ whose ultraproduct in $\hull n$
is $\eta_n(a_i)$, for each $i$, and put
        \begin{equation*}
        \hull R := \up w \pow{\seq Kw}n/\seq Iw.
        \end{equation*}
Here $\seq Iw$ is the ideal of $\pow{\seq Kw}n$ generated by 
$\seq{b_1}w,\dots,\seq{b_m}w$. The canonical surjections $\pow{\seq 
Kw}n\to \pow{\seq Kw}n/\seq Iw$ yield a surjection
$\hull n \to\hull R$ whose kernel is $I\hull n$.
On the one hand, this  shows that $\hull R$ does not depend on the 
choice of the $\seq{b_i}w$ and that we have an isomorphism 
$\varphi\colon\hull n/I\hull n \to\hull R$. On the other hand, 
composing with the homomorphism $\eta_n\colon \pow Kn\to\hull n$ we 
obtain a homomorphism $\pow Kn\to\hull R$ whose kernel contains $I$, 
and hence an induced $K$-algebra
homomorphism
        \begin{equation*}
        \eta_R\colon R=\pow Kn/I\to\hull R.
        \end{equation*}
(According to this definition $\hull {\pow Kn}=\hull n$ and
$\eta_{\pow Kn}=\eta_n$, for all $n$.)
We have a commutative diagram
$$
\xymatrix@!C@R+2em{ R \ar[r]^{\eta_R} \ar[d] & \hull{R} 
\ar@{<-}[dl]^{\varphi} \\
{\hull n/I\hull n} }
$$
where the arrow on the left is the homomorphism obtained
from $\eta_n$ by base change modulo $I$. Hence by Proposition~\ref{P:ff}
the homomorphism $\eta_R$ is faithfully flat.
In the following we identify $\hull n/I\hull n$ and
$\hull R$ via the isomorphism $\varphi$.

Let $J$ be an ideal of $\pow K{n+m}$ with $I\subseteq J$.
The natural inclusion $\pow Kn\to\pow K{n+m}$ 
induces a morphism $\alpha\colon R\to S:=\pow K{n+m}/J$ in $\ancat$. (This is the 
only $\ancat$-morphism $R\to S$.)
Choose $\seq Jw\subseteq\pow {\seq 
Kw}{n+m}$ in the same way as we constructed the $\seq Iw$; so 
their ultraproduct is $J\hull{n+m}$ and $\hull S\iso 
\hull{n+m}/J\hull{n+m}$. Since $I\hull n\subseteq J\hull{n+m}$, we have
$I_w\subseteq J_w$ for almost all $w$, by \los. 
The natural inclusions $\pow{\seq Kw}n\to \pow{\seq Kw}{n+m}$ give rise
to homomorphisms
$$\seq\alpha w\colon \pow{\seq Kw}n/\seq Iw\to \pow{\seq Kw}{n+m}/\seq 
Jw.$$ 
The ultraproduct 
of the $\seq \alpha w$ yields a $\mathbf K$-algebra homomorphism 
$\hull\alpha\colon \hull R\to\hull S$ making the diagram
        \begin{equation*}
        \mbox{
        \xymatrix@C+3em@R+2em{{\hull n} \ar[d] \ar[r] & {\hull {n+m}} \ar[d] \\
{\hull R} \ar[r]^{\hull \alpha} & {\hull S}}
        }
        \end{equation*}
commutative. Together with \eqref{mnhull} this gives a commutative diagram
\begin{equation*}%\label{RS}
\begin{aligned}
\mbox{        \xymatrix@C+4em@R+2em{ R \ar[r]^{\alpha} 
\ar[d]^{\eta_R} & {S} \ar[d]^{\eta_S}\\
{\hull R} \ar[r]^{\hull \alpha} & {\hull S} }}
\end{aligned}
\end{equation*}
as required. \qed

\medskip

This concludes the proof  of Theorem~\ref{T:maincan} for the 
subcategory $\ancat$. Before we turn to the general case, we take a 
closer look at finite maps. We  use  the following version of 
the Weierstrass Division Theorem for $\hull{n+1}$. Let 
$f,g\in\hull{n+1}$ and suppose $g$ is
regular in $X_{n+1}$ of order $d$, that is,
        \begin{equation*}
        g \equiv X_{n+1}^d(1+\varepsilon)\mod (X_1,\dots,X_n)\hull{n+1}
        \end{equation*}
with $\varepsilon\in X_{n+1}\hull{n+1}$. Then there exist unique 
$q\in\hull{n+1}$ and $r\in\hull{n}[X_{n+1}]$ such that $f=qg+r$ and 
the degree of $r$ with respect to $X_{n+1}$ is strictly less than 
$d$. (Use \los\ and the  Weierstrass Division Theorem in 
$K_w[[n+1]]$.) A polynomial $P(Y)\in A[Y]$ in a single indeterminate 
$Y$
with coefficients in a local ring $(A,\mathfrak n)$ is called a
\emph{Weierstrass polynomial} if $P(Y)$ is monic of degree $d$ and
$P\equiv Y^d\bmod{\mathfrak n}A[Y]$.

\begin{proposition}\label{P:finan}
If $\alpha\colon R\to S$ is a \emph{finite} morphism in $\ancat $ 
\textup{(}that is to say, if $S$ is module-finite over $R$\textup{)},
then the natural map  $\hull{R}\tensor_R S\to\hull{S}$ is an 
isomorphism, making the diagram
        \begin{equation*}
        \xymatrix@R+2em@C+3em{ {\hull R} \ar[d]^{\hull\alpha} \ar[r] &
{\hull R\tensor_R S} \ar@{->}[dl]^{\iso}\\
{\hull S} }
        \end{equation*}
commutative.
\end{proposition}
\begin{proof}
We keep the notation from above, so that in particular $\alpha\colon 
R=\pow Kn/I\to S=\pow K{n+m}/J$. The case $m=0$ is clear.
By an induction on $m$, we may reduce to the case $m=1$. The ideal
        \begin{equation*}
        J_1:=J\cap K[[n]][X_{n+1}]
        \end{equation*}
of $K[[n]][X_{n+1}]$ contains a monic polynomial $P$. Now $P$ (as an 
element of $K[[n+1]]$) is regular of order at most $d=$ degree of 
$P$, hence can be written as $P=uQ$ where $u\in K[[n+1]]$ is a unit 
and $Q\in K[[n]][X_{n+1}]$
is a Weierstrass polynomial. Replacing $P$ by $Q$ we may assume that $P$ is
a Weierstrass polynomial of degree $d$ contained in $J_1$. The 
natural inclusion $K[[n]][X_{n+1}]\to K[[n+1]]$ induces an embedding
        \begin{equation}\label{eq:j1}
        K[[n]][X_{n+1}]/J_1 %=\pol R{X_{n+1}}/J_1\pol R{X_{n+1}}
\to K[[n+1]]/J=S,
        \end{equation}
which is in fact an isomorphism, for given $f\in K[[n+1]]$ we can 
write $f-r=qP\in J$ where $q\in K[[n+1]]$ and $r\in K[[n]][X_{n+1}]$ 
of degree $<d$, using Euclidean Division by $P$.  The image of 
$P$ under $\eta_{n+1}$, which we continue to denote by $P$, lies in 
$\hull n[X_{n+1}]$ and is a Weierstrass polynomial of degree $d$. The 
natural inclusion
$\hull{n}[X_{n+1}] \to \hull{n+1}$ induces a map
         \begin{equation}\label{e:surj}
         \pol{\hull n}{X_{n+1}}/P\pol{\hull n}{X_{n+1}}\to
\hull{n+1}/P\hull{n+1}.
         \end{equation}
 From the uniqueness of quotient and remainder in
Weierstrass Division by $P$ it follows that \eqref{e:surj} is in fact 
an isomorphism. Since  $J=J_1 K[[n+1]]$ and thus 
$J\hull{n+1}=J_1\hull{n+1}$, we get an induced
isomorphism
        \begin{equation}\label{eq:j1d}
        \pol{\hull n}{X_{n+1}}/J_1\pol{\hull n}{X_{n+1}}\to
\hull{n+1}/J\hull{n+1} =\hull{S}.
        \end{equation}
On the other hand, since $I\subseteq J_1$ and $R=K[[n]]/I$ we have
\begin{equation}
K[[n]][X_{n+1}]/J_1 \iso R[X_{n+1}]/J_1 R[X_{n+1}]
\end{equation}
and using $I\hull{n}\subseteq J_1\hull{n}[X_{n+1}]$
and $\hull{R}\iso\hull{n}/I\hull{n}$ we get
        \begin{equation}\label{eq:j1r}
        \hull n[X_{n+1}]/J_1\hull n[X_{n+1}] \iso \hull 
R[X_{n+1}]/J_1\hull R[X_{n+1}].
        \end{equation}
Therefore, by \eqref{eq:j1} and \eqref{eq:j1d}--\eqref{eq:j1r}:
         \begin{align*}
        \hull R\tensor_RS&\iso \hull R \tensor_R R[X_{n+1}]/J_1R[X_{n+1}]\\
        &\iso \pol{\hull R}{X_{n+1}}/J_1\pol{\hull R}{X_{n+1}}\\
        &\iso  \pol{\hull n}{X_{n+1}}/J_1\pol{\hull n}{X_{n+1}} \iso \hull S.
         \end{align*}
It is straightforward to check that we have a commutative diagram as claimed.
\end{proof}

Faithful flatness of $\eta_R\colon R\to\hull R$ now yields:

\begin{corollary}\label{C:surj}
If $S$ is a finite $R$-module via $\alpha$, then
$\alpha$ is injective \iff\ $\hull\alpha$ is injective, and
$\alpha$ is surjective \iff\ $\hull\alpha$ is surjective. \qed
\end{corollary}

\subsection{Proof of Theorem~\ref{T:maincan}}\label{s:proof}
We complete the proof of Theorem~\ref{T:maincan} by defining a 
functor $C\colon \cohcat \to\ancat$ and a faithfully flat natural 
transformation $\gamma\colon \fff\to  C$. The proclaimed 
$\mathfrak D$ and $\eta$     are then realized as the composite 
functor $\mathfrak D\after C$ and the natural transformation given by 
$\eta_{C(\Lambda)}\after\gamma_\Lambda$, for each object $\Lambda$. 
In essence, $C$ will be a kind   of `completion' functor. (See also 
\S\ref{s:extsc} below.)   More precisely, let $\Lambda=(R,\mathbf 
x,k,u)$ be an object in $\cohcat $ and let $k^*$ be the algebraic 
closure of $k$ in $\complet R$. Recall from \S\ref{s:qcf} that $k^*$ 
is the unique coefficient field of $\complet R$ containing $k$. We 
view $\complet R$ as a $k^*$-algebra via the inclusion 
$k^*\subseteq\complet R$. Let $\mathbf x=\rij xn$ and let 
$\theta_\Lambda\colon \pow {k^*}n\to \complet R$ be the $k^*$-algebra 
\homo\   given by $X_i\mapsto x_i$, for $i=\range 1n$. We  
denote its kernel by $I_\Lambda$. Consequently, we have associated to 
each $\Lambda$ a \emph{Cohen presentation} $\pow {k^*}n/I_\Lambda\iso 
\complet R$ of the completion of its underlying ring.

Let $\complet u\colon\complet R\to K$ be the completion of $u$ and 
denote the restriction of $\complet u$ to $k^*$ by $u^*$. There is a 
unique local \homo\ $\pow{k^*}n\to \pow Kn$ extending $u^*\colon 
k^*\to K$ and leaving the variables invariant, which we
denote by $u^*_n$. Define the functor $C$ on objects by the rule
        \begin{equation*}
        C(\Lambda):=\pow Kn/u^*_n(I_\Lambda)\pow Kn.
        \end{equation*}
As for morphisms, let $\Lambda\to \Gamma=(S,\mathbf y,l,v)$ be a 
morphism with underlying local \homo\ $\alpha\colon R\to S$.  Since 
$k^*$ (respectively, $l^*$) is the algebraic closure of $k$ in 
$\complet R$ (respectively, of $l$ in $\complet S$) and since 
$\alpha(k)\subseteq l$, the completion $\complet\alpha\colon \complet 
R\to\complet S$ of $\alpha$ 
maps $k^*$ inside $l^*$. Let us denote the restriction of $\complet\alpha$
to a field embedding  $k^*\to l^*$ by $\alpha^*$, and the induced map 
$\pow{k^*}n\to \pow{l^*}{n+m}$ leaving the variables $X_1,\dots,X_n$ 
invariant  by 
$\alpha^*_n$. Since  $\alpha(\mathbf x)$ is an initial segment of 
$\mathbf y$, we get $\theta_\Gamma(X_i)=\alpha(x_i)\in\complet S$
for $i=\range 1n$. Therefore, we have a commutative diagram
\begin{equation}\label{morp}
\begin{aligned}\mbox{
\xymatrix@R+2em@C+3em{\pow {k^*}n \ar[d]^{\alpha^*_n} 
\ar[r]^{\theta_\Lambda} & \complet R \ar[d]^{\complet\alpha} \\
\pow {l^*}{n+m} \ar[r]^{\theta_\Gamma} & \complet S.}}
\end{aligned}
\end{equation}
  In particular,  $\alpha^*_n(I_\Lambda)\subseteq I_\Gamma$. 
Since $v\after \alpha=u$, we get $\complet v \after \complet 
\alpha=\complet u$ which in turn yields $v^*_{n+m}\after 
\alpha^*_n=u^*_n$. Hence $u^*_n(I_\Lambda)\subseteq 
v^*_{n+m}(I_\Gamma)$ under the inclusion $\pow Kn\to\pow K{n+m}$, and 
this inclusion then induces a $K$-algebra \homo\ $C(\alpha)\colon 
C(\Lambda)\to C(\Gamma)$. Note that $C(\alpha)$ is indeed a morphism 
in $\ancat $.  Because every step in this construction is carried out 
in a canonical way,   $C$ is   a functor; details are left to the 
reader. Note that $C$ is the identity on the subcategory $\ancat$ of
$\cohcat$.

To define the natural transformation $\gamma\colon \fff\to C$ we let 
$\gamma_\Lambda$ be the composite map
        \begin{equation*}
        R\to \complet R\iso \pow{k^*}n/I_\Lambda\to \pow 
Kn/u^*_n(I_\Lambda)\pow Kn=C(\Lambda)
        \end{equation*}
where the isomorphism is induced by $\theta_\Lambda$ and the last 
arrow is the base change of $u^*_n$.  Each map in this composition is 
canonically defined and faithfully flat. It is now straightforward to 
check that  $\gamma$ is the desired faithfully flat natural 
transformation.\qed

\begin{remark}
It follows from our construction that the maximal ideals of 
$C(\Lambda)$ and $\hull\Lambda$ are $\maxim C(\Lambda)$ and 
$\maxim\hull \Lambda$, respectively, where $\maxim$ is the maximal 
ideal of the underlying ring of $\Lambda$.
\end{remark}

\begin{remark}\label{R:notacf}
If we do not insist that the ultraproducts are Lefschetz 
rings, then we can let $K$ be any ultraproduct of arbitrary fields, 
and Theorem~\ref{T:maincan} above, suitably
reformulated, remains
true in this more general setting, apart from the Lefschetz condition.
\end{remark}

\subsection{Extension of scalars}\label{s:extsc}
On occasion, we need a Lefschetz extension with some additional 
properties, and to achieve this, we enlarge the category 
$\cohcat $ to a category $\extcohcat$. To this end, we need a method 
to extend scalars. Suppose that we have a quasi-coefficient field $k$
of a Noetherian local ring $(R,\maxim)$ and a local 
\homo\ $u\colon R\to L$  to a field $L$. Let $k^*$ be the algebraic 
closure of $k$ in $\complet R$ (the unique coefficient field of
$\complet R$ containing $k$). We view $\complet R$ and $L$ as 
$k^*$-algebras via respectively the inclusion $k^*\subseteq \complet 
R$ and the restriction of $\complet u\colon \complet R\to L$ to 
$k^*$. Let $\complet R_{(k,u)}$ be the completion of the Noetherian
local ring $\complet 
R\tensor_{k^*}L$ with respect to its maximal ideal $\maxim(\complet 
R\tensor_{k^*}L)=\maxim \complet R\tensor_{k^*}L$. 
We view
$\complet R_{(k,u)}$ as an $R$-algebra 
(respectively, as an $L$-algebra) via the natural map 
$R\to \complet R\to \complet R\tensor_{k^*}L\to \complet R_{(k,u)}$ (respectively,  $L\to 
\complet R\tensor_{k^*}L\to \complet R_{(k,u)}$). 
The image of $L$ in $\complet R\tensor_{k^*}L$ is a coefficient field
of  $\complet R\tensor_{k^*}L$, and
hence of  the complete Noetherian local ring $\complet R_{(k,u)}$.
The $R$-algebra $\complet R_{(k,u)}$ is faithfully flat. 
The following transfer result will be used in the next section:
%For more about this construction see \cite[\S\S6.3 and 7.39]{HHFreg}.

\begin{lemma}\label{l:scalarext}
Suppose that $\op{char} k=0$. Then,
for given $i\in\mathbb N$, the completion
$\complet R$ of $R$ satisfies \textup{(R$_i$)}
\textup{(}or \textup{(S$_i$)}\textup{)} \iff\ 
$\complet R_{(k,u)}$ does.
In particular, $\complet R$ is reduced
\textup{(}regular, normal, or Cohen-Macaulay\textup{)}
if and only if $\complet R_{(k,u)}$ has this property.
Similarly, $\complet R$ is equidimensional \iff\ $\complet R_{(k,u)}$ is.
\end{lemma}
\begin{proof}
There is
probably a more straightforward way to see this, but we argue as follows:
Since $\operatorname{char} k=0$,  the homomorphism $\complet u|k^*\colon 
k^*\to L$
is separable. Therefore the induced homomorphism $\pow {k^*}n\to \pow Ln$ is 
formally smooth \cite[Theorem~28.10]{Mats}
hence regular \cite[p.~260]{Mats}. By the Cohen Structure Theorem 
we may assume $\complet R\cong k^*[[n]]/I$ (as $k^*$-algebras)
for some ideal $I$ of $k^*[[n]]$;
then $\complet R_{(k,u)} \cong L[[n]]/IL[[n]]$ (as $L$-algebras). Now use 
\cite[Theorem~23.9, and the remark following it]{Mats} to
conclude that  $\complet R$ satisfies (R$_i$)
\textup{(}or (S$_i$)\textup{)} \iff\ $\complet R_{(k,u)}$ has this property.
The local ring $\complet R_{(k,u)}$ is complete and hence catenary. Thus
if $\complet R_{(k,u)}$ is equidimensional, then $\complet R$ is equidimensional by
\cite[Theorem~31.5]{Mats}.  Conversely, if $\complet R$ 
is equidimensional, then so is $\complet R_{(k,u)}$ by \cite[(3.25)]{HHbigCM2}.
\end{proof}

\begin{remark}\label{r:scalarext}
Suppose that $\op{char} k=0$.
If $R$ is excellent, then $R\to\complet R$ is regular, hence
$R$ satisfies (R$_i$)
(or (S$_i$)) \iff\ $\complet R$ has this property, for every $i\in\mathbb N$.
Therefore,  if $R$ is excellent, then $R$ is reduced
(regular, normal, or Cohen-Macaulay)
if and only if $\complet R_{(k,u)}$ is.
By \cite{Ratliff}, 
$\complet R$ is equidimensional \iff\ $R$ is equidimensional and universally
catenary.
\end{remark}

Suppose we are given 
another Noetherian local ring $(S,\frak n)$ with quasi-coefficient
field $l$ of characteristic zero
and local homomorphim $v\colon S\to L$, as well as a
local homomorphism $\alpha\colon R\to S$ such that  $u=v\after 
\alpha$ and $\alpha(k)\subseteq l$.
 Since $\complet\alpha(k^*)\subseteq l^*$, 
we get natural maps
        \begin{equation}\label{eq:extsc}
        \complet R\tensor_{k^*} L \map{{\complet \alpha\tensor 1}} 
\complet S\tensor_{k^*}L=\complet 
S\tensor_{l^*}(l^*\tensor_{k^*}L)\to \complet S\tensor_{l^*}L,
        \end{equation}
where the last map is induced by the map 
$l^*\tensor_{k^*}L\to L$ given by $a\tensor b\mapsto v(a)b$, for 
$a\in l$ and $b\in L$. Taking completions yields an $L$-algebra
\homo\ $\complet R_{(k,u)}\to \complet S_{(l,v)}$, which we  
denote by $\complet\alpha_L$.

\subsection{The category $\extcohcat$}\label{s:extcohcat}
Let us first look at  an object $\Lambda=(R,\mathbf x,k,u)$ in $\cohcat$.
Applying
 the above construction with respect to the  \homo\ $u\colon R\to K$, 
we get a $K$-algebra $\complet R_{(k,u)}$ which is  isomorphic with 
$C(\Lambda)$; the isomorphism is uniquely determined by $\mathbf 
x$.  Allowing more general choices for $\mathbf x$ leads to
the extension $\extcohcat $ of  $\cohcat $. Namely, for  objects we 
take the quadruples $\Lambda=(R,\mathbf x,k,u)$ where as before 
$(R,\maxim)$ is a Noetherian local ring with quasi-coefficient field 
$k$ and $u\colon R\to K$ is a local \homo, but this time $\mathbf x$ 
is a tuple in the larger ring $\complet R_{(k,u)}$ generating its
maximal ideal $\maxim\complet R_{(k,u)}$. A morphism $\Lambda\to 
\Gamma=(S,\mathbf y,l,v)$ in this extended category  is 
given by a local \homo\ $\alpha\colon R\to S$ such that $u=v\after 
\alpha$, $\alpha(k)\subseteq l$, and such that $\complet\alpha_K
\colon \complet R_{(k,u)}\to \complet S_{(l,v)}$ sends 
$\mathbf x$ to an initial segment of $\mathbf y$. It is clear that 
$\cohcat $ is a full subcategory of $\extcohcat$.

\begin{remark}\label{r:iso}
Up to isomorphism, the $k^*$-algebra 
$\complet R_{(k,u)}$ is independent of the choice of $u$, 
since every isomorphism between
subfields of $K$ can be extended to an automorphism of $K$ (but not
necessarily to an ultraproduct of automorphisms of the $K_p$).
It is also easy to see that $\complet R_{(k,u)}$ is independent
of the choice of $k$, up to local isomorphism of local rings.
\end{remark}

We extend $C$ to a functor $\extcohcat \to \ancat $ as 
follows. Let $\mathbf x=\rij xn$ and let $\complet I_\Lambda$ be the 
kernel of the $K$-algebra \homo\ $\complet\theta_\Lambda\colon \pow 
Kn\to \complet R_{(k,u)}$ with $X_i\mapsto x_i$ for $i=\range 1n$. 
We now put
        \begin{equation*}
        C(\Lambda):=\pow Kn/\complet I_\Lambda.
        \end{equation*}
It follows that $C(\Lambda)\iso \complet R_{(k,u)}$. Note that if 
$\Lambda$ is an object of the subcategory $\cohcat $, then $\complet 
I_\Lambda=u^*_n(I_\Lambda)\pow Kn$ and $\complet\theta_\Lambda$ is 
the base change of $\theta_\Lambda$ over $u^*_n$, showing that 
$C(\Lambda)$ agrees with the $K$-algebra defined previously. As for 
morphisms, let  $\alpha\colon\Lambda\to\Gamma$ be as above. We have a 
commutative diagram
\[
\begin{aligned}
\mbox{
\xymatrix@R+2em@C+3em{\pow Kn \ar[d]  \ar[r]^{\complet\theta_\Lambda} 
& \complet R_{(k,u)} \ar[d]^{\complet\alpha_K} \\
\pow K{n+m} \ar[r]^{\complet\theta_\Gamma} & \complet S_{(l,v)}}
}
\end{aligned}
\]
where the first vertical arrow is the natural inclusion. It follows 
that $\complet I_\Lambda\subseteq \complet I_\Gamma$, thus giving 
rise to a morphism $C(\alpha)\colon C(\Lambda)\to C(\Gamma)$ in 
$\ancat$. It is now straightforward to verify that $C$ is a functor. 
Furthermore, the composition
        \begin{equation*}
        \gamma_\Lambda\colon R \to \complet R_{(k,u)} \iso C(\Lambda)
        \end{equation*}
is faithfully flat and hence yields a faithfully flat natural 
transformation $\gamma\colon\fff\to C$ (extending the previously 
defined natural transformation $\gamma$). From this discussion it is 
clear that we have the following extension of Theorem~\ref{T:maincan}:

\begin{theorem}\label{T:maincanext}
There exists a functor $\mathfrak D\colon \extcohcat \to \lefcat_{\mathcal W}$ 
and a faithfully flat natural transformation $\eta\colon \fff\to 
\fff\after\mathfrak D$. \qed
\end{theorem}

\subsection{Noether normalizations}\label{s:NN}
To explain the advantages of this extended version, we need to 
discuss Noether normalizations. 
Let $(A,\maxim)$ be a complete Noetherian local ring with coefficient field 
$k$. A $k$-algebra homomorphism $k[[d]]\to A$ 
which is finite and injective
is called a \emph{Noether normalization} of $A$. (Here necessarily $d=\dim A$.)
If $\mathbf x$ is an $n$-tuple generating $\maxim$ whose first $d$ entries
form a system of parameters of $A$, then the $k$-algebra homomorphism 
$k[[d]]\to A$ given by $X_i\mapsto x_i$ for $i=1,\dots,d$ is a Noether
normalization of $A$. (See for instance 
\cite[Theorem 29.4]{Mats}.)
However,  by choosing 
$\mathbf x$ even more carefully, we can achieve this also for 
homomorphic images: 

\begin{lemma}\label{L:nnorm}
Let $(A,\maxim)$ be a complete Noetherian local ring with an uncountable 
algebraically closed coefficient field $k$ and  let $\mathcal I$ be a 
set of proper ideals of $A$. If the cardinality of $\mathcal I$ is strictly less
than that of $k$, then there exists a surjective $k$-algebra \homo\ $\theta\colon
\pow kn\to A$
with the property that for every $I\in\mathcal I$, the $k$-algebra \homo\ $\pow 
kd\to A/I$ obtained by composing the restriction of $\theta$ to the 
subring $\pow kd$   with the natural surjection  $A\to A/I$ is a 
Noether normalization of $A/I$, where $d:=\dim A/I$.
\end{lemma}
\begin{proof}
Choose generators $y_1,\dots,y_n$ of $\maxim$ and let
        \begin{equation*}
        x_i = \sum_{j=1}^n a_{ij}y_j, \qquad i=1,\dots,n, \ \text{ 
and}\ a_{ij}\in k
        \end{equation*}
be general $k$-linear combinations of the $y_j$. By 
\cite[Theorem~14.14]{Mats} there exists, for every $I\in\mathcal I$, 
a non-empty Zariski open subset $U_I$ of $k^{n\times n}$ such that 
$x_1,\dots,x_d$ (where $d=\dim A/I$) is a system of parameters modulo 
$I$ for all $(a_{ij})\in U_I$. Since the transcendence degree of $k$ 
is strictly larger than $|\mathcal{I}|$, the intersection 
$\bigcap_{I\in\mathcal I} U_I$ is non-empty. Choose $(a_{ij})$ in 
this intersection and  let $\rij xn$ be the corresponding tuple. The 
$k$-algebra \homo\  $\theta\colon \pow kn\to A$    given by 
$X_i\mapsto x_i$ for all $i$ has the required properties.
\end{proof}

Let us express the property stated in the lemma by saying that 
$\theta$ is \emph{normalizing with respect to  $\mathcal 
I$}. Let $\Lambda=(R,\mathbf x,k,u)$ be an object in $\extcohcat $ 
and let $i$ denote the embedding of $k$ in the algebraic closure of 
$u(k)$ in  $K$ induced by $u$. The  natural \homo\ $R\to \complet 
R_{(k,u)}$ factors as $R\to \complet R_{(k,i)}\to \complet R_{(k,u)}$.
We say that 
$\Lambda$ is \emph{absolutely normalizing} if 
$\complet\theta_\Lambda \colon \pow Kn\to \complet R_{(k,u)}$ is 
normalizing with respect to the set of %all minimal primes of $\complet 
%R_{(k,u)}$ and with respect to 
all ideals of the form $I\complet 
R_{(k,u)}$ with $I$ an ideal in   $\complet R_{(k,i)}$. 
(This definition will be useful in \S\ref{s:app}.)
By Lemma~\ref{L:nnorm} and noting that the cardinality of 
$\complet R_{(k,i)}$ is at most 
$ 2^{\norm R}$, we immediately get:

\begin{corollary}\label{C:nnorm}
If we choose $K$ sufficiently large 
\textup{(}e.g., so that  $2^{\norm R}<\norm 
K$\textup{)}, then there exists an absolutely normalizing object in 
$\extcohcat $ with underlying ring $R$.\qed
\end{corollary}

We say that $\Lambda$ is {\em normalizing}
if the entries of the tuple $\mathbf{x}=(x_1,\dots,x_n)$
are in the maximal ideal $\maxim\complet R$ of
$\complet R$ and if the $k^*$-algebra homomorphism
$k^*[[n]]\to\complet R$ with $X_i\mapsto x_i$
for $i=1,\dots,n$ is normalizing with respect to the collection
$\mathcal I$ consisting of the zero ideal
and all minimal prime ideals of $\complet R$.
(As before $k^*$ denotes the algebraic closure of $k$ in $\complet R$.)
Since $\complet R$ is Noetherian, $\mathcal I$ is finite, and hence:

\begin{corollary}\label{C:nnorm2}
If $k$ is uncountable, then there exists a normalizing object
in $\extcohcat$ with underlying ring $R$. \qed
\end{corollary}

\begin{remark}\label{R:nnorm}
Let us discuss now how we intend to apply Theorem~\ref{T:maincan} and 
its extension, Theorem~\ref{T:maincanext}, in practice. 
With aid of a faithfully flat Lefschetz 
extension of an equi\ch\ zero Noetherian local ring $R$,   we'll 
define in \S\ref{s:tc} a non-standard tight closure relation on $R$, 
and   in \S\ref{s:BCM}, a big \CM\ algebra for $R$. If we only are 
interested in the ring $R$ itself, then no functoriality is 
necessary, and we remarked already that the proof in that case is 
much simpler, as it only relies on Proposition~\ref{P:LPK}.

Functoriality comes in when we are dealing with several rings at the 
same time, and when we need to compare the constructions made in each of 
these rings. We explain the strategy  in the case of a single local 
\homo\ $\alpha\colon R\to S$ between equicharacteristic zero Noetherian
local rings. Choose an algebraically closed 
Lefschetz field $K$ of sufficiently large cardinality (for instance larger than
$2^{\norm R}$ and $2^{\norm S}$) and choose an embedding of the residue field
$k_S$ of 
$S$ into $K$. Denote the compostion $S\to k_S\to K$ by $v$ and let 
$u:=v\after\alpha$. Choose a quasi-coefficient field $k$ of $R$ and 
then a quasi-coefficient field $l$ of $S$ containing $\alpha(k)$. 
Finally, choose a tuple $\mathbf x$ in $R$ generating its maximal 
ideal and enlarge the tuple $\alpha(\mathbf x)$ to a  generating 
tuple $\mathbf y$ of the maximal ideal of $S$. These data yield two 
objects $\Lambda:=(R,\mathbf x,k,u)$ and $\Gamma:=(S,\mathbf y,l,v)$ 
of $\cohcat $ and $\alpha$ induces a morphism between them. 
%(This shows that the forgetful functor $\fff\colon\cohcat \to\loccat 
%$ is `surjective'.) 
We take $\hull\Lambda$ and $\hull\Gamma$ as 
the faithfully flat Lefschetz extensions of $R$ and $S$ respectively, 
and use $\hull\alpha$ to go from one to the other. Of course, in this 
way, the closure operations defined on $R$ and $S$, and similarly, 
the big \CM\ algebras associated to them, depend  on the choices 
made, but this will not cause any serious problems.  Therefore, we 
will often simply denote the Lefschetz extensions by $\hull R$ and 
$\hull S$ with $\hull\alpha\colon\hull R\to \hull S$ the \homo\ 
between them.

For certain $\alpha$,   more adequate choices for the quadruples 
$\Lambda$ and  $\Gamma$ (and hence for the Lefschetz 
extensions $\hull R$ and $\hull S$) can be made.
For instance, this is the case if 
$\alpha$ is \emph{unramified}, that is to say, if the image of the
maximal ideal 
of $R$ generates the maximal ideal in $S$ and $\alpha$ induces an 
algebraic extension of the residue fields. In that case, we can take 
$l=\alpha(k)$ and  $\mathbf y=\alpha(\mathbf x)$. It follows that 
$C(\alpha)\colon C(\Lambda)\to C(\Gamma)$ is also unramified, whence 
$\hull\alpha$ sends the maximal ideal of $\hull R$ to the 
maximal ideal of $\hull S$. We'll tacitly assume that whenever 
$\alpha$ is unramified (for instance if $\alpha$ is surjective), 
then we choose $\hull R$ and $\hull 
S$ with these additional properties.
(See also \S\ref{s:quot} below.)

In the above construction of $\Lambda$, after we chose $k$ and $u$, 
we could have chosen the tuple $\mathbf x$ with entries in 
$\complet R_{(k,u)}$, so that the resulting $\Lambda$ is only an 
object in $\extcohcat $. This has the following advantage: 
by an application of Corollary~\ref{C:nnorm}, we now may choose 
$\Lambda$ so that it is absolutely normalizing. We express 
this by saying that the corresponding Lefschetz extension $\hull R$ 
$(:=\hull\Lambda)$  is \emph{absolutely normalizing}. 
Similarly, we also say that $\hull R$ is {\em normalizing} if
$\Lambda$ is normalizing.
One easily proves that if $\alpha\colon R\to S$ is a local homomorphism
as before, then $\Lambda$ with underlying ring $R$ and 
$\Gamma$ with underlying ring $S$ can be
chosen so that
$\alpha$ is a morphism $\Lambda\to\Gamma$ and
$\Gamma$ and $\Lambda$ are absolutely normalizing.
If moreover $\alpha$ is surjective and an absolutely normalizing object
$\Lambda$ of $\extcohcat$ with underlying ring $R$ is given, then
$\Gamma$ with underlying ring $S$ chosen as above is also absolutely normalizing.
\end{remark}

Next, we extend Proposition~\ref{P:finan}. We call a morphism 
$\Lambda\to\Gamma$ in $\extcohcat $ \emph{finite} if the underlying 
\homo\ $\alpha\colon R\to S$ is finite.

\begin{proposition}\label{P:fin}
If $\alpha\colon\Lambda\to\Gamma$ is a finite morphism in $\extcohcat$, 
then $\hull\alpha$ is also finite.
If in addition  $\alpha$ 
induces an isomorphism on the residue fields, then the natural map 
$$\hull\Lambda\tensor_RS\to \hull\Gamma$$ is an isomorphism,  making 
the diagram
        \begin{equation*}
        \xymatrix@R+2em@C+3em{ {\hull{\Lambda}} \ar[d]^{\hull\alpha} \ar[r] &
{\hull{\Lambda}\tensor_R S} \ar@{->}[dl]^{\iso}\\
{\hull{\Gamma}} }
        \end{equation*}
commutative. In particular,  $\alpha$ is injective 
\textup{(}respectively, surjective\textup{)} \iff\ $\hull\alpha$ is.
\end{proposition}
\begin{proof}
In view of Proposition~\ref{P:finan} and Corollary~\ref{C:surj}, it 
suffices to show the analogous statements with $\mathfrak D$ 
replaced by the functor $C$. If $\alpha\colon (R,\maxim)\to 
(S,\mathfrak n)$ is finite, then all the maps in \eqref{eq:extsc} are 
finite and hence so is $C(\alpha)$. Assume next that $\alpha$ induces 
an isomorphism on the residue fields. By the maximality property of 
coefficient fields we have $\complet\alpha(k^*)=l^*$. 
Since the canonical map $\complet R\tensor_RS\to 
\complet S$ is an isomorphism by \cite[Theorem 8.7]{Mats}, we get a 
canonical isomorphism
        \begin{equation*}
        (\complet R\tensor_{k^*}K)\tensor_RS\iso (\complet 
R\tensor_RS)\tensor_{k^*}K\iso\complet S\tensor_{k^*} K\iso \complet 
S\tensor_{l^*} K.
        \end{equation*}
Moreover, the $\maxim(\complet S\tensor_{l^*} K)$-adic topology on 
$\complet S\tensor_{l^*} K$ is equivalent with its $\mathfrak 
n(\complet S\tensor_{l^*} K)$-adic topology, since $\maxim S$ is 
$\mathfrak n$-primary. Hence taking completions and using 
\cite[Theorem 8.7]{Mats} once more, we get a canonical isomorphism
        \begin{equation*}
         \complet R_{(k,v)}\tensor_RS\iso \complet S_{(l,v)}.
        \end{equation*}
This in turn gives rise to a canonical isomorphism 
$C(\Lambda)\tensor_RS\iso C(\Gamma)$, which fits in a analogous 
commutative diagram as the above one.
\end{proof}

\subsection{Quotients}\label{s:quot}
Given an object $\Lambda=(R,\mathbf x,k,u)$ in $\extcohcat$ and an 
ideal $I$ of $R$, we define the \emph{quotient object} $\Lambda/I$ as 
the quadruple $(R/I,\bar{\mathbf x}, k,\bar u)$, where we identify 
$k$ with its image in $R/I$, where $\bar{\mathbf x}$ denotes the 
image of $\mathbf x$ in $\complet R_{(k,u)}/I\complet R_{(k,u)}$ and 
where $\bar u$ is the factorization of $u$ through $R/I$. The residue 
map $\pi\colon R\to R/I$ gives rise to a morphism 
$\Lambda\to\Lambda/I$. It follows from Proposition~\ref{P:fin} that 
$\pi$ induces a surjective map $\hull\Lambda\to \hull{\Lambda /I}$ 
and one easily checks that its kernel is  $I\hull\Lambda$.
If $\Lambda$ is absolutely normalizing, then so is $\Lambda/I$.

%Conversely, if $R$ is of the form $R=S/J$ for some Noetherian local
%ring $S$ and some ideal $J$ of $S$, and if $S$ is complete with 
%the same embedding
%dimension as $R$, then we can construct an object $\Gamma=(S,\mathbf y,l,v)$
%of $\extcohcat$ with underlying ring $S$ such that $\Gamma/J=\Lambda$.

\subsection{Further basic properties}
Recall from \S\ref{s:pow} that $\mathbf K=\hull 0 $ is just the ultrapower  
$K^{\mathcal U}$. By 
construction, $\mathbf K$ is a coefficient field of $\hull \Lambda$,
for every $\Lambda$,
and $\hull\alpha\colon\hull\Lambda\to\hull
\Gamma$ is a morphism of analytic $\mathbf K$-algebras with respect
to $\mathcal W$ (as defined in \S\ref{E:analytic-Lefschetz}), for
every morphism $\alpha\colon\Lambda\to\Gamma$ in $\extcohcat$.
The following is an analogue of
Proposition~\ref{P:sep}.

\begin{proposition}\label{P:infal}
For each $n$ there exists an exact sequence
         \begin{equation*}
         0\to {\infal{\hull n}} \to {\hull n} \xrightarrow{\pi}
{\pow{\mathbf K}n} \to 0
         \end{equation*}
where $\pi$ is a $K[[n]]$-algebra \homo.
\end{proposition}
\begin{proof}
Recall that $\infal{\hull n}$ denotes the ideal of infinitesimals of 
$\hull n$, that is to say, the intersection of all $\maxim^d\hull n$, 
where $\maxim :=\rij Xn\pow Kn$. Define $\pi\colon \hull n\to 
\pow{\mathbf K}n$ as follows. Take an
element $f\in \hull n$ and realize it as an ultraproduct of power
series $\seq fw\in\pow{\seq Kw}n$, say of the form
         \begin{equation*}
         \seq fw= \sum_\nu \seq{a_\nu} wX^\nu
         \end{equation*}
with $\seq{a_\nu}w\in\seq Kw$, where   $\nu$ ranges
over $\nat^n$.  For each such  $\nu$ let $a_\nu\in\mathbf K$ be
the ultraproduct of the $\seq {a_\nu}w$.   Define now $\pi(f)$ as the 
power series $\sum_\nu a_\nu X^\nu$. We leave it to the reader to 
verify that  $\pi$ is a
well-defined,  surjective   $\pol Kn$-algebra \homo, and that its
kernel is equal to $\infal{\hull n}$. (The argument is the same as in
the proof of Proposition~\ref{P:sep}.) It
remains to  show that it is in fact a $\pow Kn$-algebra \homo. Let
$f\in\pow Kn$ and choose polynomials $f_l\in\pol Kn$ so that $f\equiv
f_l\bmod\maxim^l$. It follows that
         \begin{equation*}
         \eta_n(f)\equiv\eta_n(f_l) = f_l\mod\maxim^l\hull n.
         \end{equation*}
Taking the image under $\pi$ shows that
         \begin{equation*}
         \pi(\eta_n(f))\equiv f_l \equiv f\mod \maxim^l\pow{\mathbf K}n.
         \end{equation*}
Since this holds for all $l$, we get that $\pi(\eta_n(f))=f$, proving
that $\pi$ is a $\pow Kn$-algebra \homo.
\end{proof}

\begin{remark}\label{r:partial}
The ultraproduct of the $i$-th partial derivative on each $\pow{\seq 
Kw}n$, for $i=\range 1n$, is a $\mathbf K$-linear endomorphism of 
$\hull n$, which we denote  again by $\partial/\partial X_i$. It
follows that
        \begin{equation*}
        \pi\left (\frac {\partial a}{\partial X_i}\right) = \frac 
{\partial (\pi(a))}{\partial X_i}
        \end{equation*}
for each $a\in \hull n$. In particular, for every $f\in\pow Kn$ we have
        \begin{equation*}
        \varepsilon(f):=
\eta_n\left(\frac {\partial f}{\partial X_i}\right) - \frac 
{\partial (\eta_n(f))}{\partial X_i}\in \infal{\hull n}.
        \end{equation*}
The map $f\mapsto \varepsilon(f)\colon K[[n]]\to \infal{\hull n}$
is a derivation which is trivial on $K[n]$.
We do not know whether $\varepsilon(f)=0$ for all $f\in K[[n]]$.
(Note that $\Omega_{K[[n]]/K[n]}\neq 0$.)
\end{remark}

\begin{corollary}\label{c:infini}
For each $\Lambda=(R,\mathbf x,k,u)$ in $\extcohcat$ we have an 
isomorphism of $R$-algebras
         \begin{equation*}
         \hull{R}/\infal{\hull R} \iso \pow{\mathbf K}n/\complet 
I_\Lambda\pow{\mathbf K}n\iso \complet R_{(k,\eta_0\after u)},
         \end{equation*}
where $\hull R:=\hull\Lambda$.
If $\mathfrak n$ is an $\maxim$-primary ideal of $R$, then $\mathfrak 
n\hull R$ is $\maxim\hull R$-primary and
\begin{equation}\label{e:primary}
\hull R/\mathfrak n\hull R \iso
(R/\mathfrak n)^{\mathcal U}.
\end{equation}
\end{corollary}
\begin{proof}
For the first statement use that the base change modulo $\complet 
I_\Lambda$ of the $\pow Kn$-algebra \homo\ $\pi$ from 
Proposition~\ref{P:infal} yields an epimorphism
$\hull R\to \pow{\mathbf K}n/\complet I_\Lambda\pow{\mathbf K}n$.
One verifies that its kernel is precisely $\infal{\hull R}$. The 
second isomorphism is then clear since $\pow Kn/\complet 
I_\Lambda\iso \complet R_{(k,u)}$. (Recall that $\eta_0\colon K\to 
\mathbf K=K^{\mathcal W}$ is the diagonal embedding.)

Now let $\mathfrak n$ be an $\maxim$-primary ideal of $R$,
say $\maxim^l\subseteq\mathfrak n$. Then $\maxim^l\hull 
R\subseteq\mathfrak n\hull R$, hence $\mathfrak n\hull R$ is
$\maxim\hull R$-primary. To establish \eqref{e:primary}
we first treat the case that $R=\pow Kn$. By Proposition~\ref{P:infal},
$\pi$ induces an isomorphism
         \begin{equation*}
         \hull n/\mathfrak n\hull n\iso \pow{\mathbf{K}}n/\mathfrak 
n\pow{\mathbf{K}}n.
         \end{equation*}
The natural homomorphism $\pow{\mathbf{K}}n\to (\pow Kn/\mathfrak 
n)^{\mathcal U}$ has kernel $\mathfrak n\pow{\mathbf{K}}n$ (use 
\los). The general case follows from this by base change.
\end{proof}

\subsection{A note of caution--- unnested conditions}\label{s:uc}
In the following we fix a natural number $n$. For a field $L$ and 
$i\in\{1,\dots,n\}$ let us write $$\pow
L{\remove\imath} := \pow
L{X_1,\dots,X_{i-1},X_{i+1},\dots,X_n}.$$ Let $\hull{\remove\imath}$
be the ultraproduct of the $\pow{\seq Kw}{\remove\imath}$. The
natural inclusion $\hull{\remove\imath}\subseteq \hull n$ is a section
of the canonical epimorphism
         \begin{equation*}
         \hull n\to \hull n/X_i\hull n\iso \hull{\remove\imath}.
         \end{equation*}
However, it is not true in general that $\eta_n\colon \pow Kn\to
\hull n$ maps $\pow K{\remove\imath}$ inside $\hull{\remove\imath}$
for all $i$ (the exception being of course $i=n$ by
\eqref{mnhull}). This is rather surprising since after all, $\eta_n$
sends a power series $f$ to a limit of its truncations in $\hull n$
and if $f$ does not involve $X_i$ then neither does each
truncation, yet the limit element must involve
$X_i$.% (forcing this monomial to be an infinitesimal).

To prove that such   inclusions cannot hold in general, we use an
example due to Roberts in \cite{RobSol}, which was designed to be a
counterexample to a question of Hochster on solid closure. Namely,
suppose for $n=6$, we would have   inclusions
         \begin{equation}\label{eq:compli}
         \eta_6(\pow K{\remove\imath})\subseteq \hull{\remove\imath}
         \end{equation}
for $i=4,5,6$. Let $z:=X_1^2X_2^2X_3^2$ and $a_i:=X_i^3$ for
$i=1,2,3$. Given a field $L$, the monomial $z$ lies in
the \emph{solid closure} of the ideal $(a_1,a_2,a_3)\pow L3$ \iff\
         \begin{equation*}
         f:= zX_4X_5X_6+a_1 X_5X_6+a_2 X_4X_6+a_3 X_4X_5 \in\zet[6]
         \end{equation*}
viewed as an element of $\pow L6$, has a non-zero multiple inside the
$L$-subspace
         \begin{equation*}
         \pow L{\remove 4}+\pow L{\remove 5}+\pow L{\remove 6}
         \end{equation*}
of $\pow L6$. (See \cite[\S9]{HoSol}.)
With Hochster we say that this  non-zero
multiple of $f$  is \emph{special}. If  \eqref{eq:compli} holds, then for
$L:=K$ the image under $\eta_6$
of such a non-zero multiple lies in the $K$-subspace
         \begin{equation*}
         \hull{\remove4}+\hull{\remove5}+\hull{\remove6}.
         \end{equation*}
of $\hull 6$.
By \los, $f$, as an element of $\pow{\seq Kw}6$, has then a non-zero
multiple which is special for almost all $w$. This in turn means that
$z$, viewed as an element of  $\pow{\seq Kw}3$, lies in the solid
closure of $(a_1,a_2,a_3)\pow{\seq Kw}3$. By \cite[Theorem
8.6]{HoSol}  solid closure
is trivial in $\pow{\seq Kw}3$ (since $\pow{\seq Kw}3$ is regular of
positive \ch).
Hence $z$ lies in
$(a_1,a_2,a_3)\pow{\seq Kw}3$, which is clearly false.

The failure of the existence of inclusions~\eqref{eq:compli} bears a
strong resemblance to the fact that there is no Artin Approximation
for unnested conditions (see our discussion in \S\ref{R:AA}).

\section{Transfer of Structure}\label{s:app}

Throughout this section $(R,\maxim)$ denotes  an equi\ch\ zero Noetherian 
local ring, and
$K$ is a Lefschetz field
with respect to some ultraset with underlying set equal to 
the set of the prime numbers, whose components $K_p$ are \acf{s} of 
\ch\ $p$. Whenever
necessary, we assume that $K$ has cardinality $>2^{\norm R}$. We fix once and for
all an object 
$\Lambda=(R,k,\mathbf{x},u)$
of $\extcohcat $  with  underlying ring  
$R$. (We might on 
occasion require some additional properties for $\Lambda$, such as being 
absolutely normalizing.) By abuse 
of notation,  we write $\hull R$ for $\hull\Lambda$. We 
view $\hull R$ as an $R$-algebra via the faithfully flat map 
$\eta_\Lambda\colon R\to \hull R$ and often surpress this map in our 
notation.
In particular, given an ideal $I$ in $R$, we simply write 
$I\hull R$ for the ideal in $\hull R$ generated by 
$\eta_\Lambda(I)$. Moreover, we construct a Lefschetz hull for 
$R/I$ always by means of the quotient $\Lambda/I$, as explained 
in \S\ref{s:quot}. In particular,     $\hull{R/I}\iso \hull 
R/I\hull R$. The other notations introduced in \S\ref{s:LP} remain in force.

\subsection{\Sr{s}}
By construction, $\hull R$ is an ultraproduct  (with respect to some 
unspecified ultraset) of equi\ch\ complete Noetherian local rings 
$\seq Rw$ with algebraically closed residue field $\seq Kw$ (of prime 
\ch\ $p(w)$). We think of  $\seq Rw$ as an \emph\sr{} of $R$. 
Each $\seq Rw$ is of the form $\pow{\seq Kw}n/\seq Iw$, where $\seq 
Iw$ are ideals whose ultraproduct is equal to $\complet 
I_\Lambda\hull n$ (in the notation of \S\ref{s:extcohcat}).
In this section, we make more precise how the $\seq Rw$ play
the role of a reduction modulo $p$ of $R$. A similar study for affine 
$K$-algebras was carried out in \cite{SchNSTC} and the subsequent 
papers, using effective bounds and the resulting
first-order definability (as established in \cite{SvdD,SchBC}). Since 
no such tool is available in the present situation, our arguments 
are purely algebraic.
Here is a first example:

\begin{theorem}%[Transfer of Structure]
\label{T:trans}
%Let $R$ be an equi\ch\ zero Noetherian local ring   and let $\seq Rw$ be 
%an \sr\ of $R$.
\ 
\begin{enumerate}
\item\label{i:dim} Almost all $\seq Rw$ have the same dimension 
\textup{(}respectively, embedding dimension or depth\textup{)} as $R$.
\item\label{i:trans} Almost all $\seq Rw$ are regular \textup{(}respectively, \CM\ or 
Gorenstein\textup{)} \iff\ $R$ has the same property.
\end{enumerate}
\end{theorem}

Before we begin the proof, let us introduce some more notations. 
Given an element 
$a\in\hull R$ choose elements $\seq aw\in\seq Rw$ whose 
ultraproduct is $a$. We call $a_w$ an \emph\sr\ of $a$. If $\seq aw'$ is 
another choice of elements whose ultraproduct is $a$, then $\seq 
aw=\seq aw'$ for almost all $w$. We use similar terminology for tuples
of elements in $\hull R$, and given a finitely 
generated ideal  $I=\rij as\hull R$ of $\hull R$, let 
$\seq Iw:=(\seq{a_1}w,\dots,\seq{a_s}w)\seq Rw$, where $\seq{a_i}w$ 
is an \sr\ of $a_i$. The ultraproduct of the $\seq 
Iw$ is $I$, and we call  $\seq Iw$ an \emph\sr\ of $I$. 
If we choose different generators and \sr{s} of these 
generators and denote  the resulting  ideals by $\seq Iw'$, then the 
ultraproduct of the $\seq Iw'$  is again $I$ and therefore  $\seq 
Iw=\seq I w'$ for almost all $w$.  With an \emph\sr\ of an ideal 
$I$ of $R$ we mean  an \sr\ of its extension $I\hull R$ to an ideal
of $\hull R$. Note that then $\seq Rw/\seq Iw$ is an
\sr\
of $R/I$. By faithful flatness of $R\to\hull R$ we have:

\begin{lemma}\label{L:ff}
If $I$ and $J$ are ideals of $R$ with \sr{s} $\seq Iw$ and $\seq Jw$, then 
\begin{enumerate}
\item\label{i:mem} $I\hull R\cap R=I$,
\item\label{i:cap} $I\hull R\cap J\hull R=
(I\cap J)\hull R$, 
\item\label{i:div} $\bigl(I\hull R:_{\hull R}J\hull R\bigr)=(I:_RJ)\hull R$,
\end{enumerate} 
and the ideals in \eqref{i:cap} and \eqref{i:div} have
approximations $I_w\cap J_w$ and $(I_w:_{R_w}J_w)$, respectively.
\end{lemma}

Let $\maxim$ be the maximal ideal of $R$. 
As a first step in the proof of Theorem~\ref{T:trans} we show 
the following lemma, of interest in its own right:

\begin{lemma}\label{L:sop}
A $d$-tuple $\mathbf z =(z_1,\dots,z_d)\in R^d$ 
is a system of parameters for $R$
\iff\ almost every $\seq{\mathbf z}w$ is a system of
parameters for $R_w$, where 
$\seq {\mathbf z}w=(z_{1w},\dots,z_{dw})$ is an \sr{} of $\mathbf z$.
Similarly, $\mathbf{z}$ is an $R$-regular sequence \iff\ 
$\mathbf{z}$ is a $\hull{R}$-regular sequence \iff\ 
almost every $\seq{\mathbf z}w$ is an $R_w$-regular sequence.
\end{lemma}
\begin{proof}
Suppose $\mathbf z$ is a system of parameters for $R$, so $d=\dim R$. 
We claim that almost every $\mathbf{z}_w$ is a system of parameters for $R_w$. 
We have
$\maxim^r\subseteq (z_1,\dots,z_d)R$ for some $r$, and since 
this is preserved in 
$\hull R$, we get by \los\ that $\seq\maxim w^r\subseteq
(z_{1w},\dots,z_{dw})R_w$, 
for almost all $w$. This shows that almost all $R_w$ have dimension at
most $d$, and it suffices to shows that $\dim R_w=d$ for almost all $w$. 
Suppose on the contrary that $\dim R_w<d$ for
almost all $w$. We may assume,
after renumbering if necessary, that the ideal $\seq{\mathfrak 
n}w:=(\seq{z_1}w,\dots,\seq{z_{d-1,}}w)\seq Rw$ of $R_w$ is $\seq\maxim 
w$-primary for almost all $w$.  For those $w$
let $\seq rw\in\mathbb N$ be minimal such that $(\seq{z_d}w)^{\seq 
rw}\in\seq{\mathfrak n}w$.  By Noetherianity of $R$, we have for 
some $s$ that
        \[
        (\mathfrak n:z_d^s)=(\mathfrak n:z_d^r)
        \]
for all $r\geq s$, where $\mathfrak n:=\rij z{d-1}R$. By
\eqref{i:div} we get
        \begin{equation}\label{eq:ann}
        \begin{aligned}
        \bigl(\mathfrak n\hull R:_{\hull R}z_d^r\bigr) &= (\mathfrak n:_R 
z_d^r)\hull R \\
        &=  (\mathfrak n:_R z_d^s)\hull R \\
        &=  \bigl(\mathfrak n\hull R:_{\hull R}z_d^s\bigr),
        \end{aligned}
        \end{equation}
for all $r\geq s$.  Suppose $\seq rw>s$ for almost all $w$, and let
$b\in\hull R$ equal the ultraproduct of the $(\seq{z_d}w)^{\seq 
rw-s-1}$. By \los, $bz_d^{s+1}\in\mathfrak n\hull R$. By 
\eqref{eq:ann}, we have $ bz_d^s\in\mathfrak n\hull R$ and hence, by 
\los\ once more, $(\seq{z_d}w)^{\seq rw-1}\in \seq{\mathfrak n}w$ for 
almost all $w$, contradicting the minimality of $\seq rw$. Therefore, 
$\seq rw\leq s$ and hence $(\seq{z_d}w)^s\in\seq{\mathfrak n}w$, for 
almost all $w$. By \los, this yields $z_d^s\in\mathfrak n\hull R$ 
and hence  $z_d^s\in\mathfrak n$ by faithful
flatness of $R\to\hull R$, contradicting 
that $\mathbf z$ is a system of parameters for $R$.
Conversely, assume that $\mathbf{z}_w$ is a system of parameters for
$R_w$ for almost all $w$. Then $\dim R_w=d$ for almost all $w$. We have
already shown $\dim R_w=\dim R$ for almost all $w$, hence $\dim R=d$. 
Therefore it suffices to show that $(z_1,\dots,z_d)R$ is $\frak m$-primary. 
Now $(z_{1w},\dots,z_{dw})R_w$ is $\frak m_w$-primary, hence 
$\dim R_w/(z_{1w},\dots,z_{dw})R_w=0$  for almost all $w$.
The rings $S_w:=R_w/(z_{1w},\dots,z_{dw})R_w$ are
approximations to $S:=R/(z_1,\dots,z_d)R$. Thus 
$\dim S=0$ by what we have shown above, or equivalently,
$(z_1,\dots,z_d)R$ is $\frak m$-primary.

If $\mathbf z$ is $R$-regular, then $\mathbf z$ is also
$\hull{R}$-regular due to faithful flatness of $R\to\hull{R}$,
see \cite[Exercise 16.4]{Mats}. By \los, if $\mathbf{z}$ is
$\hull{R}$-regular, then almost all $\mathbf{z}_w$ are $R_w$-regular.
%Now suppose that $\mathbf x$ is an $R$-regular sequence.
%We claim that almost all 
%$(\seq{x_1}w,\dots,\seq{x_n}w)$ are $\seq Rw$-regular. Indeed, if 
%not, then there exists some $i\in\{1,\dots,n-1\}$ such that for almost all $w$
%we find $a_w\in R_w\setminus (\seq{x_1}w,\dots,\seq{x_i}w)\seq Rw$
%with $\seq aw\seq{x_{i+1,}}w\in (\seq{x_1}w,\dots,\seq{x_i}w)\seq Rw$. 
%Using \los, this gives $ax_{i+1}\in\rij xi\hull R$, where $a$ is 
%the ultraproduct of the $\seq aw$. Hence $a$ lies in the ideal
%        \begin{align*}
%        \bigl(\rij xi\hull R:_{\hull R} x_{i+1}\bigr)&= \bigl(\rij xiR:_R x_{i+1}\bigr)\hull R\\
%        &=\rij xi\hull R
%        \end{align*}
%where the first equality follows from \eqref{i:div} and the last from 
%the fact that $\mathbf x$ is $R$-regular. By \los\ once more, it 
%follows that $\seq aw\in (\seq{x_1}w,\dots,\seq{x_i}w)\seq 
%Rw$ for almost all $w$, a 
%contradiction. 
Finally, suppose that almost all $\mathbf{z}_w$ are
$R_w$-regular, and let $i\in\{1,\dots,d-1\}$ and $a\in R$ with
$az_{i+1}\in (z_1,\dots,z_i)R$. Then we have 
$a_wz_{i+1,w}\in(z_{1w},\dots,z_{iw})R_w$
for almost all $w$, hence $a_w\in (z_{1w},\dots,z_{iw})R_w$ for almost all
$w$, and therefore $a\in (z_1,\dots,z_i)\hull R$, by \los. 
Now \eqref{i:mem} yields
$a\in (z_1,\dots,z_i)R$. Similarly one shows that $1\notin (z_1,\dots,z_d)R$.
Hence $\mathbf{z}$ is $R$-regular.
\end{proof}

\begin{proof}[Proof of Theorem~\ref{T:trans}]

%\subsubsection*{Dimensions}
Suppose that $R$ has embedding dimension $e$, so that we can write 
$\maxim=\rij zeR$ for some $z_1,\dots,z_e\in R$. Hence $\seq\maxim 
w=(\seq{z_1}w,\dots,\seq{z_e}w)\seq Rw$, where $\seq{z_i}w$ is an 
\sr\ of $z_i$. If the  embedding dimension of almost 
all $\seq Rw$ would be strictly less than $e$, then after 
renumbering if necessary, 
$\seq \maxim w=(\seq{z_1}w,\dots,\seq{z_{e-1,}}w)\seq 
Rw$ for almost all $w$ (by Nakayama's Lemma). Therefore $\maxim\hull R=\rij 
z{e-1}\hull R$ by \los,   hence $\maxim =\rij z{e-1} R$ by faithful 
flatness of $R\to\hull R$, contradiction.
From Lemma~\ref{L:sop} %(or Corollary~\ref{C:Noeth}) 
it follows that $\dim R=\dim R_w$ for almost all $w$.
Now suppose $R$ has depth $d$, and let $\mathbf{z}=(z_1,\dots,z_d)$,
with $z_i\in\frak m$ for all $i$ be an 
$R$-regular sequence. 
By Lemma~\ref{L:sop} almost every approximation
$\mathbf{z}_w\in\mathfrak{m}$ of $\mathbf{z}$ 
is an $R_w$-regular sequence in $\frak{m}_w$ and
hence $R_w$ has depth at least $d$, for almost all $w$.
On the other hand, since $R$ has depth $d$, the quotient $R/\rij zdR$ 
has depth zero, that is to say, $\maxim$ is an associated prime of 
$\rij zdR$. Choose $s\notin\rij zdR$ such that $s\maxim\subseteq \rij 
zdR$ and let $\seq sw$ be an \sr\ of $s$. By \los, $\seq sw\seq\maxim 
w\subseteq (\seq{z_1}w,\dots,\seq{z_d}w)\seq Rw$ and $\seq sw\notin 
(\seq{z_1}w,\dots,\seq{z_d}w)\seq Rw$ for almost all $w$. Hence
the depth of  almost all $\seq Rw$ equals $d$.

%\subsubsection*{Singularities}
Since $R$ is regular (respectively, \CM) \iff\ its dimension is equal 
to its embedding dimension (respectively, to its depth), the desired 
transfer follows from the preservation of these invariants in the 
\sr{s}. Let $d:=\dim R$, and recall that
$R$ Gorenstein means that $R$ is \CM\ and for some (equivalently, every)
$R$-regular sequence $\mathbf{z}=\rij zd$ in $\mathfrak{m}$,
the socle of $R/\mathfrak n$ is principal, where $\mathfrak n:=\rij zdR$;
that is, there exists $a\in R$ such that $(\mathfrak 
n:\maxim)=\mathfrak n+aR$. In order to show that $R$ is Gorenstein \iff\ almost
all $R_w$ are, we may assume, by our agument above, that $R$ and hence
almost all $R_w$ are \CM.
Suppose that $R$ is Gorenstein. Let $\mathbf{z}$, $\frak{n}$, and 
$a$ as above, and 
let $\seq aw$ and $\seq{\mathfrak n}w$ be 
\sr{s} of  $a$ and $\mathfrak n$ respectively, so $\frak{n}_w$ is 
generated by an $R_w$-sequence, for almost all $w$. By \los, we get
        \begin{equation}\label{eq:gor}
        (\seq{\mathfrak n}w:\seq\maxim w)=\seq{\mathfrak n}w +\seq aw\seq Rw,
        \end{equation}
for almost all $w$. It follows that almost all $\seq Rw$ are 
Gorenstein. 
Conversely, if  almost all $\seq Rw$ are Gorenstein, then  there 
exist $\seq aw\in R_w$ satisfying  \eqref{eq:gor}. By \los,
        \begin{equation}\label{eq:gorul}
                \bigl(\mathfrak n\hull R:_{\hull R}\maxim\hull R\bigr)=\mathfrak n\hull 
R+  \ul a\hull R
        \end{equation}
where $\ul a\in\hull R$ is the ultraproduct of the $\seq aw$. 

Let $f$ 
and $g$ be elements in $(\mathfrak n:\maxim)$ but not in $\mathfrak 
n$. From \eqref{eq:gorul} it follows that $f\equiv \ul a\ul 
b\bmod\mathfrak n\hull R$ and $g\equiv \ul a\ul c\bmod\mathfrak n\hull 
R$, for some $\ul b,\ul c\in \hull R$. By faithful flatness of
$R\to \hull R$, neither 
$f$ nor $g$ belongs to $\mathfrak n\hull R$, so that $\ul b$ and $\ul 
c$ must be units in $\hull R$. In particular, $f\in g\hull 
R+\mathfrak n\hull R$ and $g\in f\hull R+\mathfrak n\hull R$. 
Therefore, again by faithful flatness,  $f\in gR+\mathfrak n$ and 
$g\in fR+\mathfrak n$. Since this holds for every choice of $f$ and 
$g$,   the socle of $R/\mathfrak n$ is principal, showing that 
$R$ is Gorenstein. 
\end{proof}

Since a Noetherian local ring is a discrete valuation ring (DVR) \iff\ 
it has
positive dimension and its maximal ideal is principal \cite[Theorem 11.2]{Mats},
we get:

\begin{corollary}\label{C:DVR}
The following are equivalent:
\begin{enumerate}
\item $R$ is a DVR;
\item almost every $R_w$ is a DVR;
\item $\hull R$ is a valuation ring. \qed
\end{enumerate}
\end{corollary}

\subsection{Flatness and Noether normalization}
Let $\Gamma$ be an object in $\extcohcat$ with underlying 
ring $S$, and $\Lambda\to\Gamma$ a morphism in $\extcohcat$ with
underlying \homo\ $\alpha\colon R\to S$.
We denote the induced morphism $\hull R\to \hull S:=\hull\Gamma$ 
by $\hull\alpha$. By definition, $\hull \alpha$ is an ultraproduct of 
$K_w$-\homo{s} $\seq\alpha w\colon \seq Rw\to \seq Sw$, where $\seq 
Sw$ is an \sr\ of $S$.

\begin{proposition}\label{P:flat}
If $\alpha\colon R\to S$ is finite, then so are almost all 
$\seq\alpha w$. If $\alpha$ moreover induces an isomorphism on the 
residue fields, then
the following are
equivalent:
\begin{enumerate}
\item\label{flat1} $\alpha$ is flat;
\item\label{flat2} $\hull\alpha$ is flat;
\item\label{flat3} almost all $\seq\alpha w$ are flat.
\end{enumerate}
\end{proposition}
\begin{proof}
The first assertion and the implication 
\eqref{flat1}~$\Rightarrow$~\eqref{flat2} are immediate by 
Proposition~\ref{P:fin}. From the commutative diagram \eqref{fun} and 
the faithful flatness of  $\eta_R:=\eta_\Lambda$ and 
$\eta_S:=\eta_\Gamma$ we get 
\eqref{flat2}~$\Rightarrow$~\eqref{flat1}. Hence remains to show that 
\eqref{flat1} and \eqref{flat3} are equivalent. We use the local 
flatness criterion
\cite[Theorem~22.3]{Mats}: \textsl{a finitely generated module $M$ over
a  local Noetherian ring $(A,\frak n)$ is flat \iff\ $\operatorname{Tor}_1^A(A/\frak n,M)=0$}.
Since $\eta_R$ is flat we have an isomorphism of $\hull R$-modules
        \begin{equation*}
        \hull R \tensor_R \operatorname{Tor}_1^R(R/\maxim,S)
\iso\operatorname{Tor}_1^{\hull R}\bigl(\hull R \tensor_R (R/\maxim),\hull R
\tensor_R S\bigr).
        \end{equation*}
Moreover $\hull R \tensor_R (R/\maxim) \iso \hull{R/\maxim}=\hull K$, and
$\hull R \tensor_R S \iso \hull S$ by Proposition~\ref{P:fin}.
The finitely generated $R$-module $S$ has a free resolution by
finitely generated free $R$-modules (since $R$ is Noetherian). Hence
by the faithful flatness of $R\to \hull R$, the finitely generated 
$\hull R$-module
$\hull S$ has a free resolution by finitely generated free $\hull R$-modules.
Since $\hull K$ is a field and hence coherent,
   Proposition~\ref{P:Tor} yields
        \begin{equation*}
        \hull R \tensor_R \operatorname{Tor}_1^R(R/\maxim,S) \iso
\up w \operatorname{Tor}_1^{R_w}(K_w,S_w).
        \end{equation*}
The Noetherian local ring $R_w$ has residue field $K_w$, and
$S_w$ is finitely generated as a module over $R_w$, for almost all $w$.
The claim now follows from the local flatness criterion and
faithful flatness of $R\to \hull R$.
\end{proof}

\begin{proposition}\label{P:Noeth}
Let $I$ be an ideal 
in $R$ with \sr{s} $\seq I w\subseteq \seq Rw$, and $d=\dim R/I$.
If  $\hull R$ is absolutely normalizing, then
the composition
        \begin{equation*}
        \pow{\seq Kw}d\subseteq\pow{\seq Kw}n\to \seq Rw\to \seq Rw/\seq I w
        \end{equation*}
\textup{(}where the first map is given by inclusion and the remaining maps are 
the natural surjections\textup{)} is a Noether normalization of $\seq Rw/\seq 
I w$, for almost all $w$.
\end{proposition}
\begin{proof}
By Remark~\ref{R:nnorm}, the natural $K$-algebra homomorphism 
$$\pow Kd\to \pow Kn \overset{\complet\theta_\Lambda}{\longrightarrow} 
C(\Lambda)\to
C(\Lambda)/I C(\Lambda)=C(\Lambda/I)$$
is injective and finite, hence a Noether normalization. By 
Proposition~\ref{P:finan},
applying $\mathfrak D$ yields a finite and injective \homo\ $\hull 
d\to \hull{R/I}$. By \los, the maps in the statement are therefore 
almost all injective and finite, since their ultraproduct is 
precisely $\hull d\to \hull{R/I}$.
\end{proof}

\begin{remark}
Suppose $I=(0)$. Then the conclusion
of the proposition holds if $\hull R$ is only assumed to be normalizing.
\end{remark}

Let us elaborate some more on the Proposition~\ref{P:Noeth}. 
Let $T:=K[[d]]$, where $0\leq d\leq n$.
We have a commutative diagram
$$\xymatrix@R+2em{T[X_{d+1},\dots,X_{n}]\ar[r]\ar[d] & 
\hull{d}[X_{d+1},\dots,X_n]\ar[d]\\
K[[n]]\ar[r]^{\eta_n} & \hull{n}.
}$$
Hence given $f\in T[X_{d+1},\dots,X_{n}]$ we may choose approximations
$f_w$ of $f$ in the subring
$T_w[X_{d+1},\dots,X_{n}]$ of  $\pow{\seq Kw}n$. 
Note that $\seq Tw:=\pow{\seq Kw}d$ is an \sr\ of $T$, so that
$f$ is  the ultraproduct of the polynomials $f_w$ of bounded degree, 
in the sense of \S\ref{s:ulpol}.
Given generators $f_1,\dots,f_r$ of an ideal $J$ of $T[X_{d+1},\dots,X_{n}]$
we let $J_w$ be the ideal of  $T_w[X_{d+1},\dots,X_{n}]$ generated by
$f_{1w},\dots,f_{rw}$. (We think of $J_w$ as an approximation of the ideal
$J$.)

%The following definition is useful:
%
%\begin{definition}%\label{d:ws}
%Given a local ring $(A,\frak n)$, a \emph{Weierstrass sequence}
%in  $A[Y_1,\dots,Y_m]$ is a sequence $(g_1,\dots,g_m)$
%where each $g_i$ is a polynomial of the form 
%$$g_i=Y_i^{d_i} + \sum_{j=0}^{d_i-1}
%a_{ij}Y_i^j \in A[Y_1,\dots,Y_{i}]$$ with $d_i\in\nat$ and $a_{ij}\in 
%(\frak n,Y_1,\dots,Y_{i-1}) A[Y_1,\dots,Y_{i-1}]$ for all $i,j$.
%(So in particular, $g_1\in A[Y_1]$ is a Weierstrass polynomial.) 
%\end{definition}
%
%\begin{remark}\label{r:ws}
%If $A=L[[X]]=L[[X_1,\dots,X_d]]$, where $L$ is a field, and if $M$ is
%an ideal of $A[Y]=A[Y_1,\dots,Y_m]$ which contains a Weierstrass
%sequence, then $M=ML[[X,Y]]\cap A[Y]$.
%(See the argument used in the proof of Proposition~\ref{P:finan}.)
%\end{remark}

Suppose that $\hull R$ is normalizing.
Recall that we denote the kernel of $\complet\theta_\Lambda$ by
$\complet I_\Lambda$, and
let $J:=T[X_{d+1},\dots,X_n]\cap \complet I_\Lambda$, 
where as above $T=K[[d]]$, with $d:=\dim R$. 
Then $J\cap T=(0)$, and the natural inclusion $T[X_{d+1},\dots,X_n]\to
K[[n]]$ induces an isomorphism
$$
T[X_{d+1},\dots,X_n]/J \to K[[n]]/\complet I_\Lambda=C(\Lambda).
$$
Hence for every ideal $M$ of $T[X_{d+1},\dots,X_n]$ containing $J$, we have
$M=MK[[n]]\cap T[X_{d+1},\dots,X_n]$. By the remark following the
proposition above, we see that then
$M_w=M_wK_w[[n]]\cap T_w[X_{d+1},\dots,X_n]$ for almost
all $w$. This fact is used in the proof of Theorem~\ref{T:normal} below.
%By Proposition~\ref{P:finan}, applying $\mathfrak D$   yields an isomorphism
%$$\hull d [X_{d+1},\dots,X_{n}]/J \hull d [X_{d+1},\dots,X_{n}]
%\to\hull{n}/I\hull{n}=\hull{R}.$$
%and almost all of the
%$K_w$-algebra homomorphisms $$T_w:=K_w[[d]]\to R_w:=K_w[[n]]/I_w$$
%whose ultraproduct is $\hull{d}\to\hull{R}$ finite and
%are injective. 
%Since $C(\Lambda)$ is finite over $T$,
%the ideal $J$ contains a Weierstrass sequence $(g_{1},\dots,g_{m})$,
%where $m=n-d$. By \los\ 
%$(g_{1w},\dots,g_{mw})$ is a Weierstrass
%sequence in $J_w$, for almost all $w$. In particular, if $M$
%is an ideal of $T[X_{d+1},\dots,X_n]$ containing $J$, then
%$M_w = M_w K_w[[n]]\cap T_w$ for almost all $w$, by Remark~\ref{r:ws}.
%Proposition~\ref{P:Noeth} implies this for $M$ which correspond to
%ideals of the form $\id\complet R_{(k,u)}$ for an ideal $\id$ of
%$\complet R_{(k,i)}$ under the isomorphism \eqref{Niso}.

\begin{remark}\label{r:nn}
Because of its importance, let us give alternative arguments for 
\eqref{i:trans} using Noether normalization. These arguments work if
$\Lambda=(R,\mathbf{x},k,u)$ 
is absolutely normalizing. (In fact, it is enough that
$\complet\theta_\Lambda$ be normalizing with respect to the zero ideal.)
First, in all three cases we may replace $R$ by 
$\complet R_{(k,u)}$  
and hence assume that $R$ is in $\ancat$.
(See Theorem~23.7, Corollary to Theorem~23.3, and Theorem~23.4, respectively, 
in \cite{Mats}.) Say $R=K[[n]]/I$
for some $n$ and
some ideal $I$ of $K[[n]]$.
The restriction of the $K$-algebra homomorphism
$K[[n]]\to R$ with $X\mapsto \mathbf{x}$
to $T:=K[[d]]$, where $d:=\dim R$, is a Noether normalization of $R$. By 
Proposition~\ref{P:Noeth}, $T_w\to R_w$  is
a Noether normalization of $R_w$, for almost all $w$. 
The proof of \cite[Theorem~29.4]{Mats} shows
that $R$ is regular \iff\ $T\to R$ is surjective. 
Hence $R$ is regular \iff\ $\hull T\to \hull R$ is surjective
(by Corollary~\ref{C:surj})
\iff\ $T_w\to R_w$ is surjective 
for almost all $w$. Therefore $R$ is regular \iff\ almost 
each $R_w$ is regular. By  \cite[Proposition~2.2.11]{Bruns-Herzog},
$R$ is \CM\ \iff\ $T\to R$ is flat.
By Proposition~\ref{P:flat} this is equivalent with the 
flatness of almost all $T_w\to \seq Rw$, which in turn 
is equivalent with almost all $\seq Rw$ being \CM. Finally, for the 
Gorenstein property, observe that  $R/\mathfrak n$ is Artinian, 
where $\frak n$ is a parameter ideal of $R$, 
and so is $\hull {R/\mathfrak n}$, as it is an 
ultrapower of $(R/\mathfrak n)\tensor_k K$ by \eqref{e:primary}.
Since being 
Gorenstein is first order definable for Artinian local rings by 
\cite{SchEC}, we get that $R/\mathfrak n$ is Gorenstein \iff\ $\hull 
{R/\mathfrak n}$ is \iff\ almost all $\seq Rw/\seq{\mathfrak n}w$ are.
Since almost every $\frak{n}_w$ is generated by an $R_w$-sequence, this
is equivalent with $R_w$ Gorenstein for almost all $w$.
\end{remark}

\subsection{Hilbert-Samuel functions}
We now want to strengthen \eqref{i:trans} and show 
that almost all $R_w$ have the same Hilbert-Samuel
function as $R$.
For this, we assume that the reader is familiar with
the fundamentals of the theory of standard bases in power series rings;
for example, see \cite{Becker}. We fix $n>0$, and we denote by $\preceq$ 
the degree-lexicographic ordering on $\nat^n$, 
that is, $\nu\preceq\mu$ \iff\ $|\nu|<
|\mu|$, or $|\nu|=|\mu|$ and $\nu\leq\mu$ lexicographically. 
Let $L$ be a field. For every non-zero 
$$f=\sum_\nu a_\nu X^\nu\in L[[n]] 
\qquad\text{(with $a_\nu\in L$ for 
all $\nu\in \nat^n$)}$$
 there exists a $\preceq$-smallest 
$\lambda\in\nat^n$ with $a_\lambda\neq 0$, and we put $c(f):=a_\lambda$ and
$v(f):=\lambda$. It is convenient to define $c(0):=0$ and 
$v(0):=\infty$ with $\infty+
\nu=\nu+\infty=\infty$ for all $\nu\in\nat^n\cup\{\infty\}$. We extend
$\preceq$ to $\nat^n\cup\{\infty\}$ by $\nat^n\prec\infty$.
Note that $v$ is a valuation on $\pow Ln$ with values in the ordered
semigroup $(\mathbb N^{n},\preceq)$, that is, for all $f,g\in \pow Ln$: 
\begin{enumerate}
\item $v(f)=\infty \Longleftrightarrow
f=0$,
\item $v(fg)=v(f)+v(g)$, and 
\item $v(f+g) \succeq \min\big\{v(f),v(g)\big\}$. 
\end{enumerate}
Given a subset $\mathfrak{s}$ of $L[[n]]$ we put
$$v(\mathfrak{s}) := \bigl\{ v(f) : f\in \mathfrak{s}\pow Ln \bigr\} \subseteq\nat^n\cup\{\infty\}$$
where $\mathfrak{s}\pow Ln $ denotes the ideal generated by $\mathfrak{s}$. Let $f,g_1,\dots,g_m\in L[[n]]$.
We call an expression
$$f=\sum_{i=1}^m q_ig_i \qquad \text{(where $q_1,\dots,q_m\in L[[n]]$)}$$
such that
$v(f)\preceq v(q_i)+v(g_i)$ for all $i$ a
\emph{standard representation} of $f$ with respect to $\mathfrak{s}=\{g_1,\dots,g_m\}$
(in $L[[n]]$). Note that then
$v(f)$ equals  the ($\preceq$-) minimum of the
$v(q_i)+v(g_i)$.
%since $v$ is a
%valuation. 
If $L\subseteq L'$ is a field extension, and $f\in L[[n]]$ has a standard
representation with respect to $\mathfrak{s}$ in $L'[[n]]$, then
$f$ has a standard
representation with respect to $\mathfrak{s}$ in $L[[n]]$. (Since $L[[n]]\to L'[[n]]$
is faithfully flat.) Moreover:

\begin{lemma}\label{l:standard}
An element $f$ of $K[[n]]$ has a standard representation with respect to
a subset $\mathfrak{s}=\{g_1,\dots,g_m\}$ of $K[[n]]$ \iff\ almost every $f_w$ has
a standard representation with respect to
$\seq{\mathfrak{s}}w:=\{g_{1w},\dots,g_{mw}\}$,  where $\seq fw$ and $\seq{g_i}w$ are
 \sr{s} of $f$ and $g_{i}$ respectively.
\end{lemma}
\begin{proof}
We may assume $f\neq 0$. 
Writing $f=f_0+\varepsilon$ where 
$f_0\in K[n]$ is homogeneous of degree $d:=|v(f)|$ and
$\varepsilon\in\maxim^{d+1}$ we see that $c(f)_w=c(f_w)$ and
$v(f)=v(f_w)$ for almost all $w$. Hence if $f=\sum_{i=1}^m q_ig_i$ is a
standard representation of $f$ with respect to $\mathfrak{s}$, then
$f_w=\sum_{i=1}^m q_{iw}g_{iw}$ is a standard representation of $f_w$ in
terms of $\seq{\mathfrak{s}}w$, where $\seq{q_{i}}w$ is an \sr\ of $q_{i}$. 
Conversely, suppose that almost every $f_w$ has  
a standard representation $f_w=\sum_{i=1}^m q_{iw}g_{iw}$ with
respect to $\seq{\mathfrak{s}}w$, where $q_{iw}\in K_w[[n]]$. Since $v$ is a
valuation, there is some $i$ such that $v(\seq fw)=v(\seq{q_{i}}w)+
v(\seq{g_{i}}w)\preceq v(\seq{q_{j}}w)+ v(\seq{g_{j}}w)$ for all $j$ and
almost all $w$. Therefore, if we let $q_j$ be the
ultraproduct of the $q_{jw}$ and $\pi$ as in Proposition~\ref{P:infal}, then
$v(f)=v(\pi(q_i))+v(g_{i})\preceq v(\pi(q_j))+v(g_{j}))$ for all $j$, showing that
$f=\sum_{i=1}^m \pi(q_i)g_i$ is a standard representation of $f$ with
respect to $\mathfrak{s}$ in $\hull K[[n]]$. Hence $f$ has a standard representation
with respect to $\mathfrak{s}$ in $K[[n]]$ by faithful flatness.
\end{proof}

Every ideal $I$ of $L[[n]]$ has a
\emph{standard basis,} that is, a finite subset 
$\mathfrak{s}$ of $I$ such that every element of $I$
has a standard representation with respect to $\mathfrak{s}$, or equivalently, such that
$v(\mathfrak{s})=v(I)$. (See \cite[Theorem on p.~219]{Becker}.)

\begin{proposition}\label{P:standard}
A subset $\mathfrak{s}$ of an ideal $I\subseteq K[[n]]$ is
a standard basis for $I$ \iff\ its \sr\ $\seq{\mathfrak{s}}w$
is a 
standard basis for the \sr\ $I_w\subseteq K_w[[n]]$ of $I$, for almost all $w$.
In particular we have $v(I)=v(I_w)$ for almost all $w$.
\end{proposition}
\begin{proof}
We use the Buchberger criterion for standard bases:
for non-zero 
$f,g\in L[[n]]$ we define $$s(f,g) := c(g)X^\mu f-c(f)X^\nu g\in L[[n]]$$ 
where $\mu,\nu$ are the multiindices in $\nat^n$
such that $X^{\mu+v(f)}=X^{\nu+v(g)}=$ the least common multiple of
$X^{v(f)}$ and $X^{v(g)}$. Then a 
finite subset $\mathfrak{e}$ of $L[[n]]$ is a standard basis
of the ideal it generates \iff\ $s(f,g)$ has a standard representation with
respect to $\mathfrak{e}$, for all $0\neq f,g\in \mathfrak{e}$ \cite[Theorem~4.1]{Becker}.
The claim follows
from this and Lemma~\ref{l:standard}, since
if $f,g\in \mathfrak{s}$ are non-zero then their \sr{s}
$f_w,g_w$ are non-zero and $s(f_w,g_w)$ is an \sr\ of $s(f,g)$ for almost all
$w$.
\end{proof}

Given a Noetherian local ring $(S,\mathfrak n)$ we use $\chi_S$ to
denote the Hilbert-Samuel function $d\mapsto \operatorname{length}(S/\mathfrak{n}^{d+1})$ of $S$.
By Corollary~\ref{c:infini} we see that for fixed $d\in\nat$, we have
$\chi_{R}(d)=\chi_{R_w}(d)$ for almost all $w$. Proposition~\ref{P:standard}
implies the following stronger version:

\begin{corollary}
For almost all $w$, we have $\chi_R = \chi_{R_w}$ \textup{(}that is,
$\chi_R(d)=\chi_{R_w}(d)$ for all $d$\textup{)}.
\end{corollary}
\begin{proof}
Since $\complet R_{(k,u)}/\maxim^{d+1}\complet R_{(k,u)} = 
(R/\maxim^{d+1})\tensor_k K$ for all $d$, we have $\chi_{R}=
\chi_{\complet R_{(k,u)}}$.
For an ideal $I$ of $L[[n]]$,
the Hilbert-Samuel functions of $L[[n]]/I$ and $L[[n]]/I'$, 
where $I'$ is the
ideal generated by all $X^\nu$ with $\nu\in v(I)$, agree.
Hence by Proposition~\ref{P:standard} we obtain that 
$\chi_{\complet R_{(k,u)}}=\chi_{R_w}$ 
for almost all $w$.
\end{proof}

In particular, almost all $R_w$ have the same Hilbert-Samuel polynomial
as $R$, hence the same multiplicity, and we see once more that almost all
$R_w$ have the same dimension and the same embedding dimension as $R$.
For later use we also show:

\begin{lemma}\label{l:perturb}
Let $f_1,\dots,f_r\in K[[n]]$ and $\varepsilon_1,\dots,\varepsilon_r\in
\infal{\hull n}$, and
consider the ideals $I=(f_1,\dots,f_r)\pow Kn$  and
$I_\varepsilon=(f_1+\varepsilon_1,\dots,f_r+\varepsilon_r)\hull n$  
with respective \sr{s} $\seq Iw$ and $I_{\varepsilon,w}$.
Then $v(I_w) \subseteq v(I_{\varepsilon,w})$ and hence
$\dim(K_w[[n]]/I_w) \geq \dim(K_w[[n]]/I_{\varepsilon,w})$, for almost all $w$.
\end{lemma}
\begin{proof}
We may assume $r>0$ and $f_i\neq 0$ for all $i$. Let $d:=\max_i |v(f_i)|$.
Then for almost all $w$ we have $\varepsilon_{iw}\in\frak{m}_w^{d+1}$
and hence $v(f_{iw}+\varepsilon_{iw})=
v(f_{iw})=v(f)$ for almost all $w$.
Let $\mathfrak{s}=\{g_1,\dots,g_m\}$ be a standard basis for $I$. Then its
\sr\ 
$\seq{\mathfrak{s}}w$ is a standard basis for $I_w$ by
Proposition~\ref{P:standard},
and thus $v(I_w)=v(\mathfrak{s}_w)$, for almost all $w$. For every 
$j\in\{1,\dots,m\}$ there exists $i\in\{1,\dots,r\}$ and $\nu\in\nat^n$
with $v(g_j)=v(f_i)+\nu$. Hence $v(g_{jw})=v(f_{iw})+\nu=v(f_{iw}+\varepsilon_{iw})+\nu$ for almost all $w$.
This shows $v(I_w)=v(\mathfrak{s}_w)\subseteq v(I_{\varepsilon,w})$ for almost all $w$.
\end{proof}

\subsection{Irreducibility}\label{s:irr}
Suppose that $\Lambda=(R,\mathbf{x},k,u)$, and
recall from the discussion before 
\ref{C:nnorm} that $i=u|k$ is the embedding of $k$ into the algebraic closure
of $u(k)$ inside $K$.
We call $R$
\emph{absolutely analytically irredu\-cible} 
if $\complet R_{(k,i)}$ is a domain. This does not depend on the
choice of $k$ and $u$. (Cf.~Remark~\ref{r:iso}.)
\emph{From now on up to and including \S\textup{\ref{c:prod}} we assume that
$\hull R$ is absolutely normalizing.}

\begin{theorem}\label{T:dom}
%Suppose that $\hull R$ is absolutely normalizing. 
The following statements are equivalent:
\begin{enumerate}
\item\label{i:dom1} $R$ is 
absolutely analytically irreducible;
\item\label{i:dom2} $\hull R$ is a domain;
\item almost all $\seq Rw$ are domains.
\end{enumerate}
\end{theorem}

We first establish some auxiliary facts needed in the proof.
Let $T$ be a domain with fraction field $F=\Frac(T)$. Let
$Y=(Y_1,\dots,Y_m)$ be a tuple of indeterminates,
and let $I$ be a finitely generated 
ideal of $T[Y]$. There exists a non-zero
$\delta\in T$ with the following property: for all domains $T'$ extending $T$,
with fraction field $F'=\Frac(T')$, and all $f\in T'[Y]$ we have
$f\in IF'[Y]$ \iff\ $\delta f\in IT'[Y]$.
(See, e.g., \cite[Corollary~3.5]{Asch03}.)
In other words, $$IF'[Y]\cap T'[Y] = \bigl(IT'[Y]:_{T'[Y]}\delta\bigr)$$
and therefore:

\begin{lemma}\label{l:contract}
If $T'$ is a domain extending $T$, with fraction field $F'$, and
$T'$ is flat over $T$,
then $$IF'[Y]\cap T'[Y] = \bigl(IF[Y]\cap T[Y]\bigr) T'[Y].$$
\end{lemma}

%Note that $I$ is a prime ideal with $I\cap T=(0)$ \iff\ $IF[Y]$ is a prime
%ideal and $IF[Y]\cap T[Y]=I$. Moreover, assuming that 
%$\operatorname{char} F=0$, if $IF[Y]$ is radical, then so is 
%$IF'[Y]$, by \cite[???]{ZS}, and hence $IT'[Y]$ is radical.

%Suppose now that $k$ is an algebraically closed subfield of $K$, and
%$T=k[[n]]$.
%Let $\psi$ be the unique extension of the inclusion $k\to K$ to
%a homomorphism $T\to K$ with kernel $(X_1,\dots,X_n)$. Then
%$(T,k,\psi)$ is an object in $\mathbf{C}_K$, and
%$T^* := \hull {T,k,\psi} = \hull n$ is an integrally closed 
%domain (since each $T_w=K_w[[n]]$ is). We identify $T$ with its image
%under the faithfully flat embedding $\eta_{(T,k,\psi)}\colon T\to T^*$.
In the following proposition and lemma
let $T=K[[d]]$ and $T^*=\hull{d}$.

\begin{proposition}\label{P:expand}
If  $I$ is a prime ideal of $T[Y]$ with $I\cap T=(0)$, 
then $I T^*[Y]$ is a prime ideal of $T^*[Y]$ with $I T^*[Y]\cap T^*=(0)$.
\end{proposition}

For the proof we need:

\begin{lemma}\label{l:regular}
The fraction field $F^*$ of $T^*$ is a regular extension of $F$.
\end{lemma}
\begin{proof}
Since $\operatorname{char} F=0$ we only need to show that $F$ is
algebraically closed in $F^*$. Let $y\in F^*$ be algebraic over $F$.
To show that $y\in F$ we may assume that $y$ is integral over $T$.
Since $T^*$ is integrally closed  it follows that
$y\in T^*$.  Let $P(Y)\in T[Y]$ be a monic polynomial of minimal
degree such that $P(y)=0$. Then
$\pi(y)$ is a zero of $P$ in $\mathbf{K}[[d]]$,
where $\pi\colon T^*\to
\mathbf{K}[[d]]$ is the surjective $K[[d]]$-algebra 
homomorphism from Proposition~\ref{P:infal} and $\mathbf{K}=\mathfrak D(K)$. 
Since $K$ is algebraically
closed it follows (using Hensel's Lemma) that $P$ has a zero in $K[[d]]=T$.
By minimality of $P$, this zero is $y$, so $y\in T$ as required.
\end{proof}

\begin{proof}[Proof of Proposition~\ref{P:expand}]
Suppose that  $I$ is prime and $I\cap T=(0)$, or equivalently, $IF[Y]$ is
prime and $IF[Y]\cap T[Y]=I$. 
By Lemma~\ref{l:regular}, $IF^*[Y]$ is a prime ideal of $IF[Y]$.
(See \cite{BA}, Chapitre V, \S{}15, Proposition~5 and \S{}17, 
Corollaire to Proposition~1.) In particular
$I T^*[Y]\cap T^*=(0)$. Since $T^*$ is flat over $T$, by Lemma~\ref{l:contract}
we have
$$IF^*[Y]\cap T^*[Y] = \bigl(IF[Y]\cap T[Y]\bigr) T^*[Y] = IT^*[Y].$$
It follows that $IT^*[Y]$ is prime. 
\end{proof}

\begin{proof}[Proof of Theorem~\ref{T:dom}]
%Using the notation introduced in \S\ref{s:ten} we have 
%$(\complet R_{\alg k})_K \iso \complet R_K$ as $K$-algebras, hence
%$\hull{\complet R_{\alg k}} \iso \hull{R}$.  Therefore we may 
%replace $R$ by $\complet R_{\alg k}$ and assume  that $R$ is complete 
%with algebraically closed coefficient field $k$.
By \los, almost all $R_w$ are domains \iff\ $\hull R$ is.
Moreover, if 
this is the case, then every subring 
of the domain $\hull R$ is also a domain.
Hence we only have to prove that if $R$ is absolutely analytically irreducible, 
then $\hull R$ is a domain. 
Let us first assume that $\complet R_{(k,u)}$ is a domain.
Put $T:=K[[d]]$ and let $J:=T[X_{d+1},\dots,X_n]\cap\complet I_\Lambda$, where
$d=\dim R$. Since $\complet\theta_\Lambda$ is a Noether normalization
of $C(\Lambda)$,
$$\complet R_{(k,u)}\iso T[X_{d+1},\dots,X_n]/J$$ and
$$\hull{R} \iso \hull{T}[X_{d+1},\dots,X_n]/J\hull{T}[X_{d+1},\dots,X_n].$$
(See Proposition~\ref{P:finan} and the discussion following Proposition~\ref{P:Noeth}.)
Now $\complet R_{(k,u)}$ is a domain \iff\ $J$ is a prime ideal,
and in this case, by
Proposition~\ref{P:expand}, the expansion $J \hull T [X_{n+1},\dots,X_{n+m}]$
of $J$ to an ideal of $\hull T [X_{n+1},\dots,X_{n+m}]$ remains prime. 
Hence $\hull R$ is a domain, as
required. The proof of Theorem~\ref{T:dom} is now completed by Lemma~\ref{L:domain} below.
\end{proof}

\begin{lemma}\label{L:domain}
If $\complet R_{(k,i)}$ is an integral domain then so is $\complet R_{(k,u)}$.
\end{lemma}
\begin{proof}
We write $l$ for the algebraic closure of $u(k)$ inside $K$.
It is easy to see that the unique extension of a 
Noether normalization $l[[d]]\to 
\complet R_{(k,i)}$ of $\complet R_{(k,i)}$ to a $K$-algebra
homomorphism $K[[d]]\to\complet R_{(k,u)}$ is a Noether normalization of 
$\complet R_{(k,u)}$. Hence the argument above, which allowed us to
transfer integrality from $\complet R_{(k,u)}$ to $\hull R$, can be used to
transfer integrality of $\complet R_{(k,i)}$ to $\complet R_{(k,u)}$, 
provided we know that the fraction field of $K[[d]]$ is a
regular extension of the fraction field of $l[[d]]$. 
This is shown as in Lemma~\ref{l:regular}.
\end{proof}

A prime ideal $\pr$ of $R$ is called 
\emph{absolutely analytically prime} if $R/\pr$ is  absolutely 
analytically irreducible, that is to say, if $\pr\complet R_{(k,i)}$ is prime.
Since $\hull R$ is absolutely normalizing, so is
$\hull{\Lambda/I}=\hull{R}/I\hull{R}$
for every ideal $I$ of $R$.  Hence
the theorem implies:

\begin{corollary}\label{C:dom}
%Assume that $\hull R$ is absolutely normalizing. 
The following statements are equivalent, for a prime ideal $\pr$ of $R$:
\begin{enumerate}
\item $\pr$ is absolutely analytically prime;
\item $\pr\hull R$ is prime;
\item almost all approximations $\seq\pr w$ of $\pr$ are prime. \qed
\end{enumerate}
\end{corollary}

\subsection{Reducedness}\label{s:unr}
A local ring $A$ is called \emph{analytically 
unramified} (or, \emph{analytically reduced}), if its completion is 
reduced (that is to say, without 
non-zero nilpotent elements). 

\begin{theorem}\label{T:unr}
The following statements are equivalent:
%when $\hull R$ is absolutely normalizing:
\begin{enumerate}
\item\label{i:unr1} $R$ is analytically unramified;
\item\label{i:unr2} $\hull R$ is reduced;
\item\label{i:unr3} almost all $\seq Rw$ are reduced.
\end{enumerate}
\end{theorem}
\begin{proof}
The implication
 \eqref{i:unr3}~$\Rightarrow$~\eqref{i:unr2} is a consequence of \los, and
the implication \eqref{i:unr2}~$\Rightarrow$~\eqref{i:unr1} is trivial.
Hence we only need to show that if $\complet R$ is reduced, then almost all
$\seq Rw$ are reduced. If $\complet R$ is reduced, then so is $\complet 
R_{(k,i)}$, by Lemma~\ref{l:scalarext}. 
Let $\pr_1,\dots,\pr_s$ be the minimal prime ideals of $\complet R_{(k,i)}$. 
Since $\complet R_{(k,i)}$ is reduced, their intersection is zero,
and hence so is 
the intersection of their \sr{s} $\seq{\pr_i}w$ for almost all $w$. 
Since $\hull{R}$ is absolutely normalizing,
 almost all $\seq{\pr_i}w$ are prime ideals by Corollary~\ref{C:dom}, 
proving that almost all $\seq Rw$ are reduced. 
\end{proof}

\begin{corollary}\label{c:radical}
Suppose that $R$ is excellent. 
For an ideal $I$ of $R$
the following are equivalent:
\begin{enumerate}
\item $I$ is radical;
\item $I\hull R$ is radical;
\item almost all approximations $I_w$ of $I$ are radical.
\end{enumerate}
In particular, we have $\sqrt{I}\hull R = \sqrt{I\hull R}$, and
$(\sqrt{I})_w = \sqrt{I_w}$ for almost all $w$. \qed
\end{corollary}

A Noetherian ring is called \emph{equidimensional} if all its minimal primes
have the same dimension.
A Noetherian local ring is called \emph{formally equidimensional}
if its completion is equidimensional.

\begin{corollary}\label{c:equidim}
%Suppose $\hull R$ is absolutely normalizing. 
If $R$ is complete and 
$k$ is algebraically closed, then the following are equivalent, for
a prime ideal $\pr$ of $R$:
\begin{enumerate}
\item\label{i:min1} $\pr$ is a minimal prime ideal of $R$;
\item\label{i:min2} $\pr\hull R$ is a minimal prime ideal of $\hull R$;
\item\label{i:min3} for almost all $w$ the approximation $\pr_w$ of $\pr$ is a 
minimal prime ideal of $R_w$.
\end{enumerate}
%In particular, $\hull R$ has only finitely many minimal prime ideals.
If $R$ is arbitrary, then  $R$ is formally equidimensional \iff\ almost all $R_w$
are  equidimensional.
\end{corollary}
\begin{proof}
The intersection of the minimal prime ideals
$\pr_1,\dots,\pr_s$ of $R$ equals the (nil-) radical of $R$.
By Corollary~\ref{C:dom} the $\pr_{i}\hull R$ are prime ideals of $\hull R$,
and almost all $\pr_{iw}$ are prime ideals of $R_w$.
By the previous corollary
and \eqref{i:cap}, the intersection
$\pr_1\hull R\cap\cdots\cap\pr_s\hull R$ equals the radical of $\hull R$,
and hence $\pr_{1w}\cap\cdots\cap\pr_{sw}$ is the radical of $R_w$ for
almost all $w$. 
This yields the equivalence of \eqref{i:min1}--\eqref{i:min3}.
It remains to show that when $R$ is arbitrary, it is formally equidimensional
\iff\ almost  all $R_w$ are equidimensional.
Using Lemma~\ref{l:scalarext} we reduce to the case that
$R$ is complete and $k$ is algebraically closed, and then the 
claim follows from the earlier statements and 
\eqref{i:dim}.
\end{proof}

Given a ring $A$ and ideals $\frak{a}_1,\dots,\frak{a}_s$ of $A$,
the canonical homomorphism $$A\to A/\frak{a}_1\times\cdots\times A/\frak{a}_s$$ is
an isomorphism \iff\ $\frak{a}_1\cap\cdots\cap \frak{a}_s=(0)$ and 
$1\in \frak{a}_i+\frak{a}_j$ 
for all $i\neq j$.
Hence by Theorem~\ref{T:unr} and Corollary~\ref{c:equidim} we get: 

\begin{corollary}\label{c:prod}
%Assume that $\hull R$ is absolutely normalizing. 
Suppose that $R$ is complete and $k$ is algebraically closed. Let
$\pr_1,\dots,\pr_s$ be the minimal prime ideals of $R$. The following
statements are equivalent:
\begin{enumerate}
\item the canonical homomorphism $R\to R/\pr_1\times\cdots\times R/\pr_s$
is bijective;
\item  the canonical homomorphism $$\hull{R} \to \hull{R}/\pr_1\hull{R}\times\cdots\times
\hull{R}/\pr_s\hull{R}$$ is bijective;
\item  the canonical homomorphism 
$$R_w\to R_w/\pr_{1w}\times\cdots\times R_w/\pr_{sw}$$ is bijective for
almost all $w$. \qed
\end{enumerate}
\end{corollary}

\begin{remark}\label{R:unr}
The proof of Theorem~\ref{T:dom} shows that if $\complet\theta_\Lambda$ 
is normalizing with respect to a prime ideal $\pr$ of $\complet R_{(k,u)}$, 
then almost all approximations of $\pr$ are prime. Hence by
the proof of Theorem~\ref{T:unr}: if 
$\complet\theta_\Lambda$ is normalizing 
for all minimal primes of an ideal $I$ of
$\complet R_{(k,u)}$, then $(\sqrt{I})_w=\sqrt{I_w}$ for almost all $w$.
(This will be used in \S\ref{s:normal} below.)
\end{remark}

\begin{remark}\label{R:unr2}
Suppose that $\Lambda$ is normalizing (see \S\ref{C:nnorm2}). 
Then Theorem~\ref{T:dom} above remains
true, with the same proof. Moreover, if $k$ is algebraically closed and
$\pr_1,\dots,\pr_s$ are
the minimal primes of $\complet R$, then
for almost all $w$, the approximations $\pr_{1w},\dots,\pr_{sw}$ are
the minimal primes of $R_w$, and Theorem~\ref{T:unr} and 
Corollary~\ref{c:equidim} also remain true.
(This will be used in Sections~6 and 7.)
\end{remark}

\subsection{Normality}\label{s:normal}
Recall that a domain is called \emph{normal} if it 
is integrally closed in its fraction field. 
By \los, $R$ is a normal domain \iff\ 
almost all $R_w$ are normal domains, and in this case $R$ is a normal
domain, by faithful
flatness of $R\to\hull R$.

\begin{theorem}\label{T:normal}
Suppose that $R$ is a complete normal domain with algebraically closed
residue field. Then
$\Lambda$ with underlying ring $R$ can be chosen such that $\hull{R}=
\hull\Lambda$ is a normal domain. 
\end{theorem}

The proof 
is based on the following criterion for normality due to
Grauert and Remmert \cite[pp.~220--221]{GR}; see also
\cite{deJong}. Let $B$ be a Noetherian domain, and 
$N(B)$ be the non-normal locus of $B$, that is, the set of all prime ideals
$\pr$ of $B$ such that $B_{\frak p}$ is not normal.

\begin{proposition}\label{P:GR}
Let $H$ be a non-zero radical ideal of $B$
such that every $\pr\in N(B)$ contains $H$, and $0\neq f\in H$. Then
$$\text{$B$ is normal} \qquad\Longleftrightarrow\qquad
\text{$fB=(fH:_BH)$.}$$
\end{proposition}

Let $A$ be a ring and $B$ an $A$-algebra of finite type,
that is, $B$ is of the form
$B=A[Y]/J$
where $J=(f_1,\dots,f_r)\pol AY$ is an ideal of the polynomial ring 
$A[Y]=A[Y_1,\dots,Y_m]$. Given a tuple $\mathbf g=\rij gs$ with entries in
$\{f_1,\dots,f_r\}$ we write $\Delta \mathbf g$ for the ideal of $A[Y]$ generated by
all the $s\times s$-minors of the $s\times m$-matrix $\left(\frac{\partial g_i}{\partial Y_j}\right)$, with the understanding that 
$\Delta\emptyset:=A$.
We let $H_{B/A}$ denote the nilradical of the ideal
in $A[Y]$ generated by $J$ and by the $\Delta \mathbf g\cdot \bigl(\mathbf g\pol AY:J\bigr)$, for $\mathbf g$ ranging over all tuples with entries in $\{f_1,\dots,f_r\}$. The image in $B$ of the ideal $H_{B/A}$
does not depend on the chosen presentation $B\iso A[Y]/J$
of the $A$-algebra $B$. (See \cite[Property~2.13]{Spi}.)
If $A$ is Noetherian and $\frak p\supseteq J$ a prime ideal of $A[Y]$, then
$B_{\frak p}$ is smooth over $A$ \iff\ $H_{B/A} \not\subseteq \frak p$.
In this case, $A\to B_{\frak p}$ is regular
\cite[Corollary~2.9]{Spi}. In particular, if $A$ is regular, then so is 
$B_{\frak p}$ \cite[Theorem~23.7]{Mats} and, since a regular ring is normal, the canonical image of $H_{B/A}$
in $B$ is then a 
non-zero radical ideal
which is contained in every element of $N(B)$. Therefore
Proposition~\ref{P:GR} implies:

\begin{corollary}\label{C:jacnorm}
Let $B$ be an integral domain, of finite type over a regular ring $A$, and let
$f$ be a non-zero element of the canonical image $H$ of $H_{B/A}$ in $B$. 
Then $B$ is normal \iff\ $fB=(fH:_B H)$. \qed
\end{corollary}

\begin{proof}[Proof of Theorem~\ref{T:normal}]
The desired
object $\Lambda$ has the form $(R,\mathbf x,k,u)$, where $k$ is an arbitary
coefficient field of $R$,  where $u\colon R\to K$ is an arbitrary local
homomorphism, and $\mathbf x$ is
determined as follows.
Choose a Noether normalization $\theta\colon K[[n]]\to \complet R_{(k,u)}$ 
of $\complet R_{(k,u)}$. 
Let $I:=\ker\theta$ and put
$J:=A[X_{d+1},\dots,X_n]\cap I$, where $d:=\dim R$ and $A:=K[[d]]$. 
Let
$B:=A[X_{d+1},\dots,X_n]/J$, so that $B \iso \complet R_{(k,u)}$ as $A$-algebras, 
and
let $H$ be the image of the ideal $H_{B/A}$ of $A[X_{d+1},\dots,X_n]$ in 
$\complet R_{(k,u)}$. 
We already remarked that $H$ does not depend on the choice
of $\theta$.
By Lemma~\ref{L:nnorm} we can choose $\theta$
normalizing for all minimal prime ideals of $H$ and for all ideals 
$\id\complet R_{(k,u)}$, where $\id$ is an ideal of $\complet R_{(k,i)}$. 
Now put
$x_i := \theta(X_i)$ for 
$i=1,\dots,n$ and $\mathbf x:=\rij xn$. It follows that
$\theta=\complet\theta_\Lambda$ (hence $I=\complet I_\Lambda$) for the thus constructed $\Lambda$,
and  $\Lambda$ is absolutely normalizing.

We claim that $\hull R=\hull\Lambda$ is a normal domain.
By Lemma~\ref{L:domain} and \cite[Theorem~23.9]{Mats},
$\complet R_{(k,u)}$ is a normal domain, and therefore,  
$\hull{R}$ is a domain, by Theorem~\ref{T:dom}.
With $\seq Aw:=\pow{\seq Kw}d$,  Proposition~\ref{P:Noeth} yields for almost all $w$ an isomorphism of $A_w$-algebras
$$R_w \iso B_w := A_w[X_{d+1},\dots,X_n]/J_w$$
where $\seq Jw$ is an \sr\ of the ideal $J$ of $A[X_{d+1},\dots,X_n]$. 
%We identify $\complet R_{(k,u)}$ with $B$ and $R_w$ with $B_w$.
We claim that $H_{B_w/A_w}$ is an approximation of $H_{B/A}$, 
for almost all $w$.
This implies that for almost all $w$ the canonical image of $H_{B_w/A_w}$ in
$B_w$ is an \sr\ of $H$. Lemma~\ref{L:ff} and Corollary~\ref{C:jacnorm} then show that almost all $B_w$ are normal, as required.

To establish the
claim, note that since a radical ideal of $A[X_{d+1},\dots,X_{n}]$ remains
radical upon extension to $K[[n]]$, the ideal
$H_{B/A}K[[n]]$ is the radical of the ideal
$$I+\sum_{\mathbf g} \Delta {\mathbf g}\cdot
\bigl(\mathbf gK[[n]]:_{K[[n]]} I\bigr)$$
where $\mathbf g$ ranges over all tuples with entries in a fixed set
of generators of $J$, and 
similarly the ideal
$H_{B_w/A_w}K_w[[n]]$ 
is the radical of
$$\seq Iw+\sum_{\mathbf g} \Delta \mathbf g_w\cdot
\bigl(\mathbf g_wK_w[[n]]:_{K_w[[n]]} I_w\bigr)$$
where $\seq Iw$ and $\mathbf g_w$ are \sr{s} of $I$ and $\mathbf g$  respectively.
Note that the ideal $\Delta \mathbf g_w$ of $A_w[X_{d+1},\dots,X_n]$
is an \sr\ of the ideal $\Delta\mathbf g$ of $A[X_{d+1},\dots,X_n]$.
It follows that $H_{B_w/A_w}K_w[[n]]$ is an \sr\ of $H_{B/A}K[[n]]$, for
almost all $w$, by
Lemma~\ref{L:ff}, Remark~\ref{R:unr}, and the choice of $\theta$.
Moreover  $$H_{B/A}=H_{B/A}\pow Kn\cap \pol A{X_{d+1},\dots,X_n}$$ 
and similarly, using the remarks preceding \S\ref{r:nn}:
$$H_{B_w/A_w}=H_{B_w/A_w}\pow {K_w}n\cap \pol {A_w}{X_{d+1},\dots,X_n}.$$ 
This yields that $H_{B_w/A_w}$ is an \sr\ of $H_{B/A}$, as claimed. 
\end{proof}

A ring $A$ is called \emph{normal} if $A_{\frak p}$ is a normal domain for
every prime ideal $\frak p$ of $A$.
If $A$ has finitely many minimal prime ideals $\pr_1,\dots,\pr_s$ then $A$ is
normal \iff\ $A\iso A/\pr_1\times\cdots\times A/\pr_s$ and each
domain $A/\pr_i$ is normal.
A local ring $A$ is called \emph{analytically 
normal} if $\complet A$ is normal.
%If $A$ is an excellent local ring, then $A$ is normal \iff\
%$A$ is analytically normal \cite[Theorem~32.2]{Mats}.

\begin{corollary}\label{c:normal}
Suppose that $R$ is analytically normal.
Then the object $\Lambda$ with underlying ring $R$ can be chosen such that
$\hull R=\hull{\Lambda}$ is normal and
almost all $R_w$ are normal. 
\end{corollary}
\begin{proof}
As in the proof of Theorem~\ref{T:unr} reduce to the case
that $R$ is complete and $k$ is algebraically closed. 
The claim now follows from 
Corollary~\ref{c:prod} and Theorem~\ref{T:normal}.
\end{proof}

Recall that Serre's condition $(\op R_i)$ for a Noetherian ring 
$A$ signifies that $A_\pr$ is regular for all 
prime 
ideals $\pr$ of $A$ of height at most $i$, see \cite[\S{}23]{Mats}.
In the transfer of property $(\op R_i)$, the fact that we do not
know whether $\eta_n$ 
commutes with partial differentiation (see Remark~\ref{r:partial})
poses a technical difficulty. We confine ourselves to showing:

\begin{theorem} \label{T:Ri}
Suppose that $R$ is equidimensional and excellent and $\hull R$ is absolutely normalizing. Then for
each $i$, if $R$ satisfies $(\op R_i)$ then so do
almost all $\seq Rw$.
\end{theorem}

To show this note that if $R$ is excellent, then $R$ 
satisfies $(\op R_i)$ \iff\ $S:=\complet R_{(k,u)}$ does 
(see Remark~\ref{r:scalarext}). 
By Corollary~\ref{c:equidim}, if $R$ is equidimensional and
$\hull R$ is absolutely normalizing, then 
almost all approximations $R_w$ 
of $R$ are equidimensional. Note that the $R_w$ are also approximations of
$S$. Now apply the following lemma to $S$:

\begin{lemma}
Suppose that $R\in\ancat$, and $R$  and almost all 
its approximations $R_w$ are
equidimensional.
Then for
each $i$, if $R$ satisfies $(\op R_i)$ then so do
almost all $\seq Rw$.
\end{lemma}

\begin{proof}
Let $f_1,\dots,f_r\in K[[n]]$ be generators of the ideal 
$I:=\complet I_\Lambda$, and
let $h$ be the height of $I$.
Let $J$ be the Jacobian ideal of $I$, that is to say 
the ideal of
$K[[n]]$ generated by $I$ and all $h\times h$-minors of the matrix with
entries $\partial f_i/\partial X_j$. By the Jacobian criterion for regularity
for power series rings in characteristc zero \cite[Theorem~30.8]{Mats}, 
given a prime ideal $\pr$ of $R$, 
the localization of $R$ at
$\pr$ is regular \iff\ $JR\not\subseteq\pr$.
Hence $R$ satisfies $(\op R_i)$ \iff\ $JR$ 
has height at least $i+1$. Since $R$ 
is equidimensional, this is equivalent with $J$ 
having height at least $h+i+1$, 
and hence with $K[[n]]/J$ having dimension at most $n-(h+i+1)$. By
\eqref{i:dim} this is in turn is equivalent with 
$\dim K_w[[n]]/J_w \leq n-(h+i+1)$ for almost all $w$, where $\seq Jw$ is an \sr\ of $J$. 
Now for every $w$ let $\seq{\widetilde J}w$ be the Jacobian ideal of $I_w$. 
By Remark~\ref{r:partial} and Lemma~\ref{l:perturb} we have
$\dim K_w[[n]]/\seq{\widetilde J}w \leq \dim K_w[[n]]/J_w$ for almost all $w$.
Hence if $R$ satisfies $(\op R_i)$, then almost all $\seq{\widetilde J}w$ have height
$\geq h+i+1$ and thus, since almost all $R_w$ are equidimensional,
almost all $\seq{\widetilde J}wR_w$ have height $\geq i+1$. Hence 
by the Jacobian criterion for
regularity for power series rings over the algebraically closed fields $K_w$
of positive characteristic \cite[Theorem~30.10]{Mats}, almost all
$R_w$ satisfy $(\op R_i)$. 
\end{proof}

\subsection{Affine \sr{s} and localization}\label{s:aff}
One of the main drawbacks of the present theory is the fact that 
there is no a priori way to compare the $\mathfrak D$-extension of a 
local   ring with the $\mathfrak D$-extension of one of its 
localizations.  For example, 
suppose that $R$ is complete and $k$ algebraically closed,
and let $\pr$ be a prime ideal of $R$. 
From Theorem~\ref{T:dom}, we know that we can choose $\hull{R}$ such that 
$\pr\hull R$ is  a prime ideal, and then $\hull R_{\pr\hull R}$ is a 
faithfully flat extension of $R_\pr$. However, it is not 
clear how this compares with a Lefschetz extension $\hull{R_\pr}$ of $R_\pr$:
there is no 
obvious \homo\ from $\hull R$ to $\hull{R_\pr}$, since the \homo\ 
$R\to R_\pr$ is not local. (This problem is already apparent in the simplest 
possible situation that $R=k[[n]]$ with $n>1$, and $\pr$ is generated by a 
single variable.) 

We have to take these considerations into account when comparing 
the affine \sr{}s defined in \cite{SchNSTC} with the present version of 
\sr{}s. Therefore, we restrict our attention to the case that 
$R=C_\maxim$ is a localization of a finitely generated $k$-algebra 
$C$ at a maximal ideal $\maxim$. Here $k$ is a Lefschetz
field, realized as an ultraproduct of algebraically closed subfields 
$\seq kp$ of $\seq Kp$, 
with respect to the same ultraset as used for $K$. We consider $k$
as a subfield of $K$ in the natural way.
Suppose  $C=\pol kn/I$  where  $I$ is an ideal of $\pol kn$. 
As explained in the 
introduction, the non-standard hull of $C$ is $$\ul C:=\ul{\pol 
kn}/I\ul{\pol kn}$$ where $\ul{\pol kn}$ is the ultraproduct of the 
$\pol{\seq kp}n$. By \cite[Corollary 4.2]{SchNSTC}, the 
ideal $\maxim\ul C$ is again prime and by definition 
\cite[\S4.3]{SchNSTC}, the non-standard hull of $R$ is then $$\ul 
R:=(\ul C)_{\maxim\ul C}.$$
If $C'=k[n']/I'$ is another $k$-algebra and $\maxim'$ a maximal ideal
of $C'$ such that $R':=C'_{\maxim'}$ is isomorphic to $R$ as $k$-algebras,
then $R_\infty \iso (R')_\infty$ as Lefschetz rings \cite[\S{}4.3]{SchNSTC}.
In particular, 
since $k$ is algebraically closed we can make a translation and assume that
$\maxim = (X_1,\dots,X_n)\pol kn$.
The embedding $\ul{\pol kn}\subseteq \ul{\pow kn}$  factors through 
$(\ul{\pol kn})_{\maxim\ul{\pol kn}}$, where we denote  the 
ultraproduct   of the $\pow{\seq kp}n$ by $\ul{\pow kn}$. Composing 
with the natural embedding $\ul{\pow kn}\subseteq\ul{\pow Kn}$ followed 
by the diagonal embedding $\ul{\pow Kn}\to
\ul{\pow Kn}^{\mathcal U}= \hull n$ yields 
a $\pol kn$-algebra \homo\
        \begin{equation*}
        \bigl(\ul{\pol kn}\bigr)_{\maxim\ul{\pol kn}}\to \hull n.
        \end{equation*}
Taking reduction modulo $I$ gives a \homo\ $\ul R\to \hull 
R$ making
\[
\begin{aligned}
%\label{e:hull}
\xymatrix@R+2em@C+3em{ {R} \ar[d] \ar[r]^{\eta_R} &
{\hull{R}} \ar@{<-}[dl]\\
{\ul R} }
\end{aligned}
\]
commutative, where $R\to \ul R$ is the canonical embedding.  
Let $\seqaff Rp$ be \sr{}s of $R$ in the affine sense, that is to 
say, $p$ ranges over the set of prime numbers and   the 
ultraproduct of the $\seqaff Rp$ is equal to $\ul R$. Recall that for almost 
all $p$, we can obtain $\seqaff
Rp$  as the localization of $\pol{\seq kp}n/\seqaff Ip$ at the
prime ideal $\seqaff\maxim p$, where $\seqaff Ip$ and $\seqaff\maxim p$ are
respective \sr{s} of $I$ and $\maxim$ in the sense of \cite{SchNSTC}. Let $p(w):=\op{char} \seq
Kw$, so $\seq k{p(w)}$ is a subfield of $\seq Kw$, for each $w$. Let
$\seq Rw$ be the completion of $\seqaff R{p(w)}\tensor_{\seq
k{p(w)}}\seq Kw$ at the ideal generated by
the $X_i$. Hence there is a canonical map 
$\seqaff R{p(w)}\to \seq Rw$ and this is faithfully flat. Alternatively, with 
the notation from \S\ref{s:extsc}, we have that
 $$
\seq Rw= (\complet{\seqaff R{p(w)}})_{(\seq k{p(w)},\seq uw)}
$$
where $\seq uw\colon \seqaff R{p(w)}\to
\seq Kw$ is the composition of the residue map $\seqaff R{p(w)}\to \seq k{p(w)}$
with
the inclusion $\seq k{p(w)}\subseteq\seq Kw$. It follows that the ultraproduct
of the $\seq Rw$ is equal to $\hull R$, 
showing that the $\seq Rw$  are \sr{}s of $R$ in the present sense. 
Moreover, if $c\in R$, then \sr{}s $\seq cw$ of $c$ in the present 
sense are obtained by taking \sr{}s $\seqaff cp$ of $c$ in the sense of 
\cite{SchNSTC} and setting $\seq cw:=\seqaff c{p(w)}$ (as an element of $\seq
Rw$).
Put succinctly,
 an \sr{} of $R$ is obtained by the process of taking an \sr\ of $R$ in 
the sense of \cite{SchNSTC}, extending scalars and completing. We  use this
below to
 compare results between the affine and the complete case.

\begin{proposition}\label{P:pure}
With the notations just introduced, the homomorphism $\ul R\to\hull R$ is  
pure, and it is flat if $R$ has dimension at most $2$.
\end{proposition}

A homomorphism $M\to N$ between modules over a ring $A$ is
\emph{pure} if it is injective (so $M$ can be regarded as a submodule of $N$)
and every finite system of linear equations
with constants in $M$ which admits a solution in $N$ admits a solution in $M$.
For a  module $M$ over a ring $A$ let $\mu(M)\in\mathbb N\cup\{\infty\}$ 
be the least number
of elements in a generating set for $M$, and put $$\mu_A(m):=
\sup\bigl\{\mu(\ker\varphi): \varphi\in\operatorname{Hom}_A(A^m,A)\bigr\}
\in\mathbb N\cup\{\infty\}\qquad\text{for all $m$.}$$ The ring $A$ is called
\emph{uniformly coherent} if $\mu_A(m)<\infty$ for all $m$.
If $A$ is a
finitely generated algebra over a field then
$A$ is uniformly coherent \iff\ $\dim A\leq 2$, and
in this case $\mu_A(m)\leq m+2$ for all $m$. (See
\cite[Corollary~6.1.21]{Glaz}.) 

\begin{lemma}\label{L:flatcoh}
For each $v$ in an ultraset $\mathcal V$ , let $\seq Cv\to \seq Dv$ be a flat
\homo, with each $\seq Cv$ a  two-dimensional algebra  over a field, and let
$\ul C\to\ul D$ be their ultraproduct. If $\ul D\to D^*$ is any elementary map,
then the composition $\ul C\to \ul D\to D^*$ is flat.
\end{lemma}
\begin{proof}
 We have to show that for every linear form $L\in\pol{\ul C}Y$ where $Y=\rij
Ym$, 
the solution set of $L=0$ in $(D^*)^m$ is generated by the solution set of 
$L=0$ in $(\ul C)^m$. Let $\seq Lv$ be an \sr\ of $L$. For each $v$, there exist
$m+2$ tuples $\seq{\mathbf a_1}v,\dots,\seq{\mathbf
a_{m+2,}}v$ with entries in $\seq Cv$ which generate the  solution set of $\seq
Lv=0$ in $(\seq Cv)^m$, by uniform coherence. These same tuples generate the solution set
of $\seq Lv=0$ in $(\seq Dv)^m$, by flatness. 
The ultraproducts $\ul{\mathbf
a_1},\dots,\ul{\mathbf a_{m+2,}}$ of these $m+2$ tuples are then solutions of
$L=0$ and generate the solution set of $L=0$ in $(\ul D)^m$, by \los. 
Since $\ul D\to
D^*$ is elementary, the $\ul{\mathbf a_i}$ also generate the solution set of
$L=0$ in $(D^*)^m$, as required.
\end{proof}

\begin{proof}[Proof of Proposition~\ref{P:pure}] 
%It suffices to show that, if $\dim
%C\leq 2$, where we keep the notation introduced above. 
We keep the notation from above. The inclusions
\begin{equation}
\label{e:pure2}
\seqaff Rp  \to K_p[[n]]/\seqaff IpK_p[[n]]
\end{equation}
are faithfully flat and hence pure. Their ultraproduct $\ul R\to
K[[n]]_\infty/IK[[n]]_\infty$ is also pure. 
%Here the $\seqaff Ip$ and $\seqaff Cp$ are approximations of   $I$ and $C$ in
%the  sense of \cite{SchNSTC}. 
The diagonal embedding 
\begin{equation}
\label{e:pure1} K[[n]]_\infty/IK[[n]]_\infty\to \bigl(\ul{\pow Kn}/I\ul{\pow
Kn}\bigr)^{\mathcal U}= \hull R
\end{equation}
is pure, hence so is the composition
\begin{equation}
\label{e:pure3}
\ul R\to\hull R.
\end{equation}
Assume next that $R$ has dimension at most $2$. We may choose a finitely
generated $k$-algebra $C$
such
 that $R=C_\maxim$ with $C$ of dimension at most $2$. It suffices to show that
 the composition $\ul C \to \ul R\to \hull R$ is  flat. Almost all
$\seqaff Cp$ have dimension at most $2$ by
\cite[Theorem~4.5]{SchNSTC}.
%, whence
%$\mu_{\seqaff Cp}(n)\leq n+2$ for all $n$ and almost all $p$. It follows that 
%$\ul
%C$ is uniformly coherent with $\mu_{\ul C}(n)\leq n+2$ for all $n$.  
Since
\eqref{e:pure1} is   elementary, and since $\seqaff Cp\to \seqaff  Rp$   and  
\eqref{e:pure2} are flat, Lemma~\ref{L:flatcoh} yields that
\eqref{e:pure3} is flat, as required.
\end{proof}

We do not know in general whether $\ul R\to \hull R$ is flat (and hence
faithfully flat).

\subsection{The non-local case}\label{S:nonlocal} Let 
$A\supseteq\mathbb Q$ be a
Noetherian ring 
of cardinality at most the cardinality of $K$. Let $\op{Max} A$ be the set of
all maximal ideals of $A$, and for every $\frak n\in\op{Max} A$ choose a
faithfully flat Lefschetz extension $\eta_{A_{\frak n}}\colon A_{\frak
n}\to\hull{A_{\frak n}}$ of the Noetherian local ring $A_{\frak n}$ of
equicharacteristic zero. The product
of the $\eta_{A_{\frak n}}$ yields a faithfully flat embedding
\begin{equation}\label{e:ff-emb} 
A\to A^*:=\prod_{\frak n\in\op{Max} A} \hull{A_{\frak n}}.
\end{equation}
In general, $A^*$ is not   a Lefschetz ring, but it is so if $A$ is semi-local.
Thus:

\begin{proposition} \label{p:semilocal}
Every semi-local Noetherian ring containing $\mathbb Q$
admits a faithfully flat Lefschetz extension. \qed
\end{proposition}

For arbitrary $A$, in spite of the fact that $A^*$ is not Lefschetz, 
it still admits a non-standard
Frobenius, so that  the constructions   in the next two
sections can be generalized to the non-local case as well; see
\S\ref{s:tcnl} for a further discussion.

\part{Applications}

The standing assumption for the rest of this paper is that 
 $(R,\maxim)$ is an
equicharacteristic zero Noetherian local ring and
$K$ 
is an algebraically closed Lefschetz field with respect to an ultraset
whose underlying set is the set of all prime numbers, whose approximations
$K_p$ are algebraically closed fields of characteristic $p$ (as in Section~\ref{s:app}).
We take $K$ of  uncountable
cardinality,  as large as necessary. (Most of the time $\norm K>2^{\norm R}$
will suffice.)
We fix a
Lefschetz extension $\hull R$
of $R$ as defined in Part~1, and we let $(\seq Rw,\seq\maxim w)$
be the corresponding
\sr\ of $R$. In
other words, we fix some $\Lambda=(R,\mathbf{x},k,u)$ in
$\extcohcat$ with underlying ring $R$ and put $\hull R:=\hull \Lambda$. Where
necessary, we'll make some additional assumptions on $\Lambda$ (for instance
so that $\hull R$ is absolutely normalizing; see \S\ref{s:NN}). If $\alpha\colon
R\to S$ is a local \homo, then we choose an object
$\Gamma$ of $\extcohcat$ so that $\alpha$ induces a
morphism $\Lambda\to\Gamma$, and hence a local homo $\hull\alpha\colon \hull R\to \hull
S$. In the sequel, we often will  use a subscript $w$ to indicate a choice of 
\sr\ of a certain object without explicitly mentioning this. For instance, $\seq Sw$
will stand for some  \sr\ of $S$, etc.
%In the next two sections 
We now discuss non-standard tight closure and big \CM\  algebras, and indicate
several applications of these notions. 
%The 
%final section contains some model-theoretic consequences of the existence of
%$\hull{\,\cdot\,}$.

\section{Non-standard Tight Closure}\label{s:tc}

Every Lefschetz ring comes with a canonical endomorphism obtained by 
taking the ultraproduct of the Frobenii on each component:
let $\frob w\colon\seq Rw\to \seq Rw$ be the Frobenius $x\mapsto
x^{p(w)}$ on $\seq Rw$, where $p(w)$ denotes the \ch\ of $\seq Rw$, and let
$\ulfrob$ be the ultraproduct of the $\frob w$, that is to say, 
      \begin{equation*}
        \ulfrob\colon \hull R\to\hull R\colon  \up w \seq aw 
\mapsto \up w \frob w(\seq aw).
        \end{equation*} We call $\ulfrob$ the \emph{non-standard Frobenius} 
on $\hull R$.     More generally,  
for each $w$ let $\seq lw$ be a positive 
integer and let $\ul l$ be its ultraproduct in the ultrapower 
$\zet^{\mathcal W}$ of $\zet$. We let $\ulfrob^{\ul l}$ denote the ultraproduct
of the $\frob w^{\seq lw}$,    and call it an \emph{ultra-Frobenius} on $\hull
R$. In this paper, we are only concerned with the 
(powers of the) non-standard Frobenius $\ulfrob$; for an application 
of ultra-Frobenii, see \cite{SchLogTerm}. Note that if $\alpha\colon R\to S$ is
a local homomorphism, then for each $\ul l$, we have  a commutative diagram $$
\xymatrix@R+2em@C+3em{ \hull R \ar[r]^{\hull\alpha} \ar[d]^{\ulfrob^{\ul l} } &
\hull S \ar[d]^{\ulfrob^{\ul l}}\\
\hull R \ar[r]^{\hull\alpha} & \hull S.}$$
%This justifies our practice of omitting reference to the ring on which
%$\frob{\ul l}$ acts.
Given an ideal $\id$ of $R$, we 
use $\ulfrob(\id)\hull R$ to denote the ideal in 
$\hull R$ generated by all $\ulfrob(a)$ with $a\in\id$ 
(and a similar notation for powers of $\ulfrob$). Note 
that in general, $\ulfrob$ does not leave the subring $R$ invariant. 
In fact, we have an inclusion
        \begin{equation}
\label{eq:inf}
        \ulfrob (\maxim) \hull R\subseteq \infal{\hull R}. 
        \end{equation}
It follows that $\ulfrob(\maxim)\hull R\cap R=(0)$, by the faithful 
flatness of $R\to\hull R$.

Below we  often make use of the important fact (easily checked using \los\ 
and \cite[Theorem~16.1]{Mats}) that the image of a $\hull R$-regular sequence in
$\hull R$ under $\ulfrob$, and 
hence under each of its powers $\ulfrob^m$, is $\hull{R}$-regular. In
particular, by Lemma~\ref{L:sop}, the image under $\ulfrob^m$ of any $R$-regular
sequence in $R$ is $\hull R$-regular.

\subsection{Non-standard tight closure}\label{s:nstc}
 Let $\id$ be an ideal of
$R$.  We say that $z\in R$  belongs to the 
\emph{non-standard tight closure} of $\id$ if there exists $c\in R$ 
not contained in any minimal prime of $R$ such that
        \begin{equation}
\label{eq:tc}
        c\ \ulfrob^m (z) \in \ulfrob^m(\id)\hull R,
        \end{equation}
for all sufficiently big $m$. We  denote the non-standard tight 
closure of an ideal $\id$ by $\tc\id$. A priori, this notion depends on the
choice of $\hull{R}$, that is to say, on the choice of $\Lambda$. If we want to
make this dependence explicit, we  write $\op{cl}_{\Lambda}(\id)$.
It is an interesting (and probably difficult) question to determine whether
different choices of $\Lambda$ give rise to the same closure operation.
Here we take a pragmatic approach:
we are primarily interested in using non-standard tight closure to prove
statements (about $R$) which do not mention it, and for this,
we are free to choose $\Lambda$ to suit our needs.

The next proposition shows that
$\tc{\,\cdot\,}$ shares some basic properties with characteristic $p$ tight
closure. %(See, for example \cite[Proposition 10.1.2]{Bruns-Herzog}.) 
We denote
the set of all elements of a ring $A$ which are not contained in a minimal prime of $A$
by $A^\circ$ (a multiplicatively closed subset of $A$).

\begin{proposition} Let $\id$ and $\mathfrak b$ be ideals of $R$. Then the
following hold:
\begin{enumerate}
\item $\tc{\id}$ is an ideal of $R$ and $\id\subseteq\mathfrak{b}$ 
implies that $\tc{\id}\subseteq\tc{\mathfrak b}$;
\item there exists $c\in R^\circ$ such that $c\ \ulfrob^m 
(\tc{\id})\hull R\subseteq \ulfrob^m(\id)\hull R$ for all sufficiently 
large $m$;
\item $\id\subseteq \tc{\id}=\tc{\tc{\id}}$;
%\item if $\id$ is non-standard tightly closed then so is $(\id:_R\mathfrak b)$.
\item $\tc{\id\cap\mathfrak b} \subseteq \tc{\id}\cap\tc{\mathfrak b}$,
$\tc{\id+\mathfrak b}=\tc{\tc{\id}+\tc{\mathfrak b}}$, and $\tc{\id\mathfrak
b}=\tc{\tc{\id}\tc{\mathfrak b}}$;
\item if $R$ is reduced and  the residue class of
$z\in R$ lies in $\op{cl}_{\Lambda/\pr}\bigl(\id(R/\pr)\bigr)$ 
for each minimal prime $\pr$ of $R$, then
$z\in\op{cl}_\Lambda(\id)$.
\end{enumerate}
\end{proposition}
\begin{proof}
The proofs of the first four properties are as in the case of tight closure
in positive characteristic. Suppose that $R$ is reduced.
Let $\pr_1,\dots,\pr_s$ be all the 
minimal prime ideals of $R$, and for each $j$ choose an element $c_j$ 
inside all minimal primes except $\pr_j$. In particular, 
$c_j\pr_j=0$. By assumption, for each $j$ there exists an element 
$d_j\notin\pr_j$ such that
        \begin{equation*}
        d_j\ulfrob^m(z)\in\ulfrob^m(\id)\hull{ R/\pr_j},
        \end{equation*} for all large $m$. By the discussion in \S\ref{s:quot},
this means that
        \begin{equation}
\label{eq:pri}
        d_j\ulfrob^m(z) \in\ulfrob^m(\id)\hull R+\pr_j\hull R,
        \end{equation}
for all large $m$. Put $c:=c_1d_1+\dots+c_sd_s$; note that $c$ does 
not lie in any minimal prime of $R$. Multiplying \eqref{eq:pri} with 
$c_j$ and taking the sum over all $j$, we get that $c\ulfrob^m(z)$ 
lies in $\ulfrob^m(\id)$, for all large $m$, showing that 
$z\in\tc\id$.
\end{proof}

Next we derive versions of some other well-known results about tight 
closure in prime characteristic.
 We say that an ideal of
$R$ is \emph{non-standard tightly closed}\/ if it is equal to its non-standard
tight closure.

\begin{theorem}\label{T:reg} If $R$ is regular, then every ideal of $R$ is
non-standard tightly closed. 
\end{theorem}
\begin{proof} The image under $\ulfrob^m$ of any regular system of parameters of
$R$ is 
$\hull R$-regular, and by \los\ and \cite[Proposition 1.1.6]{Bruns-Herzog} every
permutation of a $\hull R$-regular sequence in $\hull R$ is $\hull R$-regular.
Hence the $R$-algebra structure on $\hull R$ given by
        \begin{equation}
\label{e:frobm}
         R \to \hull R\colon a\mapsto \ulfrob^m(a)
        \end{equation}
is that of a balanced big \CM\ algebra. Since $R$ is regular, this 
implies that the homomorphism \eqref{e:frobm} 
is flat. (See the remarks preceding the proof of Proposition~\ref{P:ff}.)
Suppose towards a contradiction that $z$ lies in the non-standard 
tight closure of an ideal $\id$ in $R$ but not in $\id$. For some 
non-zero $c\in R$, we have relations \eqref{eq:tc} for $m$ 
sufficiently large. Thus
        \begin{equation*}
        c\in \bigl(\ulfrob^m(\id)\hull R:_{\hull R} \ulfrob^m(z)\bigr) = 
\ulfrob^m (\id:_Rz) \hull R
        \end{equation*} where  we used flatness of \eqref{e:frobm} for 
the last equality. Since $z\notin\id$, the colon ideal $(\id:_Rz)$ is 
contained in $\maxim$. Therefore, $c$  is zero by \eqref{eq:inf}, 
contradiction.
\end{proof}

\begin{remark} For this argument to work, it suffices that \eqref{eq:tc} 
only holds for $m=1$; the ensuing notion is the analogue of what was 
called \emph{non-standard closure} in \cite{SchNSTC}.
\end{remark}

In the next result, we require that $R$ is a homomorphic image of a \CM\ local
ring $S$, say $R= S/I$. In order to get a induced map on the Lefschetz hulls, we
tacitly assume   that $\Lambda$ is   equal to a quotient   $\Gamma/I$ for some
object
$\Gamma$  in $\extcohcat$ whose underlying ring is $S$ (see \S\ref{s:quot}).

\begin{theorem}[Colon Capturing]\label{T:CC} 
Suppose that $R$ is a homomorpic image of a \CM\ local ring and that 
$R$ is equidimensional. If 
$\mathbf{z}=\rij zd$ is a system of parameters of $R$, then for each 
$i=1,\dots,d$, we have an inclusion
        \[
%\label{eq:cc}
        \bigl(\rij z{i-1}R:_R z_i\bigr)\subseteq 
\op{cl} \bigl(\rij z{i-1}R\bigr).
        \]
\end{theorem}
\begin{proof}
Write $R=S/I$ with $S$ a \CM\ local ring and
consider $\mathbf z$ already as
a tuple in $S$. Suppose $I$ has height $e$. By prime avoidance, we can find
$y_1,\dots,y_e\in I$, such that for each
$i$, the ideal $J+\rij ziS$ has height $e+i$, where  $J:=\rij yeS$. In
particular, $(y_1,\dots,y_e,z_1,\dots,z_d)$ is  a system of parameters in $S$,
whence $S$-regular. By
the Unmixedness Theorem (see for instance \cite[Theorem 17.6]{Mats}), the ideal
$J$ has no embedded associated primes. We can now use the same argument as in the
proof of
\cite[Theorem 8.1]{SchNSTC}, to get $c\in S$ not contained in any
minimal prime of $I$ and $N\in\nat$ such that 
\begin{equation}\label{eq:cI}
cI^N\subseteq J.
\end{equation}
Let $a\in S$ be such that its image in $R$ lies in  $\bigl(\rij z{i-1}R: z_i\bigr)$, and
hence
$az_i$ lies in $I+\rij z{i-1}S$. For a fixed $m$, applying $\ulfrob^m$ yields
$$
\ulfrob^m(a)\ulfrob^m(z_i)\in \ulfrob^m(I)\hull S+
\bigl(\ulfrob^m(z_1),\dots,\ulfrob^m(z_{i-1})\bigr)\hull S.
$$
Multiplying this with $c$ and using that $\ulfrob(I)\hull S\subseteq I^N\hull
S$, we get from \eqref{eq:cI} that
$$
c\ulfrob^m(a)\ulfrob^m(z_i)\in J\hull S+
\bigl(\ulfrob^m(z_1),\dots,\ulfrob^m(z_{i-1})\bigr)\hull S.
$$
By the remark before \S\ref{s:nstc}, the sequence
$$\bigl(y_1,\dots,y_e,\ulfrob^m(z_1),\dots,\ulfrob^m(z_d)\bigr)$$ is $\hull
S$-regular, so that  the previous equation can be simplified to
$$
c\ulfrob^m(a)\in J\hull S+
\bigl(\ulfrob^m(z_1),\dots,\ulfrob^m(z_{i-1})\bigr)\hull S.
$$
By our choice of $\Lambda$ 
we have $\hull R = \hull S/I\hull S$.
Taking the reduction modulo $I\hull S$ we get    equations exhibiting   $a$ as an element of the
non-standard tight closure of the ideal $\rij z{i-1}R$. (Note that the image of $c$ lies in $R^\circ$.)
\end{proof}

\begin{remark}\label{R:CC} 
Every complete Noetherian local ring $R$ is a homomorphic image of a \CM\
(in fact, regular) local ring by Cohen's Structure Theorem,
hence Colon Capturing holds for $R$.
If we were able to prove that
        \begin{equation}
\label{eq:complettc}
        \tc\id\overset ?=\tc{\id\complet R}\cap R,
        \end{equation}
for every ideal $\id$ in $R$, then we get Colon Capturing for every 
equidimensional and universally catenary Noetherian local ring $R$.  Note that the
inclusion $\subseteq$ in
\eqref{eq:complettc} is immediate. On the other hand, even for tight 
closure in \ch\ $p$, the other inclusion is   still open. 
Below (see Lemma~\ref{l:cc}), we prove Colon Capturing for
complete reduced $R$ with algebraically closed residue field. 
\end{remark}

Using the previous theorem, we get a direct proof of the celebrated 
Hochster-Roberts Theorem \cite{HR}. %(In the generality stated below, 
%it was first proved in \cite[Theorem 2.3]{HHbigCM2}.) 
A ring \homo\ $A\to B$ is called \emph{cyclically pure} if $\id=\id B\cap 
A$, that is to say, if $A/\id\to B/\id B$ is injective, for every ideal 
$\id$ of $A$. A cyclically pure homomorphism $A\to B$ between local
rings $A$ and $B$ is automatically local. Moreover:

\begin{lemma}\label{l:HR}
Let  $A$ and $B$ be  Noetherian local rings
with respective completions $\complet A$ and $\complet B$.
The completion $\complet A\to \complet B$ of a
cyclically pure homomorphism $A\to B$
is  cyclically pure.
\end{lemma}
\begin{proof}
The homomorphism $B\to\complet B$ is faithfully flat, hence
cyclically pure; thus
the composition $A\to B\to\complet B$ is cyclically pure. 
So from now on we may suppose that $B=\complet B$.
It suffices to show that $\complet A\to B$ is 
injective, since the completion of $A/\id$ is equal to $\complet 
A/\id\complet A$ for any ideal $\id$ in $A$. 
Let $a\in\complet A$ be 
such that $a=0$ in $B$, and for each $i$ choose
$a_i\in A$ such that $a\equiv 
a_i\bmod{\pr}^i\complet A$, where $\pr$ is the maximal ideal of $A$.
Then $a_i$ lies in ${\pr}^iB$, hence
by cyclical purity, in ${\pr}^i$. 
Therefore $a\in{\pr}^i\complet A$ for all $i$, showing that $a=0$ 
in $\complet A$ by Krull's Intersection Theorem.
\end{proof}

\begin{theorem}[Hochster-Roberts]\label{T:HR} 
If there exists a cyclically pure 
\homo\ $R\to S$ into a regular local ring $S$, then $R$ is \CM.
\end{theorem}
\begin{proof} 
By Lemma~\ref{l:HR} we reduce to the case that $R$ and $S$ are complete. 
Let $\rij zd$ be a 
system of parameters in $R$.
We need to show that $\rij zd$ is $R$-regular.
To this end  assume that
        \begin{equation*}
        az_i\in \id:=\rij z{i-1}R,
        \end{equation*} for some $i$ and some $a\in R$.  
Since $R$ is a complete domain, we can apply Theorem~\ref{T:CC},
to get that 
$a\in\op{cl}_{\Lambda}(\id)$, for a suitable choice of $\Lambda$ with
underlying ring $R$.  So for some $c\neq 0$ in $R$ we have
relations \eqref{eq:tc} for all sufficiently large $m$.
Now $R\to S$
induces a \homo\ $\hull R\to 
\hull {S}$. Applying this \homo\ to \eqref{eq:tc} we get that $a$ lies in 
$\tc{\id S}$.
(Note that $c$ remains non-zero in $S$ since $R\to S$ 
is injective.)
Hence $a\in\id S \cap R=\id$ by Theorem~\ref{T:reg} and cyclic purity. 
\end{proof}

\begin{remark}\label{R:Freg} 
We say that  $R$ is \emph{weakly non-standard F-regular} if
every  ideal  is  non-standard tightly closed, for every choice of
$\Lambda$ with underlying ring $R$.
The argument in the proof above actually gives two independent
results. Firstly, if $S$ is weakly non-standard F-regular, then $S$ is 
\CM. Secondly, if $R\to S$  is cyclically pure and  $S$ is weakly 
non-standard F-regular, then so is $R$.
\end{remark}

For some more proofs of this theorem, 
see Remarks~\ref{R:CCgentc} and \ref{R:HRBCM}  below. 
By the argument in the beginning of
the proof of \cite[(2.3)]{HHbigCM2}, the theorem implies the
following global version; for further discussion, see   
Conjecture~\ref{C:bout} in the next section.
%Recall that 
%a ring homomorphism $A\to B$ is \emph{pure} if $M\to M\otimes_A B$ is
%injective for every $A$-module $M$.

\begin{corollary} \label{C:HR}
If $A\to B$ is a pure \homo\ of
Noetherian rings 
containing $\mathbb Q$ and if $B$ is regular, then $A$ is \CM. \qed
\end{corollary}

The \emph{integral closure} of an ideal $J\subseteq S$ of a ring $S$ 
will be denoted by $\overline{J}$. It is the set of all $z\in S$ which are
\emph{integral over $J$,} that is, which satisfy a relation
        \begin{equation}
\label{e:int}
        z^d+a_1z^{d-1}+\dots+a_{d-1}z+a_d=0
        \end{equation}
with $a_i\in J^i$ for each $i$. See \cite[\S10.2]{Bruns-Herzog}  for a proof 
that $\overline{J}$ is an ideal of $S$, and other basic properties of 
$\overline{J}$. The following is a useful characterization of 
integral closure:

\begin{lemma}\label{L:val} Let $S$ be a Noetherian local ring and $J$ an ideal
of $S$. An 
element $z\in S$ is integral over $J$ \iff\ $z\in JV$ for every local 
homomorphism $S\to V$ to a discrete valuation ring $V$ whose kernel 
is a minimal prime of $S$.
\end{lemma}

See \cite[Lemma 3.4]{Huneke-3} for the proof in the case where $S$ is 
a domain; the general case easily reduces to this one; see for 
instance \cite[Lemma 3.2]{Hochster-TIC}.

Before we state the next property of tight closure, we make a general remark:

\begin{lemma}\label{l:int} Let $J$ be an ideal of a ring $S$ and suppose that
$z\in S$ satisfies an integral relation \eqref{e:int}. Then $J^{d-1}z^N\in J^N$
for all $N\in\nat$.
\end{lemma}
\begin{proof} 
We claim that $z^{d+k}\in J^{k+1}$ for all $k\in\nat$. 
We show this by induction on $k$, the case $k=0$ being trivial. For the
inductive step note that 
by \eqref{e:int} we have $$
z^{d+k+1} = -(a_1z^{d+k}+\cdots+a_{k+1}z^d + a_{k+2}z^{d-1} + \cdots +
a_dz^{k+1}).$$
Since $a_{i}z^{d+k+1-i}\in J^{k+2}$ for $i=1,\dots,k+1$ 
(by the inductive hypothesis) and $a_{i}\in J^{i}\subseteq J^{k+2}$
for $i=k+2,\dots,d$, we get that $z^{d+k+1}\in J^{k+2}$ as required. Now clearly
$J^{d-1}z^N\subseteq J^N$ if $N<d$, and by the claim we get
$$J^{d-1}z^N=J^{d-1}z^{d+k}\subseteq  J^{d-1}J^{k+1}=J^N$$ for all $N\geq d$, where
$k:=N-d$. 
\end{proof}

\begin{theorem}[Brian\c{c}on-Skoda]\label{T:BS} For every ideal $\id$ of $R$ we
have $\tc \id\subseteq \overline{\id}$. 
Moreover, if $\id$ has positive height and is generated by at most 
$m$ elements, then the integral closure of $\id^m$ is contained in 
$\tc \id$.
\end{theorem}
\begin{proof}
Let $z\in\tc \id$; so we have a  relation~\eqref{eq:tc} for some 
$c\in R^\circ$  and all sufficiently large $m$. In order to prove 
that $z\in\overline{\id}$, we apply  Lemma~\ref{L:val}. Let $V$ be a \DVR\ and
let $R\to
V$ be a local homomorphism with kernel a minimal prime of $R$. This induces a
\homo\ $\hull R\to\hull V$, and applying this \homo\  to the 
relations~\eqref{eq:tc} shows that $z\in\tc{ \id V}$. (Note that by 
assumption $c\neq 0$ in $V$.) By Theorem~\ref{T:reg}, the latter 
ideal is just $\id V$, and we are done.

Suppose now that $\id$ has positive height and is generated by $\leq m$
elements, and let $z$ lie in the integral closure of $\id^m$. Then 
$z$ satisfies  a relation
        \begin{equation*}
        z^d+a_1z^{d-1}+\dots+a_d=0
        \end{equation*} with $a_i\in\id^{im}$. %Take \sr{s} $\seq Rw$, $\seq
%zw$, $\seq\id w$ 
%and $\seq{a_i}w$ of $R$, $z$, $\id$ and $a_i$ respectively. 
By \los, 
we have for almost all $w$ an integral relation
        \begin{equation*}
        \seq zw^d+\seq{a_1}w\seq zw^{d-1}+\dots+\seq {a_d}w=0
        \end{equation*} with $\seq{a_i}w\in\seq\id w^{im}$ for all $i$. For
those $w$, we get 
for all $N$ that
        \begin{equation*}
        \seq\id w^{m(d-1)}\seq zw^N\subseteq \seq\id w^{Nm}
        \end{equation*} by Lemma~\ref{l:int}. 
For $N$ equal to the $l$th power of the \ch\ of $R_w$ we get  $\seq\id 
w^{Nm}\subseteq\frob w^l(\seq\id w)\seq Rw$, since $\seq\id w$ is 
generated by at most $m$ elements by \los. Hence taking 
ultraproducts, we get
        \begin{equation*}
        \id^{m(d-1)} \ulfrob^l(z)\subseteq \ulfrob^l(\id)\hull R.
        \end{equation*} Since this holds for all $l$ and since we assumed that
$\id$ has 
positive height (hence $R^\circ\cap \id^{m(d-1)}\neq\emptyset$), we 
get that $z\in\tc\id$.
\end{proof}

\begin{remark} The same argument together with \cite[Remark 5.8.2]{HuTC} proves 
under the hypothesis of the theorem that the integral closure of 
$\id^{m+l}$ lies in the non-standard tight closure of $\id^{l+1}$, 
for all $l$.
\end{remark}

\begin{remark}\label{R:norm} It follows that $\tc{\id}=\overline \id$ for each
principal ideal $\id$ in $R$. 
Hence a domain $R$ is normal \iff\ every principal ideal is equal to 
its non-standard tight closure. In particular, using 
Remark~\ref{R:Freg}, we see that a cyclically pure subring of a 
regular local ring (and more generally, a weakly non-standard 
F-regular local ring) is normal.
\end{remark}

We immediately obtain the following classical version of the \BS\ Theorem from
\cite{LS}. (For the ring of convergent power series over 
$\mathbb C$ this was first proved in \cite{BS}; see \cite[\S5]{HuTC} 
or \cite{SchBS} for some more background.)

\begin{theorem}[Brian\c{c}on-Skoda for regular rings] If $A$ is a regular ring
containing $\mathbb Q$ and $\id$   an 
ideal of $A$  generated by at most $m$ elements, then the 
integral closure of $\id^m$ is contained in $\id$. In particular, if $f$ is a
formal power series in $n$ variables over a field 
of \ch\ zero with $f(0)=0$, then $f^n$ lies in the ideal generated by 
the partial derivatives of $f$.
\end{theorem}
\begin{proof} Since this is a local property, we may assume that $A$ is local.
By 
Theorem~\ref{T:BS}, the integral closure of $\id^m$ is contained in 
$\tc\id$, hence in $\id$,  by Theorem~\ref{T:reg}. It is an exercise on the
chain rule to show, using Lemma~\ref{L:val}, 
that $f$ lies in the integral closure of the ideal $J$ generated by 
the partials of $f$. (See \cite[Exercise 5.1]{HuTC}.) 
%or \cite[Fact 
%5.1]{SchBS}.) 
Hence $f^n$ lies in $\overline {J^n}\subseteq J$ by our 
first assertion.
\end{proof}

\subsection{Tight closure---non-local case}\label{s:tcnl} 
Although of minor use,
one can extend the notion of non-standard 
tight closure to an arbitrary Noetherian $\mathbb Q$-algebra $A$ as 
follows.  For every maximal ideal $\frak n$
of $A$ choose a Lefschetz hull $\hull{A_{\frak n}}$ of the
equicharacteristic zero Noetherian local ring $A_{\frak n}$, and write
$\op{cl}_{\frak n}$ for the ensuing notion of non-standard tight closure for
ideals of $A_{\frak n}$.
We define the \emph{non-standard tight 
closure} of an ideal $\id$ of $A$ as the intersection
        \begin{equation*}
        \tc\id:= \bigcap_{{\frak n}\in\op{Max}A} \op{cl}_{\frak n}(\id A_{\frak n})\cap A.
        \end{equation*} 
We invite the 
reader to check that this is indeed a closure operation, admitting 
similar properties as in the local case: for instance, the analogues 
of Theorems~\ref{T:reg} and \ref{T:BS} hold.
% and with the proper 
%notion of parameters in a non-local setup, so does Theorem~\ref{T:CC}. 

If $A^*$ is the product of all $\hull{A_{\frak n}}$ as in
\eqref{e:ff-emb}, then each of its factors admits the action of a non-standard
Frobenius. Let us denote the product of these Frobenii again by $\ulfrob$. We
can now define directly a tight closure operation on ideals in $A$ by mimicking
the definition in the local case, that is to say: $z\in A$ belongs to the
`global' non-standard tight closure of an ideal $\id$ if there is some $c\in
A^\circ$ such that $c\ulfrob^m(z)\in\ulfrob^m(\id)A^*$, for all sufficiently
large $m$. It is immediate that an element in the `global' non-standard tight
closure of $\id$ belongs to $\tc\id$ as defined above. In case $A$ is semi-local, the
converse also holds, but this is no longer clear for arbitrary $A$, for we do not have
 yet an appropriate notion of uniform test elements
for non-standard tight closure (see also Proposition~\ref{P:test} below). This is
presumably not an easy problem, and we will not further investigate it here.

\subsection{Comparison with affine non-standard tight closure}\label{s:affnstc}
We confine ourselves to the geometric 
case, that is, where $R$ is the local ring at a closed point on a 
scheme of finite type over an algebraically closed Lefschetz field 
$k\subseteq K$ as in \S\ref{s:aff}.
In such a ring,
 non-standard tight closure was defined in \cite{SchNSTC} in a similar 
fashion, using the non-standard hull $\ul R$ instead of $\hull R$. 
More precisely, an element $z\in R$ lies in the (affine) non-standard 
tight closure of an ideal $\id$ of $R$ if there exists $c\in R^\circ$ such that
        \begin{equation}
\label{eq:nstc}
        c\ \ulfrob^m(z)\in\ulfrob^m(\id)\ul R
        \end{equation}
for all sufficiently 
large $m$, where we also write $\ulfrob$ for the non-standard 
Frobenius on the Lefschetz ring $\ul R$. As discussed in 
\S\ref{s:aff}, we have a natural embedding $\ul R\to\hull R$, and 
this is compatible with the 
non-standard Frobenii defined on each ring. In particular, taking 
the image of the relations~\eqref{eq:nstc} via this \homo\ shows that 
$z\in\tc\id$  in the present sense. Conversely, suppose 
there exists $c\in R^\circ$ such that
\eqref{eq:tc} holds in $\hull R$ for all sufficiently large $m$. 
By \los, for those $m$ we have that
        \begin{equation}
\label{eq:tcm}
        \seq cw\ \frob w^m(\seq zw)\in\frob w^m(\seq\id w)\seq Rw
        \end{equation}
where $\seq cw$, $\seq zw$, $\seq\id w$ and $\seq Rw$ are \sr{s} of 
$c$, $z$, $\id$ and $R$ respectively. By our discussion in 
\S\ref{s:aff}, we can realize  these \sr{s} as follows. If  $\seqaff 
cp$, $\seqaff zp$, $\seqaff\id p$ and $\seqaff Rp$ are   \sr{s} of $c$, $z$, 
$\id$ and  $R$ in the sense of \cite{SchNSTC}, then we may take $\seq 
Rw$ to be the completion of $\seqaff R{p(w)}\tensor_{\seq k{p(w)}}\seq Kw$ and
$\seq cw$, $\seq zw$
 and $\seq \id w$ the corresponding image of $\seqaff c{p(w)}$, $\seqaff 
z{p(w)}$ and $\seqaff\id{p(w)}$ in this completion. (Recall that $p(w)=
\op{char}K_w$.)
Therefore, by faithful flatness, relation 
\eqref{eq:tcm} already holds in the subring $\seqaff R{p(w)}$, for 
almost all $w$, hence for almost all \ch{s}. Taking ultraproducts of 
this relation in almost all $\seqaff Rp$ yields \eqref{eq:nstc}, and 
since this is true for any sufficiently large choice of $m$, we 
showed that $z$ lies in the non-standard tight closure of $\id$ in 
the sense of \cite{SchNSTC}. In conclusion, we showed that for 
localizations of finitely generated $k$-algebras at maximal ideals, 
both notions of tight closure coincide.

\subsection{Generic tight closure}
We finish this section with studying   a related  closure operation, 
which also played an important role in the affine case. Let 
$\id$ be an ideal of $R$.
We say that an element $z\in R$ lies in the \emph{generic
tight 
closure} of  $\id$ if $\seq zw$ lies in the (positive characteristic)
tight closure of $\seq 
\id w$ for almost all $w$.
We denote the generic tight closure of $\id$ by $\op{cl}^*(\id)$. 
Again, this depends on the choice of $\Lambda$ with underlying ring $R$;
if we want to stress this dependence, we write $\op{cl}_\Lambda^*(\id)$.
From \cite[Appendix 1]{HuTC}  recall Hochster-Huneke's notion of
tight closure in equicharacteristic $0$.
Here and below, given a ring $S$ and a prime $p$ we let
$S(p):=S\tensor_{\zet} \mathbb{F}_p$, and for an ideal $I$ of $S$ we
let $I(p)$ be the image of $I$ in $S(p)$ under the map $z\mapsto z(p):=
z\tensor 1\colon S\to S(p)$.

\begin{definition}\label{d:tc}
 An element
$z$ of $R$ is in the (\emph{equational}) \emph{tight closure} $\id^*$ 
of $\id$ if there exists a finitely generated
subring $S$ of $R$ with $z\in S$ such that 
$z(p)$
is in the (characteristic $p$) tight closure 
of $(\id\cap S)(p)$ in $S(p)$, 
for all but finitely many primes $p$.
\end{definition}

Let $\mathbf{y}=(y_1,\dots,y_m)\in R^m$, and let $J$ be the kernel of
the ring homomorphism $$\zet[Y]=\zet[Y_1,\dots,Y_m]\to R$$ given by $Y_j\mapsto y_j$ for all $j$.
We get an induced embedding $\zet[Y]/J\to R$, and we
identify $\zet[Y]/J$ with its image,
the subring $S$ of $R$ generated by $\mathbf{y}$. Given a prime $p$
we then have $$S(p) = S\tensor_{\zet} \mathbb{F}_p=\mathbb{F}_p[Y]/J(p)$$ where
$J(p)$ is the image of $J$ under the canonical surjection 
$\zet[Y]\to\mathbb{F}_p[Y]$. 
We let $S_\infty$ denote the ultraproduct of the $S(p)$ with respect to
the same ultraset that builds $K$ (and whose underlying set is the set of prime
numbers). The canonical maps $S\to S(p)$ combine to give a ring homomorphism
$S\to S_\infty$. Composing with the diagonal embedding $S_\infty\to 
S_\infty^{\mathcal U}$, where
$\mathcal U$ is the ultraset constructed in \S\ref{s:LP}, we obtain an $S$-algebra
structure on $S_\infty^{\mathcal U}$. We also get an $S$-algebra structure on
$\hull R$  via the restriction of $\eta_R$ to $S$.

\begin{lemma}
There exists an $S$-algebra homomorphism $\varphi\colon S_\infty^{\mathcal U}\to 
\hull R$.
\end{lemma}
\begin{proof}
For every $w$ let $\tilde S_w$ be the subring of $R_w$ generated by
the approximations
$\mathbf{y}_w=(y_{1w},\dots,y_{mw})$ of 
$\mathbf{y}$, and let $\tilde S_\infty$ be the
ultraproduct of the $\tilde S_w$. If $P(Y)\in J$, so $P(\mathbf{y})=0$, then
$P(\mathbf{y}_w)=0$ for almost
all $w$. Therefore, since $J$ is finitely generated, we have for almost all $w$
a surjection $ S\bigl(p(w)\bigr)\to \tilde S_w$
with $y_{j}(p(w))\mapsto y_{jw}$ for all $j$. Let $\varphi_w\colon
S\bigl(p(w)\bigr)\to\seq Rw$ denote the
composition
of this surjection with the embedding  $\tilde S_w\subseteq \seq Rw$ and let
$\varphi$ be the 
ultraproduct of the $\varphi_w$. One easily checks that $\varphi$ is an  $S$-algebra \homo\ $S_\infty^{\mathcal U}\to 
\hull R$.
\end{proof}

\begin{remark}\label{R:modp}
This means in particular that for every $z\in S$,  the $z_w:=\varphi_w(z(p(w)))$
are an \sr\ of $z$. Indeed,
$z=\varphi(z)$ is by construction the ultraproduct of the $\varphi_w(z(p(w)))$.
\end{remark}

\begin{corollary}
For every ideal $\id$ of $R$, we have $\id^*\subseteq \op{cl}^*(\id)$.
\end{corollary}
\begin{proof}
Let $z\in\id^*$, and choose 
$S=\zet[\mathbf{y}]$, where
$\mathbf{y}=(y_1,\dots,y_m)\in R^m$, which contains $z$ and such that
$z(p)$ is in the tight closure of $(\id\cap S)(p)$ in $S(p)$, 
for all but finitely
many $p$. Then $z\bigl(p(w)\bigr)$ is in the tight closure  of 
$\seq Iw:=(\id\cap S)\bigl(p(w)\bigr)$ in 
$S\bigl(p(w)\bigr)$, for almost all $w$. 
By \cite[Theorem 2.3]{HuTC}, almost each $\seq zw:=\varphi_w\bigl(z\bigl(p(w)\bigr)\bigr)$
 is in the tight closure of $\varphi_w(\seq Iw)\seq Rw$. By
Remark~\ref{R:modp}, the $\seq zw$ and the $\varphi_w(\seq
Iw)\seq Rw$ are \sr{s} of $z$ and $(\id\cap S)R$ respectively. In particular,
if $\id_w$ is an \sr\ of $\id$, then almost each $\seq
zw$ lies in the tight closure of $\id_w$, showing that $z\in\op{cl}^*(\id)$.
\end{proof}

The  relation between generic tight closure and non-standard tight closure
is more subtle. We 
need a result on test elements.
(See \cite[Chapter 2]{HuTC} for the notion of test element.)

\begin{proposition}\label{P:test} 
Suppose that $\hull R$ is   normalizing and $R$ is absolutely
analytically irreducible. 
There exists an element of $\complet R$ almost all of
whose approximations are test elements.
\end{proposition}
\begin{proof} 
The assumption on $\hull R$ implies that the \homo\
$T_0:=\pow {k^*}d\to \complet R$ given by $X_i\mapsto x_i$ is a Noether
Normalization, where $d:=\dim R$ and $k^*$ is the algebraic closure of $k$ in
$\complet R$ (whence a coefficient field of $\complet R$). Moreover, this \homo\
induces (by extension of scalars) the restriction of
 $\complet\theta_\Lambda\colon K[[n]]\to \complet R_{(k,u)}$ to $T:=\pow Kd$
  (see \S\ref{s:NN}). Let $c$ be a non-zero element in the 
\emph{relative Jacobian} $J_{\complet R/T_0}$. (Recall that $J_{\complet R/T_0}$
is the  $0$-th Fitting ideal of the relative module of K\"ahler differentials 
$\Omega_{\complet R/T_0}$.)   By Remark~\ref{R:unr2}, 
almost each $\seq Rw$ is a
domain, and by Proposition~\ref{P:fin}, almost each $\seq Rw$ is a finite extension  of $\seq
Tw$, of degree $e$.  
In particular, for almost all $w$, the field 
of fractions of $\seq Rw$ is separably algebraic over the field of fractions of 
$\seq Tw$. Since $J_{\complet R/T_0}\subseteq J_{\complet R_{(k,u)}/T}$, almost
each $\seq cw$ is  a non-zero element of $J_{\seq Rw/\seq Tw}$, hence a test 
element for $\seq Rw$  by \cite[Exercise 2.9]{HuTC}.
\end{proof}

\begin{remark}
From this we can also derive the same result for 
$R$ analytically unramified with $k$
algebraically closed and $\hull R$ normalizing.
Namely,
 let $\pr_1,\dots,\pr_s$ be the minimal prime ideals of 
$\complet R$.
By Remark~\ref{R:unr2}
the approximations $\pr_{1w},\dots,\pr_{sw}$ are the minimal prime
ideals of $R_w$, for almost all $w$, and  almost all
$\seq Rw$ are reduced.
For each 
$j$, choose $t_j\in\complet R$ 
inside all minimal prime ideals except $\pr_j$. By Proposition~\ref{P:test}, there
exists $c_j\in\complet R\setminus \pr_j$ whose \sr\  is a test 
element for almost all $\seq Rw/\pr_{jw}$. Using \los, one
shows that $c:=c_1t_1+\dots+c_st_s$ has the desired properties.
(See for instance \cite[Exercise 2.10]{HuTC} for more details.)
\end{remark}

\begin{theorem}\label{T:gen} 
Suppose that $R$ is complete  and $\hull R$ is  normalizing. If $R$
is either absolutely analytically irreducible or otherwise reduced with  
$k$ algebraically closed, then
$\op{cl}^*(\id)\subseteq\tc{\id}$ for every ideal $\id$ of $R$.
\end{theorem}
\begin{proof} 
Let $z\in\op{cl}^*(\id)$, that is, $z_w$ is in the tight closure of 
$\id_w^*$ for almost all $w$.
By either Proposition~\ref{P:test} or the remark following it, there exists an
element $c$ of $R$ whose
\sr\ $\seq 
cw$ is a test element in $\seq Rw$, for almost all $w$. Hence   for 
almost $w$ and for all $m$:
        \begin{equation*}
        \seq cw\frob w^m(\seq zw)\in \frob w^m(\seq \id w)\seq Rw.
        \end{equation*}
 Taking  ultraproducts, we get for all $m$ that
        \begin{equation*}
        c\ \ulfrob^m(z)\in\ulfrob^m(\id)\hull R
        \end{equation*} showing that $z\in\tc\id$.
\end{proof}

For the Hochster-Huneke notion of tight closure in equicharacteristic zero, 
Colon Capturing is only known to be true 
in locally excellent rings.
Since Colon Capturing holds for every complete Noetherian local ring 
of positive \ch, hence for every \sr\ of $R$, \los\ in conjunction with
Lemma~\ref{L:sop} immediately yields:

\begin{lemma}[Colon Capturing for generic tight closure]\label{l:cc}
If $(z_1,\dots,z_d)$ is a system of parameters of $R$, then 
        \[
%\label{eq:cc}
        \bigl(\rij z{i-1}R:_R z_i\bigr)\subseteq 
\op{cl}^*\bigl(\rij z{i-1}R\bigr)
        \]
for each $i=1,\dots,d$. \qed
\end{lemma}

In particular, combining this lemma with Theorem~\ref{T:gen} yields Colon
Capturing  for non-standard tight closure in case $R$ 
is reduced and complete, with  
algebraically closed $k$ and $\hull R$ normalizing.

\begin{remark}\label{R:CCgentc}
It follows from Theorem~\ref{T:trans} 
that every ideal in an equi\ch\ zero regular local ring is equal to 
its generic tight closure. Together with Lemma~\ref{l:cc}, we 
get an even easier proof of the Hochster-Roberts Theorem (including 
the global version of Corollary~\ref{C:HR}), using $\op{cl}^*$ in place
of $\op{cl}$.
\end{remark}

\section{Balanced Big \CM\ Algebras}\label{s:BCM}

 Recall that an 
$R$-algebra $B$ is called a \emph{balanced big \CM\ $R$-algebra}, if 
any system of parameters of $R$ is a $B$-regular sequence. (If we only 
know this for a single system of parameters, we call $B$ a \emph{big 
\CM\ $R$-algebra}.)  The key result on big \CM\ algebras was proved 
by Hochster-Huneke in \cite{HHbigCM}: if $S$ is an excellent local 
domain  of prime \ch\ $p$, then its absolute integral closure $S^+$ 
is a balanced big \CM\ algebra. (Incidentally, this is false in equi\ch\ 
zero if $\dim S\geq 3$, see \cite{HHbigCM}.) The \emph{absolute integral closure} $A^+$ of a
domain $A$ is defined to be the integral closure of $A$ in 
an algebraic closure of its field of fractions. (We put $A^+:=0$ 
if $A$ is not a domain.) In \cite{SchBCM}, this is used  to give a 
canonical construction of a balanced big \CM\ algebra  for a local 
domain $S$ essentially of finite type over $\mathbb C$, by taking the 
ultraproduct of the $\seq Sp^+$, where $\seq Sp$ is an \sr\ of $S$ in 
the sense of \cite{SchNSTC}. The $\seq Sp$ are local 
domains, for almost all $p$, by \cite[Corollary 4.2]{SchNSTC}, so 
that the construction makes sense. In view of the restrictions 
imposed by Theorem~\ref{T:dom}, we cannot directly generalize this to 
arbitrary domains. 
We first consider the case that
$\hull\Lambda$ is a domain, or equivalently, that
almost all \sr{s}  $\seq Rw$ of $R$ are domains.
This is the case if $R$ is absolutely irreducible
and $\Lambda$ is absolutely
normalizing or normalizing
(by Theorem~\ref{T:dom} and Remark~\ref{R:unr2}, respectively), but also
if $R$ is a DVR (by Corollary~\ref{C:DVR}).

\begin{definition}\label{D:bigCM}
        \begin{equation*}
        \BCM\Lambda :=\up w \seq Rw^+.
        \end{equation*}
We often write $\BCM R$ for $\BCM\Lambda$, keeping in mind that
$\BCM R$  depends on the choice of $\Lambda$. 
%Where
%necessary, we will assume without further warning that 
%$\Lambda$ has the required
%additional properties for invoking results from \S\S\ref{s:irr}--\ref{T:Ri}.
\end{definition}

The canonical \homo\ $\eta_R\colon R\to \hull R$ induces a \homo\ 
$R\to \BCM  R$, turning $\BCM  R$ into an $R$-algebra. (Note that this 
is no longer an integral extension.) Since the $\seq Rw$ are complete 
(hence Henselian), the $\seq Rw^+$ are local, whence so is 
$\BCM R$. Moreover,  the canonical \homo\ $R\to \BCM R$ is local.

\begin{theorem}\label{T:BCMcomp}
The $R$-algebra
$\BCM  R$ is a balanced big \CM\ algebra.
If $\alpha\colon\Lambda\to\Gamma$ 
is a morphism in $\extcohcat$
with underlying ring homomorphism $R\to S$, where 
$\hull\Gamma$ is a domain,
then there exists a \textup{(}non-unique\textup{)} \homo\ 
$\widetilde\alpha\colon \BCM  R\to \BCM  S$ giving rise to a commutative 
diagram
\begin{equation}\label{BRS}
        \begin{aligned}
        \mbox{
        \xymatrix@R+2em@C+3em{R \ar[r]^\alpha \ar[d] & 
S       \ar[d] \\
        \BCM R \ar[r]^{\widetilde \alpha} &{\BCM S}.}
        }
        \end{aligned}
\end{equation}
Moreover, if $\alpha$ is finite, injective, and induces an
isomorphism on the residue fields, then $\BCM  R=\BCM  S$.
\end{theorem}
\begin{proof}
Let $\mathbf z$ be a system of 
parameters in $R$ with \sr{s}  $\seq{\mathbf z}w$. By 
Lemma~\ref{L:sop} almost each $\seq{\mathbf z}w$ is a  system of 
parameters in $\seq Rw$, hence is  $\seq Rw^+$-regular by \cite{HHbigCM}. By 
\los,  $\mathbf z$ is $\BCM  R$-regular.
 From the \homo\ $\hull\alpha\colon \hull R\to\hull S$ we get \homo{s} 
$\seq Rw\to\seq Sw$ for almost all $w$, where $\seq Sw$ is an \sr\ of 
$S$. These extend (non-uniquely) to \homo{s} $\seq Rw^+\to\seq Sw^+$ 
whose ultraproduct is the required $\widetilde  \alpha$. 
If $\alpha$ is finite, 
injective, and induces an isomorphism on the residue fields,
then almost all $\seq Rw\to\seq Sw$ are finite and injective by 
Proposition~\ref{P:fin}, and hence $\seq 
Rw^+=\seq Sw^+$. The last assertion is now clear.
\end{proof}

\begin{remark}
Incidentally, the argument at the end of the proof shows that 
there is essentially only one ring in each 
dimension $d$ playing the role of a big \CM\ algebra: 
Indeed, suppose that the restriction of ${\complet\theta}_\Lambda$ to $K[[d]]$ 
is a Noether normalization of $\complet R_{(k,u)}$,
where $d=\dim R$.
(This is satisfied, 
for example, if $\Lambda$ is absolutely normalizing.) 
Then
$\BCM R$ is isomorphic (non-canonically) to $\BCM{K[[d]]}$.  
\end{remark}

\begin{corollary}\label{C:ff}
If $R$ is regular, then the $R$-algebra $\BCM R$ is faithfully flat.
\end{corollary}
\begin{proof}
We already  mentioned that a balanced big \CM\ algebra over a regular 
local ring is automatically flat; see the remarks before 
Lemma~\ref{l:balanced}. 
Since  $R\to \BCM R$ 
is  local, it is therefore faithfully flat.
\end{proof}

\begin{remark}\label{R:HRBCM}
This gives us a second direct proof of the Hochster-Roberts Theorem 
(Theorem~\ref{T:HR}): with notation from the theorem, we may reduce 
again to the case that $R$ and $S$ are complete and that $R$ has 
algebraically closed residue field. Suppose $\rij zd$ is a system of 
parameters in $R$ and let $az_i\in \id:=\rij z{i-1}R$. Since $\rij 
zd$ is $\BCM R$-regular by Theorem~\ref{T:BCMcomp}, we get $a\in 
\id\BCM R$. Choose absolutely normalizing objects $\Lambda$ and $\Gamma$
of $\extcohcat$ with underlying rings $R$ and $S$, respectively, such that
$R\to S$ becomes an $\extcohcat$-morphism, hence
induces a
\homo\ $\BCM R\to \BCM S$ which makes 
diagram~\eqref{BRS} commutative. Then $a\in \id\BCM S$. Since 
$S\to\BCM S$ is faithfully flat by Corollary~\ref{C:ff}, we get $a\in 
\id S$ and hence, by cyclical purity, $a\in \id$.
\end{remark}

As in positive \ch, the ring $\BCM R$ has many additional properties 
(which fail to hold for the big \CM\ algebras in equi\ch\ zero 
constructed by Hochster-Huneke in \cite{HHbigCM2}). For instance, 
$\BCM  R$ is absolutely integrally closed, hence in particular 
quadratically closed, and therefore, the sum of any number of prime 
ideals is either the unit ideal or again a prime ideal (same argument 
as in \cite[\S3]{SchBCM}).
Moreover:

\begin{proposition}\label{P:sur}
The canonical map $\op{Spec}\BCM R\to \op {Spec} R$ is surjective.
\end{proposition}
\begin{proof}
Let $\pr$ be a prime ideal in $R$ and let $\mathfrak q$ be a prime 
ideal in $\complet R_{(k,i)}$ lying over $\pr$. 
%Let $\seq Rw$ and 
%$\seq{\mathfrak q} w$ be \sr{s} of $R$ and $\mathfrak q$ 
%respectively. 
By Theorem~\ref{T:dom},  almost all \sr{s} $\seq {\mathfrak q} 
w$ of $\frak q$ 
are prime ideals. Since $\seq Rw\subseteq\seq Rw^+$ is integral, 
there exists a prime ideal $\seq {\mathfrak Q}w$ in $\seq Rw^+$ whose 
contraction to $\seq Rw$ is $\seq{\mathfrak q} w$. The ultraproduct 
of the $\seq{\mathfrak Q}w$ is then a prime ideal in $\BCM R$ whose 
contraction to $R$ is $\pr$.
\end{proof}

\subsection{Big \CM\ algebras--- general case}
We now define $\BCM R=\BCM\Lambda$ for
an arbitrary equi\ch\ zero Noetherian local ring  $(R,\maxim)$, under the
assumption that $\Lambda$ is absolutely normalizing:
        \begin{equation*}
        \BCM R:= \bigoplus_{\mathfrak P} \BCM{\complet R_{(k,i)}/\mathfrak P}
        \end{equation*}
where $\mathfrak P$ runs over all prime ideals of $\complet R_{(k,i)}$ of 
maximal dimension (that is to say, such that $\op{dim}(\complet 
R_{(k,i)}/\mathfrak P)=\op{dim}R$). Note that this agrees with 
our former definition in case $\hull\Lambda$ (and hence $\complet R_{(k,i)}$)
is a domain.
Clearly, $\BCM R$ inherits an 
$R$-algebra structure via the $\complet R_{(k,i)}/\mathfrak P$-algebra 
structure on each summand. We claim that $\BCM R$ is a balanced big 
\CM\ algebra. Indeed, if $\mathbf z$ is a system of parameters in 
$R$, then it remains so in $\complet R_{(k,i)}$ and hence in each $\complet 
R_{(k,i)}/\mathfrak P$ since the $\mathfrak P$ have maximal dimension. 
Therefore, by Theorem~\ref{T:BCMcomp},  for each $\mathfrak P$, the 
sequence $\mathbf z$ is $\BCM{\complet R_{(k,i)}/\mathfrak P}$-regular, 
hence $\BCM R$-regular. All the properties previously stated in the case
that $\hull\Lambda$ is a domain remain true 
in this more general setup.

As in the Hochster-Huneke construction, there is a weak form of 
functoriality. We need a definition taken from \cite{HHbigCM2} (see 
also \cite[\S9]{HuTC}).

\begin{definition}
We say that a local \homo\ $R\to S$ of Noetherian local rings is 
\emph{permissible} %(or, that $S$ is a \emph{permissible 
%local $R$-algebra}), 
if for each   prime ideal $\mathfrak q$  
in $\complet S$  of maximal 
dimension, we can find a prime ideal $\pr$   
in $\complet R$  of 
maximal dimension such that $\pr \subseteq\mathfrak 
q\cap\complet R$. A $\extcohcat$-morphism is called
{\em permissible} if its underlying ring homomorphism is permissible.
\end{definition}

As remarked in \cite[\S9]{HuTC}, any local \homo\ with source an 
equidimensional and universally catenary local ring  is permissible.
Moreover:

\begin{lemma}
If $\Lambda\to\Gamma=(S,\mathbf{y},l,v)$ is a permissible $\extcohcat$-morphism
then 
the homomorphism $\complet R_{(k,i)}\to \complet S_{(l,j)}$ is permissible.
\end{lemma}
\begin{proof}
Recall that $i$ and $j$ denote the respective embedding of $k$ and $l$ into the
algebraic closures $\tilde{k}$ and $\tilde{l}$  of $u(k)$ and $v(l)$ inside
$K$. Let 
$\mathfrak Q$ a  prime ideal  of maximal dimension in $\complet S_{(l,j)}$ 
and let $\mathfrak q$ be its contraction to $\complet S$. We have 
inequalities
        \begin{equation}\label{eq:dim}
        \op{dim}(\complet S_{(l,j)})= \op{dim}\bigl(\complet S_{(l,j)}/\mathfrak 
Q\bigr)\leq \op{dim} (\complet S/\mathfrak q)\leq \op{dim}(\complet S)
        \end{equation}
where the middle inequality follows from \cite[Theorem 15.1]{Mats}, 
since the closed fiber is trivial. As $\complet S_{(l,j)}$ has the same 
dimension as $\complet S\tensor_l \tilde{l}$ 
and therefore as $\complet S$, all inequalities in 
\eqref{eq:dim} are equalities, so that $\mathfrak q$ is a  prime 
ideal of maximal dimension. By assumption, there is a   prime ideal 
$\pr$  in $\complet R$  of maximal dimension contained in $\mathfrak 
q$. By faithful flatness, 
$\complet R_{(k,i)}/\pr\complet R_{(k,i)}$ has dimension 
$\op{dim}(\complet R)=\op{dim}(\complet R_{(k,i)})$. 
Since $\complet R/\pr$ is universally catenary and equidimensional, 
so is $\complet R_{(k,i)}/\pr\complet R_{(k,i)}$. Therefore, if  $\mathfrak P$ is 
a minimal prime of $\pr\complet R_{(k,i)}$ contained in $\mathfrak Q$, then 
it has maximal dimension, as required.
\end{proof}

We turn to the definition of $\BCM R\to \BCM S$ for  a 
permissible \homo\ $R\to S$:

\begin{corollary}
Given a permissible $\extcohcat$-morphism $\alpha\colon\Lambda\to\Gamma$  
with $\Gamma$ absolutely normalizing, there exists a 
homomorphism $\widetilde\alpha\colon \BCM{\Lambda}=\BCM{R}\to \BCM{\Gamma}=
\BCM{S}$ making \eqref{BRS} commutative.
\end{corollary}
\begin{proof}
By the lemma, for each 
prime ideal $\mathfrak Q$ in $\complet S_{(l,j)}$ of maximal dimension 
 we can choose  a 
prime ideal $\mathfrak Q'$ of maximal dimension in $\complet R_{(k,i)}$ 
such that $\mathfrak Q'\subseteq \mathfrak Q$. Fix one such prime ideal 
$\mathfrak Q'$ for each $\mathfrak Q$. The \homo\ $$\complet 
R_{(k,i)}/\mathfrak Q'\to\complet S_{(l,j)}/\mathfrak Q$$ induces by 
Theorem~\ref{T:BCMcomp} a \homo\
        \begin{equation*}
        j_{\mathfrak Q}\colon \BCM{\complet R_{(k,i)}/\mathfrak Q'} \to 
\BCM{\complet S_{(l,j)}/\mathfrak Q}.
        \end{equation*}
Define $\BCM R\to \BCM S$ now by sending a tuple $(a_{\mathfrak P})$ 
with $a_{\mathfrak P}\in\BCM{\complet R_{(k,i)}/\mathfrak P}$ and 
$\mathfrak P$ a prime ideal in $\complet R_{(k,i)}$ of maximal dimension, 
to the tuple $(j_{\mathfrak Q}(a_{\mathfrak Q'}))$, where $\mathfrak 
Q$ runs over all prime ideals in $\complet S_{(l,j)}$ of maximal dimension. 
It is easy to see that this gives rise to a commutative 
diagram~\eqref{BRS}.
\end{proof}

It is also easy to see that if  $\alpha\colon\Lambda\to\Gamma$  
is a permissible
$\extcohcat$-morphism where $\hull\Gamma$ is a domain, then there exists a 
homomorphism $\widetilde\alpha$ 
making \eqref{BRS} commutative. Calling $\Gamma$ {\em permissible}
if $\hull\Gamma$ is a domain or absolutely normalizing
(so $\BCM\Gamma$ is defined), we therefore have:

\begin{corollary}\label{C:perm}
Given a permissible $\extcohcat$-morphism $\alpha\colon\Lambda\to\Gamma$  
between permissible $\extcohcat$-objects, there exists a 
homomorphism $\widetilde\alpha\colon  \BCM{\Lambda}\to \BCM{\Gamma}$
making \eqref{BRS} commutative. \qed
\end{corollary}

To show the strength of the existence of big \CM\ algebras, let us 
give a quick proof of the Monomial Conjecture.

\begin{corollary}[Monomial Conjecture]\label{C:MC}
Given a system of parameters  $\rij zd$ in the equi\ch\ zero 
Noetherian local ring $R$, we have for all $t\in\nat$ that
        \begin{equation}\label{eq:MC}
        (z_1z_2\cdots z_d)^t\notin (z_1^{t+1},\dots,z_d^{t+1})R.
        \end{equation}
\end{corollary}
\begin{proof}
The sequence $\rij zd$ is $\BCM R$-regular and so \eqref{eq:MC} holds 
in $\BCM R$, hence {\em a fortiori}\/ in $R$.
\end{proof}

The above proof   does rely on the result of Hochster and Huneke that 
absolute integral closure in positive \ch\ yields big \CM\ algebras. 
A more elementary argument is obtained by using Lemma~\ref{L:sop}
together with the observation that the Monomial Conjecture admits an 
elementary proof in positive \ch\ \cite[Remark 
9.2.4(b)]{Bruns-Herzog}. An equally quick proof, which we will not 
produce here, relying also on the weak functoriality property of 
$\mathfrak B$, can be given for the Vanishing Theorem of maps for Tor 
\cite[Theorem 4.1]{HHbigCM2}.

\subsection{\pc}
As in \cite{SchBCM}, we can use our construction of a big 
\CM\ algebra to define yet another closure operation on ideals of $R$
as follows. Suppose that $\Lambda$ is permissible, and
let $\id$ be an ideal of $R$. The 
\emph{\pc}    of  $\id$ in $R$ is by definition
        \begin{equation*}
        \id^+:= \id\BCM R\cap R.
        \end{equation*}
We next show that the  analogues of Theorems~\ref{T:reg}, 
\ref{T:CC} and \ref{T:BS} hold for $\id^+$ in place of $\tc{\id}$. As 
for the last property in the next theorem, persistence,  it is not 
immediately clear that it also holds for non-standard tight closure. 
We also remind the reader that if $R$ is equidimensional and 
universally catenary (for instance,  an excellent domain), then  every
local $R$-algebra is permissible.

\begin{theorem}\label{T:per}
Let $\id$ be an ideal of $R$.
\begin{enumerate}
\item\label{i:reg}  If $R$ is regular, then $\id=\id^+$.
\item\label{i:cc}  If $\rij zd$ is a system of parameters in $R$, 
then $$\bigl((z_1,\dots,z_{i-1}):_R z_i\bigr)\subseteq (\rij z{i-1}R)^+$$ 
for all 
$i$ \emph{(Colon Capturing)}.
\item\label{i:BS} We have $\id^+\subseteq 
\overline\id$, and if $\id$ is generated by $m$ elements, then
        \begin{equation*}
        \overline{\id^{l+m}}\subseteq (\id^{l+1})^+
        \end{equation*}
  for all $l$ \emph{(\BS)}.
\item\label{i:per} If  $\Lambda\to\Gamma=(S,\dots)$ is
a permissible morphism between
permissible objects in $\extcohcat$, then $\id^+S\subseteq 
(\id S)^+$ \emph{(Persistence)}.
\end{enumerate}
\end{theorem}
\begin{proof}
For \eqref{i:reg}, observe that $\complet R_{(k,i)}$ is again regular
(Lemma~\ref{l:scalarext}), so 
that the composition $R\to \complet R_{(k,i)}\to \BCM R$ is faithfully 
flat, by   Corollary~\ref{C:ff}, hence cyclically pure. For 
\eqref{i:cc}, let $I:=\rij z{i-1}R$ and suppose $az_{i}\in I$. Since 
$\rij zd$ is $\BCM R$-regular, we get $a\in I\BCM R$, and hence $a\in 
I^+$. 
%In fact, the argument yields that $x_{i+1}$ is a non-zero 
%divisor modulo $I^+$. 
The argument in \cite[\S6.1]{SchBCM} (the 
affine case) can be copied almost verbatim to prove the second 
assertion in \eqref{i:BS}; for the first assertion, we use 
Lemma~\ref{L:val} together with \eqref{i:reg} in the same way as in 
the proof of Theorem~\ref{T:BS}. (Note that $R\to V$ is automatically 
permissible, where $V$ is as in Lemma~\ref{L:val}, and every 
$\extcohcat$-object with
underlying ring $V$ is permissible,
so that we get a 
\homo\ $\BCM R\to \BCM V$, by Corollary~\ref{C:perm}.)  
Persistence is immediate from  weak 
functoriality of $\mathfrak B$.
\end{proof}

Conjecturally, in \ch\ $p$, plus closure and tight closure coincide. 
A \ch\ zero analogue of this is that \pc\ and generic tight closure 
should be the same. We have at least the following analogue 
of \cite[Theorem 5.12]{HHbigCM2}.  (The second statement relies on 
Smith's work \cite{SmParId}).

\begin{proposition}\label{P:plus}
Suppose  $R$ is formally equidimensional.
For each ideal $\id$ of $R$, we have $\id^+\subseteq \op{cl}^*(\id)$.  
If $\id$ is generated by a system of parameters, then $\id^+=\op{cl}^*(\id)$.
\end{proposition}
\begin{proof}
We give the proof in the case that $\hull R$ is absolutely normalizing, the
case that $\hull R$ is a domain being similar (and simpler).
In view of Lemma~\ref{l:scalarext}, passing from $R$ to $\complet R_{(k,i)}$ 
reduces the problem to the case that $R$ is complete and equidimensional, 
with $k$
algebraically closed. (Note that
both \pc\ and generic tight closure commute with such an extension of scalars).
By Corollary~\ref{c:equidim}, 
almost all $\seq Rw$ are  equidimensional, and
their minimal primes $\seq{\pr_j}w$ are \sr{s} of the 
minimal primes $\pr_j$ of $R$.  By definition, $\BCM 
R$ is the direct sum of the $\BCM{R/\pr_j}$. Suppose $z\in\id^+$, so 
that $z\in \id\BCM{R/\pr_j}$ for each $j$. Hence $\seq zw\in
\seq\id w(\seq Rw/\seq{\pr_j}w)^+$ for all $j$ and almost all $w$. 
If $B$ is an integral extension of a Noetherian domain $A$ of 
positive characteristic and $I$ is an ideal of $A$, then
$IB\cap A$ is contained in the tight closure of $I$
\cite[Theorem 1.7]{HuTC}. Thus almost 
all $\seq zw$ lie  in the tight closure of $\seq\id w(\seq 
Rw/\seq{\pr_j}w)$, hence in the tight closure of $\seq\id w$ (since 
this holds for all minimal primes). This means that $z\in
\op{cl}^*(\id)$.  

Suppose that $\id$ is generated by a system of parameters.
By 
Lemma~\ref{L:sop}, almost all $\seq\id w$ are generated by a 
system of parameters, and this remains true in the homomorphic images 
$\seq Rw/\seq{\pr_j}w$. By \cite{SmParId}, the tight closure of 
$\seq\id w(\seq Rw/\seq{\pr_j}w)$ is contained in
$\seq\id w(\seq Rw/\seq{\pr_j}w)^+$. Taking ultraproducts yields
$\op{cl}^*(\id)\subseteq
\id\BCM R$.
\end{proof}

\begin{remark}
Suppose that $R$ is complete and $\hull R$ is normalizing.
If $R$ is either absolutely analytically irreducible or reduced,
equidimensional with $k$ algebraically closed, then the previous result
in combination with
Theorem~\ref{T:gen}  yields an inclusion
$\id^+\subseteq \tc\id$.
\end{remark}

\subsection{Comparison with big \CM\ algebras for affine local domains}
We want to compare the present construction with the one from 
\cite{SchBCM} discussed in the introduction of this section.
We restrict ourselves 
once more to the case that $R$ is the localization of a finitely 
generated $k$-algebra at a maximal ideal, with $k$ an algebraically closed
Lefschetz field contained in $K$ as in \S\ref{s:aff}; we continue to
use the notations introduced there.
Let $\seqaff Rp$ 
 denote an \sr\ of $R$ in the sense of \cite{SchNSTC}.  Recall that 
the approximations $\seq Rw$ in the sense of the present paper
are defined as 
$$R_w := \text{completion of } \seqaff R{p(w)}\tensor_{\seq k{p(w)}}\seq Kw, \qquad\text{where 
$p(w)=\op{char} \seq Kw$.}$$
Suppose that $R$ is an integral domain; then almost every $\seqaff Rp$ is
a domain.
In general, $\complet R$ and the $\seq Rw$,  though reduced 
and equidimensional,  will no longer be domains. So let 
$\pr_1,\dots,\pr_s$ be the minimal primes of $\complet R$.
Suppose that $\Lambda$ is absolutely normalizing. It follows from
Corollary~\ref{c:equidim} that for
almost all $w$, the $\seq{\pr_j}w$ are the minimal 
prime ideals of $\seq Rw$, and from Theorem~\ref{T:trans}, that they 
have maximal dimension. By definition, $\BCM R$ is the direct sum of 
all $\BCM{\complet R/\pr_j}$.  The ultraproduct $B(R)$ of the $(\seqaff Rp)^+$
is a big \CM\ $R$-algebra; see \cite{SchBCM}. For each $w$ and each 
$j$, the composition $$\seqaff R{p(w)}\to \seq Rw\to \seq Rw/\seq{\pr_j}w$$ is
injective and can
be extended (non-uniquely) to a \homo\
$$\bigl(\seqaff R{p(w)}\bigr)^+\to (\seq Rw/\seq{\pr_j}w)^+.$$
 By construction
        \begin{equation*}
        \up w \bigl(\seqaff R{p(w)}\bigr)^+ \iso  B(R)^{\mathcal U}
        \end{equation*}
where $\mathcal U$ is the ultraset from \S\ref{s:LP}. Therefore, the 
composition of the diagonal embedding with the sum of the 
ultraproducts of the \homo{s} $$\bigl(\seqaff R{p(w)}\bigr)^+\to (\seq 
Rw/\seq{\pr_i}w)^+$$ yields a \homo\ $B(R)\to \BCM R$. The reader can 
verify that this fits in a commutative diagram
        \begin{equation*}
        \xymatrix@R+2em@C+3em{ {R} \ar[d]^{} \ar[r] &
{B(R)} \ar@{->}[dl]^{}\\
{\BCM R} }
        \end{equation*}
%We mention without proof that one can also ensure a form of weak 
%functoriality. 
In \cite{SchBCM}, the \emph{$B$-closure} of an ideal 
$\id$ of $R$ is defined as the ideal $\id B(R)\cap R$. Clearly, we have 
$\id B(R)\cap R\subseteq \id^+$ and we suspect that both are equal. For this 
to be true, it would suffice to show that the \homo\ $B(R)\to \BCM R$ is 
cyclically pure. (Note, however, that $B(R)$  is not local.) We leave it  to
the reader to verify that the discussion in \S\ref{s:affnstc} also applies to generic tight
closure, that is to say,   the two notions, the present one and the
`affine' one from \cite{SchNSTC}, coincide for 
localizations of finitely generated $k$-algebras at maximal ideals. Using this
together with  Proposition~\ref{P:plus} and \cite[Corollary 4.5]{SchBCM}, we
get an equality $\id B(R)\cap R= \id^+$ for $\id$ an ideal generated by
a system of parameters.

%\begin{proposition}
%Suppose that
%$R$ is normal and equidimensional. For suitable choice of $\Lambda$ we 
%then have
%$\id B(R)\cap R = \id^+$ for every ideal $\id$ of $R$.
%\end{proposition}
%\begin{proof}
%Since $R$ is excellent, it is analytically normal. Therefore, with $\Lambda$
%as in the proof of
%Corollary~\ref{c:normal}, almost all $\seq Rw$ are normal, and so is $\complet
%R_{(k,i)}$ by Lemma~\ref{l:scalarext}.  The result now follows from the next
%lemma applied to almost all approximations $\seqaff R{p(w)}\to \seq Rw$.
%\end{proof}
%
%
%
%\begin{lemma}
%Let $A\to B$ be a cyclically pure extension of domains with $B$ normal.
%Then $IA^+\cap A=IB^+\cap A$ for every ideal $I$ of $A$.
%\end{lemma}
%\begin{proof}
%Let $p$ be the characteristic of $A$ if $\op{char} A>0$, and $p=1$ otherwise.
%Suppose that $f_0,f_1,\dots,f_n\in A$ and $y_1,\dots,y_n\in B^+$ satisfy
%$f_0=f_1y_1+\cdots+f_ny_n$; we need to find $z_1,\dots,z_n\in A^+$ with
%$f_0=f_1z_1+\cdots+f_nz_n$. 
%Let $K$ be the fraction field of $B$, with perfect closure $K^{1/p^\infty}$.
%By a trace argument we see that $B^+\cap K^{1/p^\infty}\to B^+$ is pure. Hence
%we may assume that there exists $e\geq 0$ with
%$y_j^{p^e}\in B^+\cap K$ for all $j$, that is, 
%$y_j^{p^e}\in B$ for all $j$ (since
%$B$ is normal). Since $A\to B$ is cyclically pure, we find
%$w_1,\dots,w_n\in A$ such that $f_0^{p^e}=f_1^{p^e}w_1+\cdots+f_n^{p^e}w_n$;
%hence $z_j:=w_j^{p^{-e}}\in A^+$ have the required property.
%\end{proof}

\subsection{Rational singularities}

%So far all \emph{pure} applications (that is to say, applications 
%which do not directly involve the construction of $\hull R$ or $\BCM 
%R$) can be proven by different means. Therefore, 

The main merit of the present approach to tight closure in 
equicharacteristic zero and to the construction of balanced big 
Cohen-Macauly modules, via $\hull R$, is its flexibility. We 
want to finish with a brief discussion of one possible 
application of our construction of $\BCM R$, 
which we formulate in two Conjectures.

Let us return to the situation of the Hochster-Roberts Theorem, that 
is to say, a cyclically pure \homo\ from a Noetherian local ring $R$ 
into a regular local ring $S$. We already showed that $R$ (and also 
its completion $\complet R$) is \CM\ and normal (see 
Theorem~\ref{T:HR} and Remarks~\ref{R:norm} and \ref{R:HRBCM}). In 
case $R$ and $S$ are of finite type over $\mathbb C$, Boutot has 
shown in  \cite{Bou}, using deep Vanishing Theorems, that $R$ has 
rational singularities. In fact, he proves an even stronger result in 
that he only needs to assume that $S$ has rational singularities. 
Recall that an equi\ch\ zero excellent local domain $R$ has 
\emph{rational singularities} (or, more correctly, is 
\emph{pseudo-rational}) if it is normal, analytically unramified and 
\CM, and the canonical embedding
        \begin{equation*}
        H_0(W,\omega_W)\to H_0(X,\omega_X)
        \end{equation*}
is surjective (it is always injective),  where $W\to X:=\op{Spec}R$ 
is a desingularization, and where in general, $\omega_Y$ denotes the 
canonical sheaf on a scheme $Y$.

In the affine case, the methods of the second author (via 
non-standard tight closure in \cite{SchRatSing}, and via big \CM\ 
algebras in \cite{SchBCM}) yield more elementary arguments for the 
fact  that a cyclically pure subring of an affine regular ring has 
rational singularities.  Moreover, in the second paper, a more 
general version is proven, where  $S$ is only  assumed to have 
rational singularities and be Gorenstein. However, for this stronger 
version, one needs a result of Hara in \cite{HaRat}, which itself 
uses deep Vanishing Theorems. In any case, we expect that one can 
generalize Boutot's result by removing the condition that the rings 
are finitely generated over a field. (Note that no Vanishing Theorems 
are known to hold for arbitrary excellent schemes.)

\begin{conjectuur}\label{C:bout}
Every equi\ch\ zero excellent local ring $R$ which admits a cyclically pure 
\homo\ into a regular local ring $S$ is pseudo-rational.
\end{conjectuur}

In fact, we suspect that an excellent local domain is pseudo-rational 
if there exists a system of parameters $\mathbf z$ such that $\mathbf 
zR=(\mathbf zR)^+$ (= generic tight closure of $\mathbf zR$, by 
Proposition~\ref{P:plus}). It is   clear by \eqref{i:reg} how this 
implies the Conjecture.  If $R$ is in addition $\mathbb 
Q$-Gorenstein, then in the affine case   it has log-terminal 
singularities by \cite[Theorem B and Remark 3.13]{SchLogTerm}. (Here 
again we can weaken the assumption on $S$ to be only log-terminal, 
provided we use Hara's result; see that article for the terminology.) 
In view of this, we postulate the following generalization.

\begin{conjectuur}
Every equi\ch\ zero excellent local $\mathbb Q$-Gorenstein ring
which admits a cyclically pure \homo\ into a regular local ring
has log-terminal singularities.
\end{conjectuur}

The conjecture would follow from \eqref{i:reg}, if one can show that 
a $\mathbb Q$-Gorenstein excellent local ring $R$ in which each ideal is 
equal to its \pc, or, equivalently, for which $R\to\BCM R$ is 
cyclically pure, has log-terminal singularities.

\providecommand{\bysame}{\leavevmode\hbox to3em{\hrulefill}\thinspace}

\bibliographystyle{amsplain}
\bibliography{TransferAA}

%\bibitem{SchBounds}
%\bysame, \emph{},
%  Connections between Model Theory and Algebraic and Analytic Geometry
%  (A.~Macintyre, ed.), Quaderni di Mathematica, vol.~6, 2000, pp.~43--93.

\end{document}